# FUZZY LINGUISTIC
# TOPOLOGICAL SPACES

W. B. Vasantha Kandasamy
Florentin Smarandache
K.Amal

**2012**

# FUZZY LINGUISTIC TOPOLOGICAL SPACES

# CONTENTS

# PREFACE

In this book the authors introduce a new type of fuzzy topological space called the fuzzy linguistic topological space. This topological space consists of fuzzy linguistic terms together with zero. This topological space depends on the problem, attributes chosen and the metric defined on it. This is the main difference between the usual topological space and the fuzzy linguistic topological space.

Further this topological space also may be graph-connected under a particular problem and under a specific expert. Likewise the fuzzy topological space may be lattice-connected under one expert and may be disconnected under another expert depending ton the problem under investigation.

We have used the space elements and have constructed fuzzy linguistic matrices and fuzzy linguistic polynomials. These algebraic structures are used in the study of fuzzy linguistic models, defined and described in chapter four.

The striking features of these models are that the models can be used with ease by any expert, be it a socio-scientist, technologist or a medical scientist. However the fuzzy linguistic membership can be converted into fuzzy membership and vice versa.

This book has five chapters. Chapter one is introductory in nature. Fuzzy linguistic spaces are introduced in chapter two. Fuzzy linguistic vector spaces are introduced in chapter three. Chapter four introduces fuzzy linguistic models. The final chapter suggests over 100 problems and some of them are at research level.

We thank Dr. K.Kandasamy for proof reading and being extremely supportive.

W.B.VASANTHA KANDASAMY
FLORENTIN SMARANDACHE
K. AMAL

**Chapter One**

# INTRODUCTION

In this book we introduce a new notion called fuzzy linguistic space or fuzzy linguistic set L, which consists of only fuzzy linguistic terms like low, high, very low, just low, very high, just high, 0 and so on.

That is we call a 'word' as a fuzzy linguistic term if it is an adjective describing the situation in a very sensitive and a minute way. For any set X we give the fuzzy set associated with X the membership values from the set [0, 1]. Here we give for the set X the fuzzy linguistic terms as its membership values. L denotes the fuzzy linguistic space or set.

For instance if X is a set of concepts say $C_1, C_2, \ldots, C_n$, we define $f : X \to L$ (f a map from X to the fuzzy linguistic space or set) if $f(C_i) = m_i \in L$ where $C_i \in X$. Thus $C_i$ has the fuzzy linguistic membership $m_i$, if $m_i = 0$ we say $C_i$ has no fuzzy linguistic membership. So in the usual fuzzy membership, if [0, 1] is replaced by L we get the fuzzy linguistic membership. Likewise if in a fuzzy matrix $M = (m_{ij})$; $m_{ij} \in [0, 1]$ the values $m_{ij}$ is replaced by values in L then we call M, a fuzzy linguistic matrix.

$$M = \begin{bmatrix} \text{low} & 0 & \text{very high} & \text{just low} \\ \text{medium} & \text{high} & \text{just high} & 0 \\ 0 & \text{low} & \text{medium} & \text{high} \\ \text{high} & \text{just high} & 0 & \text{low} \end{bmatrix}$$

is an example of a fuzzy linguistic matrix low, high, very high, just high, medium, 0, just low $\in$ L.

Now we can define distance between two fuzzy linguistic terms say distance between low and high is far, very low and very high is very far, low and just low are very near and so on.

We venture to develop such spaces, mainly because while dealing with fuzzy models a socio scientist or a technologist not very well versed with mathematics may feel very comfortable and confident while working with fuzzy linguistic models. Thus it will be a boon to researchers not well versed with mathematics.

Further we build fuzzy linguistic models using the fuzzy linguistic sets which can be used both as mathematical models as well as linguistic models. Study in this direction is carried out in this book and its applications are given.

**Chapter Two**

# FUZZY LINGUISTIC SPACES

In this chapter the authors for the first time introduce the notion of fuzzy linguistic spaces. We have varying types of fuzzy linguistic spaces. In the first place we use the term to show that linguistic terms like more, less, little less, very less, little more, much more, little, very much etc. are "fuzzy" for these linguistic terms are qualified and such qualification is indeed what matters. When we use in fuzzy models the fuzzy values are values from the interval [0,1], so a numerical value between 0 and 1 is associated with; however naturally a non mathematician is not at home with it. To cater to the needs of a non mathematician we the authors define in this book fuzzy linguistic spaces and define operations on them.

Before we define this notion we illustrate them by some examples.

Consider L = {more, less, little less, little more, very less, much more, some what less, some what more, not that less, not that more, too less, too much, too little}. We call L a comparable fuzzy linguistic space.

For we see any two terms in L can be compared or at times they are partly equivalent. For we see 'too much' and 'much more' are partly equivalent or equivalent.

For a person can say this is 'too much' and another person may use the term 'much more' but may be the reply for the same gesture. Further the notion of fuzzy linguistic space's elements make sense only under an 'action' or an 'instance' or a 'problem' or some 'study' or some 'gesture'. As such the space L has no meaning or any relevance, except a collection of adjectives or fuzzy linguistic terms.

We see we have three types of fuzzy linguistic spaces. At the outset all fuzzy linguistic spaces are finite, the notion of unperceived infinite has no meaning in this set up. Further all fuzzy linguistic spaces are assumed to contain '0'.

We say type one if in the fuzzy linguistic space every linguistic term is comparable with another we say type two if some are comparable and some uncomparable.

We say type three if the fuzzy linguistic space is divided into subsets and each of the subsets are comparable.

We say a subset is comparable only if it has more than one term different from zero and the terms are comparable.

We further assume all linguistic spaces contain the '0' element 'zero' element. Thus we would say fuzzy linguistic space associated with a 'problem' or 'action' or 'gesture'.

For example the 'problem' is the students performance in studies. The fuzzy linguistic space associated with this problem are L = {0, good, excellent, very good, very bad, bad, average, worst, medium, little bad, better, best, above average, below average, etc}, L is a comparable set and certainly universally the performance of a students will fall in L. It is pertinent to mention that equivalent terms are also present in L. For we see excellent is the same as best, very bad and worst are not equivalent terms.  We can define linguistic operations of them, "better" of two terms or "maximum" of two terms likewise "minimum" or 'lesser' of the two terms.  So mathematically we can define 'min' operation of a fuzzy linguistic space and it is a closed binary operation on L, provided every element in L is comparable or partially comparable.  We can also say if the binary operation "better' or "maximum" is used and assume the binary operation on L is a closed binary operation.

Now we will illustrate them by some examples.

*Example 2.1:* Let L = {0, good, poor, average, bad, very bad, worst, very poor, very good} be a linguistic space associated with the performance of students in an exam.  We can define 'min' operation or 'max' operation both are defined on L.

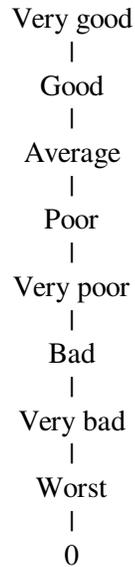

Very good
|
Good
|
Average
|
Poor
|
Very poor
|
Bad
|
Very bad
|
Worst
|
0

```
min {good, poor}     = poor
min {poor, very bad} = very bad
max {good, poor}     = good
max {poor, very bad} = poor
```

min {0, very good} = 0
max {0, very good} = very good.

We have a chain or total order and even we can call this {L, max, min} as a fuzzy linguistic chain lattice.

***Example 2.2:*** Suppose we have a fuzzy linguistic space associated with the speed of the drivers in a specific route (on the road) driving a car. Let L denote it L = {0, average, very fast, slow, very slow, fast, fairly fast, fairly slow, below average and above average}. We see every fuzzy linguistic term in L is comparable and we have both 'min' and 'max' operation defined on L. Further max {fast, average} = fast, min {average, fairly slow} = fairly slow, max {average, fairly slow} = average, min {fast, average} = average.

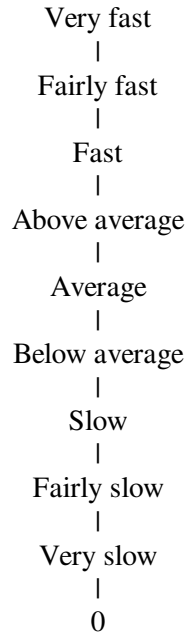

Very fast
|
Fairly fast
|
Fast
|
Above average
|
Average
|
Below average
|
Slow
|
Fairly slow
|
Very slow
|
0

We can say {L, max, min} is a fuzzy linguistic chain lattice. {L, max} and {L, min} are fuzzy linguistic semilattices.

We now give yet another example.

***Example 2.3:*** Let L be the fuzzy linguistic space associated with the "product" produced by an industry.

L = {0, bad, very bad, average, good, very good, worst, best, above average, excellent}; we can use the max or min operation on them.

min {best, excellent} = best and excellent
max {best, excellent} = best and excellent
min {very good, best} = very good
max {very good, best} = best

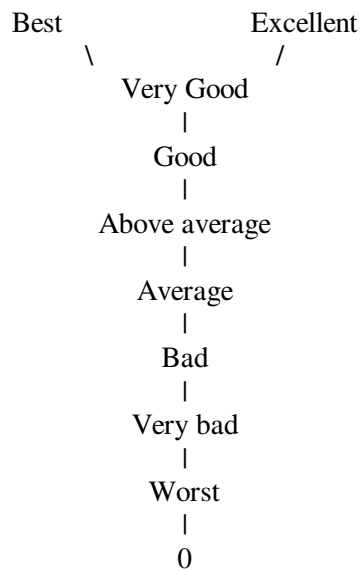

Clearly we cannot give a lattice structure or cannot have a semilattice structure under max operation. So in such situations we cannot give the fuzzy linguistic space any lattice or semilattice structure under max operation.

Thus one is always not sure to state that a fuzzy linguistic space will have a lattice structure or a semilattice structure. We can be guaranteed of this only when every pair elements of a fuzzy linguistic space are comparable.

**THEOREM 2.1:** *Let L be the fuzzy linguistic space. If every pair of elements is comparable, then L is a fuzzy linguistic chain lattice.*

Proof is direct and hence is left as an exercise to the reader.

**COROLLARY 2.1:** If L is a fuzzy linguistic space, L is not a fuzzy linguistic lattice if L has more than one greatest element.

**COROLLARY 2.2:** If L is a fuzzy linguistic space, L is a fuzzy linguistic lattice if L is a partially ordered set and with one and only one greatest element for '0' is always assumed to be the least element (provided negative fuzzy linguistic terms are not used in L).

*Note:* If two elements in L do not have a minimum than we assume '0' to be the minimum. Similarly if two elements in L do not have a maximum we assume the greatest element if it exists in L as a maximum otherwise; the operation max. is not defined on L (only under the assumption negative small or negative large, that is negative fuzzy linguistic terms are not used).

Now we say a fuzzy linguistic space L is a fuzzy linguistic lattice if for every pair of elements in L have the min and max in L (L is assumed to be a comparable fuzzy linguistic space).

Suppose L be a fuzzy linguistic space, we can define metric on L as follows which will be known as "fuzzy linguistic metric" or "fuzzy distance" and denoted by $F_m^l$ or $F_l^m$. Let L be a fuzzy linguistic space; for consider an associated fuzzy linguistic ordered space M of L where M = {0, near, very near, over lap, more over lap, partly near, little far, just near, very far, far, partly far, not comparable} in which the fuzzy linguistic terms are the ones which shows comparison between them; M can also have any other equivalent linguistic terms which compares the 'fuzzy distance' that is the shade of difference between two terms, for instance average and just average are two fuzzy linguistic terms the fuzzy distance $F_m^l$ (average, just average) = over lap, that is the fuzzy distance between the fuzzy linguistic terms, average and just average is just a over lap over each other where as the distance between the two fuzzy linguistic terms "above average" and "average" is over lap with each other in fuzzy linguistic metric. That is we can compare any two comparable fuzzy linguistic terms and find the fuzzy linguistic distance between them. Likewise 'poor' and 'very good' be two fuzzy linguistic terms, the fuzzy linguistic distance between them is 'far'. Suppose we have say terms like "beautiful" and "very far" we see no sense can be made for we can compare "beautiful" with "just beautiful" and say 'near' each other as fuzzy linguistic terms but if "very far" is a distance concept we cannot do any comparison. However "very far from being beautiful" can be compared with beautiful and say $F_m^l$ (beautiful, very far from being beautiful) = very far and so on.

So we say a fuzzy linguistic space L endowed with a fuzzy linguistic metric $F_l^m$ or $F_m^l$ as a "fuzzy linguistic topological space" with the fuzzy linguistic metric $F_l^m$ ($F_m^l$); i.e., (L, $F_l^m$) (L, $F_m^l$) is a fuzzy linguistic topological space.

We give first a few examples of them.

*Example 2.4:* Let L be the fuzzy linguistic space associated with temperature obtained from an experiment in a laboratory. L = {negative large, negative medium, negative small, zero, positive small, positive medium, positive large}. Let M = {0, larger, small, large, overlap, just small, very small and very large} be the distance in fuzzy linguistic terms associated with L. To find, $F_l^m$ (Negative large, negative medium} = over lap.

$F_l^m$ (negative large, negative small) = very small

$F_l^m$ (negative medium, negative small) = over lap.

$F_1^m$ (negative large, 0) = small

$F_1^m$ (negative medium, 0) = just small

$F_1^m$ (negative small, 0) = very small

$F_1^m$ (negative large, positive small) = just large

$F_1^m$ (negative large, positive medium) = large

$F_1^m$ (negative large, positive large) = largest.

$F_1^m$ (negative medium, positive small) = large

$F_1^m$ (negative medium, positive medium) = larger

$F_1^m$ (negative medium, positive large) = very large

$F_1^m$ (negative small, positive small) = small

$F_1^m$ (negative small, positive large) = just large

$F_1^m$ (negative small, positive medium) = large

$F_1^m$ (positive small, 0) = very small

$F_1^m$ (positive medium, 0) = just small

$F_1^m$ (positive large, 0) = small

$F_1^m$ (positive small, positive medium) = over lap

$F_1^m$ (positive small, positive large) = small

$F_1^m$ (positive medium, positive large) = over lap

$F_1^m$ (positive small, positive small) = 0.

It is pertinent to mention here that in this example we have a fuzzy linguistic metric which declares over lap.

We may have a collection of linguistic terms which have no over lap associated with a problem. In such case the fuzzy linguistic distance will not include the term over lap. Further it is in the hands of the investigator to choose the metric and define the metric. The flexibility is useful.

We will illustrate this by a simple example.

*Example 2.5:*  Let L be a fuzzy linguistic space associated with the fuzzy linguistic terms:  L = {good, bad, 0, average, best, worst} and M = {very near, near, far, very far} be the fuzzy linguistic metric.

Now we can find $F_1^m$ (good, bad) = far,

$F_1^m$ (best, worst) = very far and so on.

Thus in this space the 'over lap' does not occur. (L, $F_1^m$) is a fuzzy linguistic topological space.

If the topological fuzzy linguistic space has no over lap between two fuzzy terms we call it as discrete topological fuzzy linguistic topological space and if the 'overlap' occurs in the fuzzy linguistic topological space we call the space a non discrete topological space. We further wish to state that all our fuzzy linguistic topological spaces also vary with time for concepts are associated with the fuzzy linguistic terms and the concepts vary with time or its influence over the other concept.

We say also a fuzzy linguistic topological space which is not discrete as a fuzzy linguistic overlapping topological space. Further we define a fuzzy linguistic topological space to be connected straight or straight connected or chain lattice connected if the fuzzy linguistic terms of that fuzzy linguistic topological space forms a fuzzy linguistic chain lattice. We say a fuzzy linguistic topological space is lattice connected if the fuzzy linguistic terms forms a fuzzy linguistic lattice.

We have the following theorem the proofs of which are direct and hence left as an exercise to the reader.

**THEOREM 2.2:**  *Let (L, M) be a fuzzy linguistic topological discrete space.  Every subspace of (L, M) is a discrete fuzzy linguistic topological subspace.*

**THEOREM 2.3:**  *Let (L, M) be a fuzzy linguistic topological space which is over lapping.  Every subspace of (L, M) need not*

*be a over lapping fuzzy linguistic topological subspace. (L, M) can also have discrete fuzzy linguistic topological subspaces.*

The proof of these two theorems are straight forward and hence left as an exercise to the reader.

**THEOREM 2.4:** *Every straight fuzzy linguistic (chain) lattice connected fuzzy linguistic topological space is a fuzzy linguistic lattice connected fuzzy linguistic topological space and not conversely.*

This proof is also simple for an interested reader can prove and also construct counter examples.

Now we just recall that for we have a fuzzy linguistic space L associated with some concepts / attributes. These concepts or attributes take all the fuzzy linguistic values / terms in that space depending on the other or related concepts / attributes in that problem. When we draw the graph using these concepts as vertices / nodes and edges as the fuzzy linguistic terms we get the graphs which we define as a fuzzy linguistic graph associated with the problem / the fuzzy linguistic space L.

For instance L is a fuzzy linguistic space associated with a problem and the concepts / attributes related with the problem be $C_1, C_2, \ldots, C_{10}$. The fuzzy linguistic nodes be $C_1, C_2, \ldots, C_{10}$ and the fuzzy linguistic graph is given in the following figure.

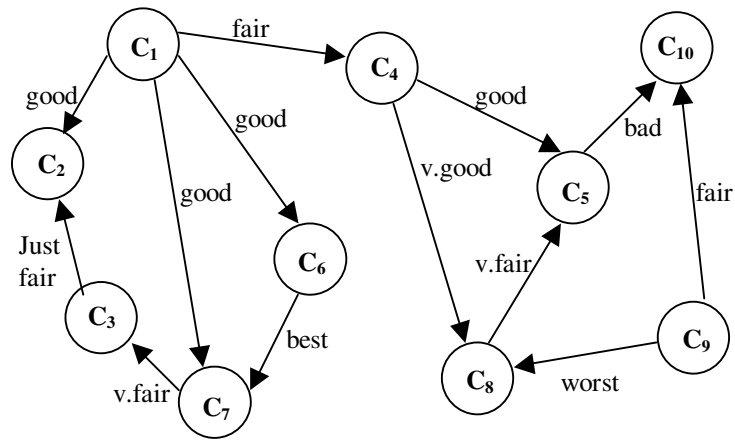

This sort of graphs where $C_1, C_2, \ldots, C_{10}$ are fuzzy linguistic nodes and fuzzy linguistic edges are fair, good, bad, v-good, v-fair, worst and so on. This sort of graph will be known as the fuzzy linguistic graph associated with the fuzzy linguistic space L. The fuzzy linguistic matrix associated with the fuzzy linguistic graph is as follows:

|        | $C_1$ | $C_2$ | $C_3$ | $C_4$ | $C_5$ | $C_6$ | $C_7$ | $C_8$ | $C_9$ | $C_{10}$ |
|--------|-------|-------|-------|-------|-------|-------|-------|-------|-------|----------|
| $C_1$  | 0 | good | 0 | fair | 0 | good | good | 0 | 0 | 0 |
| $C_2$  | 0 | 0 | 0 | 0 | 0 | 0 | 0 | 0 | 0 | 0 |
| $C_3$  | 0 | just fair | 0 | 0 | 0 | 0 | 0 | 0 | 0 | 0 |
| $C_4$  | 0 | 0 | 0 | 0 | good | 0 | 0 | v.good | 0 | 0 |
| $C_5$  | 0 | 0 | 0 | 0 | 0 | 0 | 0 | 0 | 0 | bad |
| $C_6$  | 0 | 0 | 0 | 0 | 0 | 0 | best | 0 | 0 | 0 |
| $C_7$  | 0 | 0 | v.fair | 0 | 0 | 0 | 0 | 0 | 0 | 0 |
| $C_8$  | 0 | 0 | 0 | 0 | v.fair | 0 | 0 | 0 | 0 | 0 |
| $C_9$  | 0 | 0 | 0 | 0 | 0 | 0 | 0 | 0 | 0 | fair |
| $C_{10}$ | 0 | 0 | 0 | 0 | 0 | 0 | 0 | 0 | 0 | 0 |

Such type of fuzzy linguistic matrices will be studied in the following chapters. Further vector spaces using fuzzy linguistic matrices are defined and their properties derived.

Now we proceed onto use these matrices in several models like Fuzzy Cognitive Linguistic Models etc.

Now using these graphs (fuzzy linguistic graphs) we can define graph connected fuzzy linguistic topological spaces and graph disconnected fuzzy linguistic topological spaces. We say a fuzzy linguistic topological space to be graph connected if the graph connected with the fuzzy linguistic space taking concepts as vertices and fuzzy linguistic terms as edges is connected.

If the fuzzy linguistic graph is disconnected or not connected we say the fuzzy linguistic topological space is not graph connected or the fuzzy linguistic topological space is graph disconnected. So we in case of fuzzy linguistic topological spaces two types of connectedness can be defined.

Further our graph connectedness of a fuzzy linguistic topological spaces for a particular problem will vary with the experts opinion. Thus unlike the usual topological spaces we cannot say about the graph connectedness of the fuzzy linguistic topological spaces.

We will illustrate this situation by an example.

Suppose $C_1$, $C_2$, $C_3$, $C_4$, $C_5$, $C_6$ and $C_7$ are the seven attributes associated with the problem, we will the fuzzy linguistic spaces L associated with $C_1$, $C_2$, …, $C_7$ can give different fuzzy linguistic graphs depending on experts.

For instance we can have the following two fuzzy linguistic graphs.

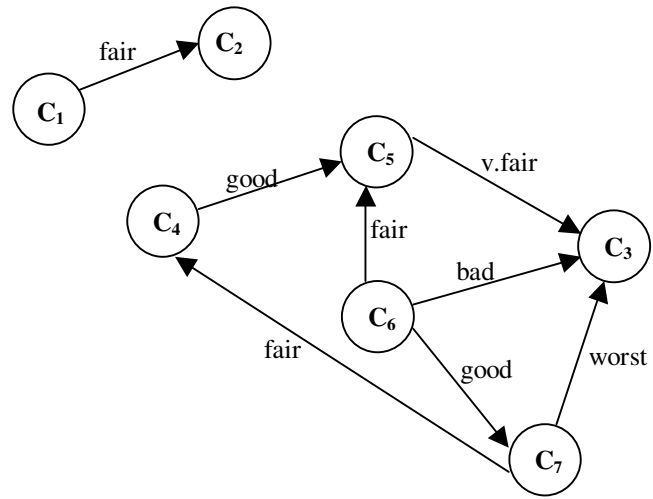

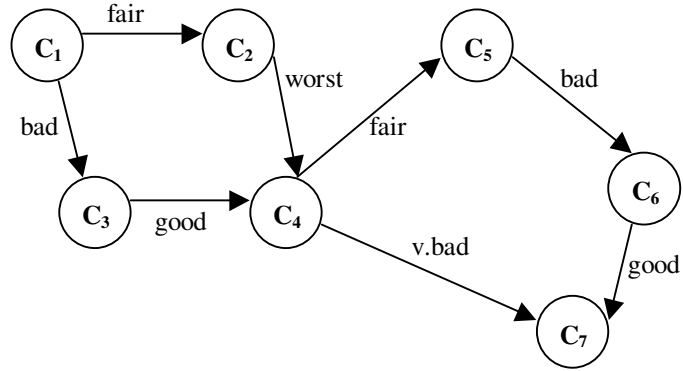

We see the first fuzzy linguistic graph is disjoint or disconnected where as the second fuzzy linguistic graph is connected.  Thus the fuzzy linguistic topological space associated with the first graph is a disconnected fuzzy linguistic topological space graph but the latter fuzzy linguistic topological space is a graph connected fuzzy linguistic topological space.

Thus we see we can have only some of the properties associated with the topological space to be true in case of fuzzy linguistic topological spaces as the fuzzy linguistic topological spaces are finite and their properties are time variant and expert variant.  Thus with this in mind we develop the fuzzy linguistic models associated with them.

Further these fuzzy linguistic spaces L are also associated basically with a problem and this space L and the fuzzy linguistic terms enjoyed by each of these concepts / attributes and also its relation with other concepts. Thus unlike usual topological spaces these fuzzy linguistic topologies are highly selective and much dependent on the problem under investigation, time, interdependent on the attributes, choice of fuzzy linguistic terms and above all the need or the interpretation of the expert.

Let P be a social problem at hand. Let L be the fuzzy linguistic terms associated with the concepts/attributes of the problem P.

Let $\{G_i\}_{i=1}^n$ be the collection of all fuzzy linguistic graphs given by n experts. If L is a graph connected fuzzy linguistic topological space; we say L is a n-strong fuzzy linguistic topological space, if for every graph the fuzzy linguistic topological space is graph connected, we say L is a super strong fuzzy linguistic topological graph connected space. If for this problem P, none of the fuzzy linguistic graphs associated with P is connected we call the fuzzy linguistic topological space is said to be graphically totally disconnected. In chapter four we will illustrate this situation by some examples.

**Chapter Three**

# FUZZY LINGUISTIC SET VECTOR SPACES

In this chapter we for the first time define the notion of set fuzzy linguistic spaces or fuzzy linguistic set vector space. We assume L denotes a fuzzy linguistic finite space. L is also a fuzzy linguistic lattice. Both our assumptions are permissible as we are going to talk about a fuzzy linguistic space associated only with a problem. Once the space is associated with a problem certainly we see all the linguistic terms involved with that problem is comparable. Hence a fuzzy linguistic space is also a fuzzy linguistic lattice. All problems are mainly social problems or problems which can be described using fuzzy linguistic terms.

**DEFINITION 3.1:** *Let L be a fuzzy linguistic space associated with a problem P (P may be social problem or a technological problem or a mathematical problem where the terms can be described using fuzzy linguistic terms). Let 'min' be a operation defined on L. L is a semigroup under 'min'. Also $0 \in L$ by the very definition of fuzzy linguistic space.*

*If L is assumed to be lattice L under the max operation to be a semigroup then {L, min, max} is a fuzzy linguistic lattice. {L, min} is a fuzzy linguistic semigroup and {L, max} is a fuzzy linguistic semigroup..*

Now at times we may not be in a position to give a lattice structure to the fuzzy linguistic spaces.

We give examples of fuzzy linguistic spaces.

**Example 3.1:** Suppose one is interested in the experiment of studying the temperature variation in a chemical plant. The fuzzy linguistic space L associated with the problem is {0, negative large, negative small, negative medium, positive large, positive small, positive medium and largest}.

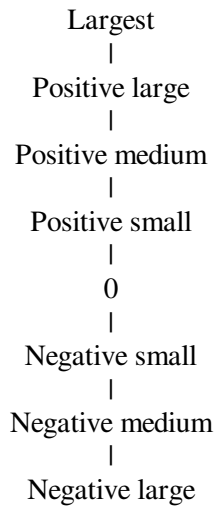

Largest
|
Positive large
|
Positive medium
|
Positive small
|
0
|
Negative small
|
Negative medium
|
Negative large

L is a chain lattice with negative large as the with least element.

**Example 3.2:** Suppose the study of the problem that causes cancer is carried out by an expert. If the concepts are taken as smoking tobacco, chewing tobacco, environment is present with tobacco smoke, alcohol, poverty, use of pesticides are taken as some of the concepts. The fuzzy linguistic terms are taken as

{much, very much, more often, often, very often, most often, at times, some times, frequently, more frequently, less often, less frequently, 0, never, occasionally, always, most of the time} = L. We see 0 is the least element. Most of the time is the greatest element of L and L is a linguistically partially ordered set for it; '≤' denotes the partial order.

often ≤ very often
often ≤ most often

however very often and most often cannot be linguistically ordered; never ≤ frequently and so on.

Thus L is a lattice but L is not a chain lattice.

However we do not get fuzzy linguistic spaces which are not partially ordered whenever the fuzzy linguistic terms are associated with a problem. Always '0' acts as the least element provided we do not involve in L negative terms. If negative terms are involved largest negative term would be the least element.

Now it is pertinent to mention that for any social problem under study the important factor to be remembered is the fuzzy linguistic terms associated with the problem may involve many concepts / attributes and with each of these attributes we have a set of linguistic terms associated with it.

Thus for us the terms can occur as a linguistic array. It may be a column array or a rectangular array or a square array or a row array. What ever be it we have to define some operations on them to arrive at a conclusion / solution. To this end we first define the notion of fuzzy linguistic matrices and define operations on them.

**DEFINITION 3.2:** *Let L be a fuzzy linguistic space. Let $B = \{(x_1, x_2, …, x_n)$ where $x_i \in L; 1 \leq i \leq n\}$. We define B to be fuzzy linguistic row matrix.*

We will give examples of them.

*Example 3.3:* Let X = (very fair, fair, fair, just fair, 0) be a fuzzy linguistic row matrix.

*Example 3.4:* Let X = (low, very low, low, 0, just low, high, 0, very low) be a fuzzy linguistic row matrix.

*Example 3.5:* Let M = (0, fast, not that fast, super fast, 0, fast, not fast) be a fuzzy linguistic row matrix.

Now we show we can define operations on fuzzy linguistic row matrices of same order provided the entries of each of these matrices is from the same fuzzy linguistic row matrices of same order and provided the entries of each of these matrices is from the same fuzzy linguistic space L.

**DEFINITION 3.3:** *Let L be a fuzzy linguistic space,*
$L_R = \{X = (a_1, a_2, \ldots, a_n) \mid a_i \in L, 1 \leq i \leq n\}$; *$L_R$ under min (or max) operation is a semigroup called the semigroup of row matrix fuzzy linguistic space or the associated semigroup of the fuzzy linguistic space L.*

We will just illustrate them by some examples.

*Example 3.6:* Let L be a fuzzy linguistic space.
$M = \{(x_1, x_2, x_3, x_4) \mid x_i \in L, 1 \leq i \leq 4\}$ be a fuzzy linguistic space semigroup associated with M under max operation.

Let X = (strong, 0, weak, very weak) and Y = (fairly strong, weak, very weak, strong) be in M.

Max {X, Y} = {(strong, 0, weak, very weak), (fairly strong, weak, very weak, strong)} = (strong, weak, weak, strong). This is the way max operation is performed on M. Thus {M, max} is a fuzzy linguistic semigroup of row fuzzy linguistic matrices.

*Example 3.7:* Let L be a fuzzy linguistic space.
$T = \{(a_1, a_2, a_3, a_4, a_5, a_6) \mid a_i \in L; 1 \leq i \leq 6\}$; T under min is a fuzzy linguistic row matrix semigroup.

Consider X = (low temperature, 0, negative high, negative medium, 0, high) and Y = (just low, positive medium, 0, medium, low, high) be two fuzzy linguistic row matrices in T. min {X, Y} = {(low temperature, 0, negative high, negative medium, 0, high), (just low, positive medium, 0, medium, low, high)} = (low temperature, 0, 0, negative medium, 0, high). Thus {T, min} is a fuzzy linguistic row matrix semigroup.

**DEFINITION 3.4:** *Let L be a fuzzy linguistic space.*

*M = {($a_1$, $a_2$, …, $a_n$) | $a_i \in L$, $1 \leq i \leq n$} be a fuzzy linguistic row vector semigroup under min operation (or max operation, or used in the mutually exclusive sense). Let V $\subseteq$ M; be a proper subset of M, if V under the min operation (or max operation) is a semigroup then we call V to be a subsemigroup of the fuzzy linguistic row matrices of M.*

We give some examples.

*Example 3.8:* Let L be a fuzzy linguistic space.

M = {($x_1$, $x_2$, $x_3$, $x_4$, $x_5$, $x_6$, $x_7$, $x_8$, $x_9$) | $x_i \in L$; $1 \leq i \leq 9$} be a fuzzy linguistic row matrix semigroup under min.

P = {($y_1$, $y_2$, $y_3$, $y_4$, 0, 0, 0, 0, 0) | $y_i \in L$; $1 \leq i \leq 4$} $\subseteq$ M, P is a subsemigroup of fuzzy linguistic row matrices.

W = {(0, $x_1$, 0, $x_2$, 0, $x_3$, 0, $x_4$, 0) | $x_i \in L$; $1 \leq i \leq 4$} $\subseteq$ M is a subsemigroup of fuzzy linguistic row matrices under min operation. Infact both P and W are also ideals of M under min operation defined as the fuzzy linguistic row matrix ideal of M.

*Example 3.9:* Let L be a fuzzy linguistic space. L = {0, low, very low, medium low, high, medium, high, very high, medium high, just low, just high} be the cost estimation in a company. Suppose max is the binary operation defined on the fuzzy linguistic matrix M = {($x_1$, $x_2$, $x_3$) | $x_i \in L$, $1 \leq i \leq 3$}.

Suppose W = {($y_1$, $y_2$, $y_3$) | $y_i \in$ {high, medium high, very high, medium}, $1 \leq i \leq 3$} $\subseteq$ M, W is a fuzzy linguistic row matrix subsemigroup of M under max operation. Further W is also an ideal of fuzzy linguistic row matrix semigroup M. However if S = {($x_1$, $x_2$, $x_3$) | $x_i \in$ {0, low, very low, just low,

medium low} ⊆ M, S is only a subsemigroup of M and not an ideal of M.

We now proceed onto define zero divisors in these semigroups. We can define zero divisors only in case of fuzzy linguistic row matrix semigroups under min operation on them.

Suppose L be a fuzzy linguistic space,
$S = \{(x_1, x_2, \ldots, x_n) \mid x_i \in L\}$ be a fuzzy linguistic row matrix semigroup under min operation. S has zero divisors. For take n = 5, then if $X = (x_1, 0, 0, x_2, 0)$ and $Y = (0, y_1, y_2, 0, y_3)$ are in S then we see min $\{X, Y\} = \{(0, 0, 0, 0, 0)\}$.

Thus we have zero divisors. If the fuzzy linguistic row matrix semigroup S has 'max' operation on it we see S cannot have zero divisors.

Next we proceed onto define fuzzy linguistic column matrix and the fuzzy linguistic column matrix semigroup under 'min' or 'max' operation.

Let L be a fuzzy linguistic space.

Suppose

$$X = \begin{bmatrix} m_1 \\ m_2 \\ \vdots \\ m_n \end{bmatrix} \text{ where } m_i \in L, 1 \leq i \leq n,$$

then we define L to be fuzzy linguistic column matrix.

We can define max (or min) of two fuzzy linguistic column matrices of same order which is as follows.

If $X = \begin{bmatrix} m_1 \\ m_2 \\ \vdots \\ m_n \end{bmatrix}$ and $Y = \begin{bmatrix} s_1 \\ s_2 \\ \vdots \\ s_n \end{bmatrix}$ then

$$\max(X, Y) = \begin{bmatrix} \max(m_1, s_1) \\ \max(m_2, s_2) \\ \vdots \\ \max(m_n, s_n) \end{bmatrix} \text{ (or}$$

$$\min(X, Y) = \begin{bmatrix} \min(m_1, s_1) \\ \min(m_2, s_2) \\ \vdots \\ \min(m_n, s_n) \end{bmatrix}).$$

We see the concept of min or max can be defined only if both the column matrices X and Y are of same order otherwise we cannot define any such notion.

Let L be a fuzzy linguistic space.

$$S = \left\{ \begin{bmatrix} y_1 \\ y_2 \\ \vdots \\ y_n \end{bmatrix} \middle| y_i \in L, 1 \leq i \leq n \right\}$$

be the collection of $n \times 1$ fuzzy linguistic column vectors or column matrices. On S we can define the binary operation min (or max 'or' used in the mutually exclusive sense).

Thus (S, min) (or (S, max)) is a semigroup called the fuzzy linguistic column matrix semigroup.

We see all properties described and defined in case of fuzzy linguistic row vectors can be easily extended in case of fuzzy linguistic column vectors.

We will describe this by some examples.

*Example 3.10:* Let L be a fuzzy linguistic space associated with the performance of school students in mathematics exam. L = {0, good, very good, very poor, poor, average, above average, below average, fair, very fair}.

Let

$$M = \left\{ \begin{bmatrix} x_1 \\ x_2 \\ x_3 \\ x_4 \\ x_5 \end{bmatrix} \,\middle|\, x_i \in L, 1 \leq i \leq 5 \right\},$$

M under 'min' operation is a semigroup.

$$\text{Take } X = \begin{bmatrix} 0 \\ \text{good} \\ 0 \\ \text{poor} \\ \text{average} \end{bmatrix} \text{ and } Y = \begin{bmatrix} \text{very good} \\ \text{good} \\ 0 \\ \text{very poor} \\ \text{below average} \end{bmatrix} \text{ in M;}$$

we have

$$\min\{X, Y\} = \begin{bmatrix} 0 \\ \text{good} \\ 0 \\ \text{very poor} \\ \text{below average} \end{bmatrix}.$$

Thus M is closed with respect to the binary operation min. Suppose the expert or the teacher wants to try with max operation on M, then

$$\max\{X, Y\} = \begin{bmatrix} \text{very good} \\ \text{good} \\ 0 \\ \text{poor} \\ \text{average} \end{bmatrix},$$

thus one can see the best cut off by using max operation and the possible correct situation by defining min operation on M.

*Example 3.11:* Let L be a fuzzy linguistic space associated with the performance of the workers in a firm. Let

$$P = \left\{ \begin{bmatrix} x_1 \\ x_2 \\ x_3 \\ x_4 \\ x_5 \\ x_6 \\ x_7 \end{bmatrix} \,\middle|\, x_i \in L \right\}$$

= {poor, 0, very poor, better, best, below average, good, very good, fair, average, above average, $1 \leq i \leq 7$} belonging to seven units / divisions of the firm; say, manufacturing, packaging, testing, marketing, purchasing, selling and transporting. The experts gives only these eleven fuzzy linguistic terms. If operation max is defined then P with max we will get the best out come of the seven units of the firm. If min is the operation taken on P then this will give the least or the worst possible performance of the firm in each of these seven units. They can also compare by assigning both max and min operations on P.

*Example 3.12:* Let L be the fuzzy linguistic space associated with six wards of hospital viz; ICU, general ward, children ward, women ward, cancer ward and TB ward. The fuzzy linguistic terms associated with each of these wards are L = {sick, very sick, recovering, slowly recovering, fastly recovering, serious, no recovery, deterioting, very serious, some what serious}.

$$\text{Let } M = \left\{ \begin{bmatrix} w_1 \\ w_2 \\ w_3 \\ w_4 \\ w_5 \\ w_6 \end{bmatrix} \right\} w_i \in L, 1 \leq i \leq 6,$$

that is patients status in the ward $w_i$}. We can use min or max depending on the expert to get the best situation or the worst scenario.

Thus we can have fuzzy linguistic column matrix semigroups.

Now we proceed onto describe the concept of fuzzy linguistic rectangular matrix semigroup under max or min.

Let L be a fuzzy linguistic space.

Suppose

$$M = \begin{bmatrix} m_{11} & m_{12} & ... & m_{1n} \\ m_{21} & m_{22} & ... & m_{21} \\ \vdots & \vdots & & \vdots \\ m_{m1} & m_{m2} & ... & m_{mn} \end{bmatrix} \text{ where } m_{ij} \in L;$$

$1 \leq i \leq m$ and $1 \leq j \leq n$.

We call M a fuzzy linguistic rectangular matrix.

We will illustrate them with examples.

*Example 3.13:* Let L be a fuzzy linguistic space.

$$P = \left\{ \begin{bmatrix} a_1 & a_2 & a_3 & a_4 & a_5 \\ a_6 & a_7 & a_8 & a_9 & a_{10} \\ a_{11} & a_{12} & a_{13} & a_{14} & a_{15} \\ a_{16} & a_{17} & a_{18} & a_{19} & a_{20} \end{bmatrix} \middle| a_i \in L, 1 \le i \le 20 \right\}$$

denotes the collection of all $4 \times 5$ fuzzy linguistic matrices.

We can define either min or max operation on P.

With min operation on P we see P is a fuzzy linguistic rectangular semigroup and also under the max operation P is a fuzzy linguistic rectangular semigroup.

*Example 3.14:* Let L be a fuzzy linguistic space.

$$M = \left\{ \begin{bmatrix} a_1 & a_2 & a_3 \\ a_4 & a_5 & a_6 \\ a_7 & a_8 & a_9 \\ a_{10} & a_{11} & a_{12} \end{bmatrix} \middle| a_i \in L, 1 \le i \le 12 \right\}$$

be the fuzzy linguistic semigroup under 'max'.

*Example 3.15:* Let L be a fuzzy linguistic space.

$$T = \left\{ \begin{bmatrix} a_1 & a_2 & a_3 & a_4 & a_5 & a_6 \\ a_7 & a_8 & a_9 & a_{10} & a_{11} & a_{12} \\ a_{13} & a_{14} & a_{15} & a_{16} & a_{17} & a_{18} \end{bmatrix} \middle| a_i \in L, 1 \le i \le 18 \right\}$$

be a fuzzy linguistic semigroup of rectangular matrices with min operation.

***Example 3.16:*** Let L be a fuzzy linguistic space.

$$P = \left\{ \begin{bmatrix} a_1 & a_2 & a_3 & a_4 \\ a_5 & a_6 & a_7 & a_8 \\ a_9 & a_{10} & a_{11} & a_{12} \\ a_{13} & a_{14} & a_{15} & a_{16} \\ a_{17} & a_{18} & a_{19} & a_{20} \\ a_{21} & a_{22} & a_{23} & a_{24} \\ a_{25} & a_{26} & a_{27} & a_{28} \\ a_{29} & a_{30} & a_{31} & a_{32} \end{bmatrix} \middle| \; a_i \in L, \; 1 \le i \le 32 \right\}$$

be a fuzzy linguistic rectangular semigroup.

All these matrices are determined using the inter relation between two attributes this will be illustrated in chapter IV.

We can define fuzzy linguistic square semigroup with entries from the fuzzy linguistic space, which is a matter of routine.

We proceed onto give examples of them.

***Example 3.17:*** Let L be a fuzzy linguistic space and

$$M = \left\{ \begin{bmatrix} a_1 & a_2 & a_3 \\ a_4 & a_5 & a_6 \\ a_7 & a_8 & a_9 \end{bmatrix} \middle| \; a_i \in L, \; 1 \le i \le 9 \right\}$$

be a fuzzy linguistic square matrix; M is a semigroup both under min (or max).

*Example 3.18:* Let L be a fuzzy linguistic space.

$$P = \left\{ \begin{bmatrix} a_1 & a_2 & a_3 & a_4 \\ a_5 & a_6 & a_7 & a_8 \\ a_9 & a_{10} & a_{11} & a_{12} \\ a_{13} & a_{14} & a_{15} & a_{16} \end{bmatrix} \middle| a_i \in L, 1 \leq i \leq 16 \right\}$$

be a fuzzy linguistic square matrix semigroup under min.

*Example 3.19:* Let L be a fuzzy linguistic space.

$$M = \left\{ \begin{bmatrix} a_1 & a_2 \\ a_3 & a_4 \end{bmatrix} \middle| a_i \in L, 1 \leq i \leq 4 \right\}$$

be a fuzzy linguistic square matrix semigroup under max.

Now we show how these are used in the rectangular and square matrices.

Let L be a fuzzy linguistic space. Let $C_1, C_2, \ldots, C_n$ be some attributes and L be the fuzzy linguistic space L associated with each of the attributes.

Now if we consider the square $n \times n$ matrices which depicts the fuzzy linguistic matrix then the diagonal terms of the square $n \times n$ linguistic matrix is zero.

Likewise if L is a fuzzy linguistic space and $C_1, \ldots, C_n$ and $R_1, R_2, \ldots, R_m$ ($m \neq n$) are attributes / concepts such that both $R_i$ and $C_j$ take from time to time the linguistic terms from the fuzzy linguistic space L.

These will be discussed in the chapter on applications.

Now we have seen examples of fuzzy linguistic spaces and fuzzy linguistic matrices associated with them and discussed the algebraic structure on them.

Next we proceed onto define set vector spaces and semigroup vector space of fuzzy linguistic terms.

**DEFINITION 3.5:** *Let S be a set of fuzzy linguistic terms. V another set of fuzzy linguistic terms. We say V is a set fuzzy linguistic vector space over the set S if for all v $\in$ V and for all s $\in$ S, min (s, v) belongs to V.*

We give examples of them.

*Example 3.20:* Let
S = {good, 0, very good, bad, very bad, fair, very fair} be a linguistic set. V = {0, good, very good, fairly good, bad, fair, very fair, very bad}; we see V is a set fuzzy linguistic vector space over S.

*Example 3.21:* Let V = {low temperature, medium temperature, very low temperature, very high temperature, high temperature, zero temperature} be a fuzzy linguistic set.

Take S = {0, low temperature, high temperature} another fuzzy linguistic set. We see V is a set fuzzy linguistic vector space over S under min.

However S is not a set fuzzy linguistic vector space over the set V under min.

We say in all these cases, V is a set fuzzy linguistic vector space over the fuzzy linguistic set S with min operation.

We have a theorem which give a one way condition for a fuzzy linguistic set to be a set fuzzy linguistic vector space over another fuzzy linguistic set.

**THEOREM 3.1:** *Let V be a fuzzy linguistic set and S be another fuzzy linguistic set. If S is a proper subset of V then V is a set fuzzy linguistic vector space over S under the appropriate operation min or max depending both on S and V.*

Proof is straight forward hence left as an exercise.

However S in general is not set fuzzy linguistic vector space over V.

We now proceed onto define set fuzzy linguistic vector subspace of a set fuzzy linguistic space over a set S.

**DEFINITION 3.6:** *Let V be a set fuzzy linguistic vector space over the fuzzy linguistic set S. Let W ⊆ V; W be a proper subset of V, if W itself is a set fuzzy linguistic vector space over the set S, then we call W to be a set fuzzy linguistic vector subspace of V over the fuzzy linguistic set S.*

We will illustrate this situation by some examples.

*Example 3.22:* Let V be a fuzzy linguistic set associated with a teacher's teaching say V = {good, not that good, bad, very bad, very good, some what bad, average, fairly good, high, very low, nil, to some extent, above average, below average} and

S = {good, bad, below average} be a fuzzy linguistic set. V is a set fuzzy linguistic vector space over S.

Take
W = {not that good, bad, very bad, below average, fairly good} ⊆ V, W is a set fuzzy linguistic vector subspace of V over the fuzzy linguistic set S.

Consider
W = {high, very low, good, to some extent, nil} ⊆ V, W is a set fuzzy linguistic vector subspace of V over the set S under min.

We see as in case of usual vector spaces it may not be always possible to write the set fuzzy linguistic vector spaces as a direct sum of set fuzzy linguistic vector subspaces of V over the set S.

However it is always possible to write a set fuzzy linguistic vector space over a set as a pseudo direct sum of set fuzzy linguistic vector subspaces over the set S.

We just indicate we can as in case of set vector spaces over a set S define the notion of subset fuzzy linguistic vector subspace of a set fuzzy linguistic vector space over a subset.

The task of defining it is left as an exercise to the reader.

We illustrate this by some examples.

*Example 3.23:* Let V = {good, bad, very good, very bad, fair, very fair, average, below average, above average, best, worst} be attributes regarding the students performance is a school. V is a set fuzzy linguistic vector space over the set S = {bad, fair, good}. Consider $W_1$ = {good, bad, very good, very bad} $\subseteq$ V, $W_2$ = {very bad, fair, very fair} $\subseteq$ V, $W_3$ = {very fair, average, very bad, above average} $\subseteq$ V, $W_4$ = {very bad, above average, worst} $\subseteq$ V and $W_5$ = {very bad, worst, best fair} $\subseteq$ V be set fuzzy linguistic vector subspaces of V over the fuzzy linguistic set S. Clearly $W_i \cap W_j \neq \phi$ if $i \neq j$, $1 \leq i, j \leq 5$.

Further V $\subseteq$ $W_1$ + $W_2$ + $W_3$ + $W_4$ + $W_5$, thus V is a pseudo direct sum of set fuzzy linguistic vector subspaces of V over the set S.

Now we give examples of set fuzzy linguistic matrix vector spaces over the fuzzy linguistic set S.

*Example 3.24:* Let L be a fuzzy linguistic space associated with some problem.

$$V = \left\{ (a_1, a_2, a_3, a_4), \begin{bmatrix} a_1 \\ a_2 \\ a_3 \\ a_4 \\ a_5 \\ a_6 \end{bmatrix}, \begin{bmatrix} a_1 & a_2 & a_3 & \ldots & a_{10} \\ a_{11} & a_{12} & a_{13} & \ldots & a_{20} \end{bmatrix} \right.$$

$$\left. a_i \in L,\ 1 \leq i \leq 20 \right\}$$

be a set fuzzy linguistic vector space over the set L.

Consider

$$W_1 = \{(a_1, a_2, a_3, a_4) \mid a_i \in L,\ 1 \leq i \leq 4\} \subseteq V,$$

$$W_2 = \left\{ \begin{bmatrix} a_1 \\ a_2 \\ a_3 \\ a_4 \\ a_5 \\ a_6 \end{bmatrix} \;\middle|\; a_i \in L,\ 1 \leq i \leq 6 \right\} \subseteq V \text{ and}$$

$$W_3 = \left\{ \begin{bmatrix} a_1 & a_2 & a_3 & \ldots & a_{10} \\ a_{11} & a_{12} & a_{13} & \ldots & a_{20} \end{bmatrix} \;\middle|\; a_i \in L,\ 1 \leq i \leq 20 \right\} \subseteq V,$$

$W_1$, $W_2$ and $W_3$ are set fuzzy linguistic vector subspaces of V over the set L.

Further $W_i \cap W_j = \phi$ if $i \neq j$, $1 \leq i, j \leq 3$ and $V = W_1 + W_2 + W_3$. Thus V is a direct sum of set fuzzy linguistic vector subspaces over L.

*Example 3.25*: Let

$$V = \left\{ \begin{bmatrix} a_1 & a_2 & \ldots & a_6 \\ a_7 & a_8 & \ldots & a_{12} \\ a_{13} & a_{14} & \ldots & a_{18} \\ a_{19} & a_{20} & \ldots & a_{24} \end{bmatrix}, \begin{bmatrix} a_1 & a_4 & \ldots & a_{28} \\ a_2 & a_5 & \ldots & a_{29} \\ a_3 & a_6 & \ldots & a_{30} \end{bmatrix}, \begin{bmatrix} a_1 & a_2 & a_3 & a_4 & a_5 \\ a_6 & a_7 & a_8 & a_9 & a_{10} \\ \vdots & \vdots & \vdots & \vdots & \vdots \\ a_{96} & a_{97} & a_{98} & a_{99} & a_{100} \end{bmatrix} \right.$$

$$a_i \in L, \ 1 \le i \le 100\};$$

L a fuzzy linguistic space associated with some problem.

V is a set fuzzy linguistic space of fuzzy linguistic matrices over the set S = L.

Consider

$$W_1 = \left\{ \begin{bmatrix} a_1 & a_2 & \ldots & a_6 \\ a_7 & a_8 & \ldots & a_{12} \\ a_{13} & a_{14} & \ldots & a_{18} \\ a_{19} & a_{20} & \ldots & a_{24} \end{bmatrix} \middle| a_i \in L, \ 1 \le i \le 24 \right\} \subseteq V,$$

$$W_2 = \left\{ \begin{bmatrix} a_1 & a_4 & \ldots & a_{28} \\ a_2 & a_5 & \ldots & a_{29} \\ a_3 & a_6 & \ldots & a_{30} \end{bmatrix} \middle| a_i \in L, \ 1 \le i \le 30 \right\} \subseteq V$$

and

$$W_3 = \left\{ \begin{bmatrix} a_1 & a_2 & a_3 & a_4 & a_5 \\ a_6 & a_7 & a_8 & a_9 & a_{10} \\ \vdots & \vdots & \vdots & \vdots & \vdots \\ a_{96} & a_{97} & a_{98} & a_{99} & a_{100} \end{bmatrix} \middle| a_i \in L, \ 1 \le i \le 100 \right\} \subseteq V$$

be the set fuzzy linguistic vector subspace of V over the set L.

We see $W_i \cap W_j = \phi$, $1 \leq i, j \leq 3$ and $V = W_1 + W_2 + W_3$ is a direct sum of set fuzzy linguistic vector subspaces of V over the fuzzy linguistic set L.

*Example 3.26:* Let

$$V = \left\{ (a_1, a_2, a_3), \begin{bmatrix} a_1 \\ a_2 \\ a_3 \\ a_4 \end{bmatrix}, \begin{bmatrix} a_1 & a_2 & a_3 \\ a_4 & a_5 & a_6 \\ a_7 & a_8 & a_9 \end{bmatrix}, \begin{bmatrix} a_1 & a_2 & a_3 & a_4 & a_5 \\ a_6 & a_7 & a_8 & a_9 & a_{10} \\ a_{11} & a_{12} & a_{13} & a_{14} & a_{15} \\ a_{16} & a_{17} & a_{18} & a_{19} & a_{20} \\ a_{21} & a_{22} & a_{23} & a_{24} & a_{25} \end{bmatrix} \right.$$

$$\left. a_i \in L, 1 \leq i \leq 25 \right\}$$

be a set fuzzy linguistic vector space over the fuzzy linguistic set L.

Let

$$W_1 = \left\{ (a_1, a_2, a_3), \begin{bmatrix} a_1 & a_2 & a_3 \\ a_4 & a_5 & a_6 \\ 0 & 0 & 0 \end{bmatrix} \middle| a_i \in L, 1 \leq i \leq 6 \right\} \subseteq V,$$

$$W_2 = \left\{ \begin{bmatrix} a_1 \\ a_2 \\ a_3 \\ a_4 \end{bmatrix}, \begin{bmatrix} a_1 & a_2 & a_3 \\ a_4 & a_5 & a_6 \\ a_7 & a_8 & a_9 \end{bmatrix} \middle| a_i \in L, 1 \leq i \leq 9 \right\} \subseteq V$$

and

$$W_3 = \left\{ \begin{bmatrix} a_1 & a_2 & a_3 \\ a_4 & a_5 & a_6 \\ a_7 & a_8 & a_9 \end{bmatrix}, \begin{bmatrix} a_1 & a_2 & a_3 & a_4 & a_5 \\ a_6 & a_7 & a_8 & a_9 & a_{10} \\ a_{11} & a_{12} & a_{13} & a_{14} & a_{15} \\ a_{16} & a_{17} & a_{18} & a_{19} & a_{20} \\ a_{21} & a_{22} & a_{23} & a_{24} & a_{25} \end{bmatrix} \middle| a_i \in L, \right.$$

$$\left. 1 \leq i \leq 25 \right\} \subseteq V$$

be set fuzzy linguistic subspaces of V over the fuzzy linguistic set L.

We see $W_i \cap W_j \neq \phi$, $1 \leq i, j \leq 3$ and $V \subseteq W_1 + W_2 + W_3$, V is only a pseudo direct sum of set fuzzy linguistic subspaces of V over the fuzzy linguistic set L.

Now having seen examples of these structures we now proceed onto define the notion of set linear transformation of set fuzzy linguistic vector spaces.

In the first place if we wish to define set linear transformation of two set fuzzy linguistic vector spaces V and W then both V and W must be defined only over the same set fuzzy linguistic set otherwise the set linear transformation cannot be defined.

**DEFINITION 3.7:** *Let V and W be any two set fuzzy linguistic vector spaces over the fuzzy linguistic set L. A map T from V to W is said to be a set linear transformation if*

$$T(v) = w$$
$$T(min(s, v)) = min\{s, w\} = min\{s, T(v)\}$$

*for all $v \in V$ and $s \in S$ and $w \in W$.*

*If we take $W = V$ itself then the set linear transformation is defined as the set linear operator on V.*

We will illustrate first this situation by some examples.

*Example 3.27*: Let

$$V = \left\{ (a_1, a_2, \ldots, a_6), \begin{bmatrix} a_1 \\ a_2 \\ \vdots \\ a_{12} \end{bmatrix}, \begin{pmatrix} a_1 & a_2 & \ldots & a_{20} \\ a_{21} & a_{22} & \ldots & a_{40} \\ a_{41} & a_{42} & \ldots & a_{60} \end{pmatrix} \middle| \; a_i \in L, 1 \leq i \leq 6; \right.$$

L a fuzzy linguistic space} and

$$W = \left\{ \begin{bmatrix} a_1 \\ a_2 \\ \vdots \\ a_6 \end{bmatrix}, (a_1, a_2, \ldots, a_{12}), \begin{pmatrix} a_1 & a_2 & \ldots & a_{10} \\ a_{11} & a_{12} & \ldots & a_{20} \\ a_{21} & a_{22} & \ldots & a_{30} \end{pmatrix} \middle| \; a_i \in L, \right.$$

$$1 \le i \le 6 \}$$

be two set fuzzy linguistic vector spaces defined over the linguistic set $S = L$.

Define a map $T : V \to W$ by

$$T((a_1, a_2, \ldots, a_6)) = \begin{bmatrix} a_1 \\ a_2 \\ \vdots \\ a_6 \end{bmatrix},$$

$$T\left( \begin{bmatrix} a_1 \\ a_2 \\ \vdots \\ a_{12} \end{bmatrix} \right) = (a_1, a_2, \ldots, a_{12}) \text{ and}$$

$$T\left( \begin{pmatrix} a_1 & a_2 & \ldots & a_{20} \\ a_{21} & a_{22} & \ldots & a_{40} \\ a_{41} & a_{42} & \ldots & a_{60} \end{pmatrix} \right) =$$

$$\begin{pmatrix} a_1 & a_3 & a_5 & a_7 & a_9 & a_{11} & a_{13} & a_{15} & a_{17} & a_{19} \\ a_{21} & a_{23} & a_{25} & a_{27} & a_{29} & a_{31} & a_{33} & a_{35} & a_{37} & a_{39} \\ a_{41} & a_{43} & a_{45} & a_{47} & a_{49} & a_{51} & a_{53} & a_{55} & a_{57} & a_{59} \end{pmatrix}.$$

T is a set linear transformation of V to W over the fuzzy linguistic set L.

*Example 3.28:* Let

$$V = \left\{ \begin{bmatrix} a_1 & a_2 \\ a_3 & a_4 \end{bmatrix}, \begin{bmatrix} a_1 \\ a_2 \\ \vdots \\ a_8 \end{bmatrix}, (a_1, a_2, \ldots, a_8), \begin{bmatrix} a_1 & a_2 & a_3 & a_4 \\ a_5 & a_6 & a_7 & a_8 \\ a_9 & a_{10} & a_{11} & a_{12} \end{bmatrix} \;\middle|\; a_i \in L, \right.$$

$1 \le i \le 12$; L a fuzzy linguistic space}

be a set fuzzy linguistic vector space over the fuzzy linguistic set L.

Define a map T from V to V as follows.

$$T\left( \begin{bmatrix} a_1 & a_2 \\ a_3 & a_4 \end{bmatrix} \right) = (a_1, 0, a_2, 0, a_3, 0, a_4, 0),$$

$$T((a_1, a_2, \ldots, a_8)) = \begin{bmatrix} a_1 \\ a_2 \\ \vdots \\ a_8 \end{bmatrix},$$

$$T\left( \begin{bmatrix} a_1 \\ a_2 \\ \vdots \\ a_8 \end{bmatrix} \right) = (a_1, a_2, \ldots, a_8) \text{ and}$$

$$T\left( \begin{bmatrix} a_1 & a_2 & a_3 & a_4 \\ a_5 & a_6 & a_7 & a_8 \\ a_9 & a_{10} & a_{11} & a_{12} \end{bmatrix} \right) = \begin{bmatrix} a_1 & 0 & 0 & a_2 \\ 0 & 0 & 0 & 0 \\ a_3 & 0 & 0 & a_4 \end{bmatrix}.$$

T is a set fuzzy linear operator on the set fuzzy linguistic space V over L.

We can also define the notion of inverse of a set linear transformation T where T is a set linear transformation from V into W, V and W set fuzzy linguistic vector spaces over a fuzzy linguistic set L. If T is invertible the inverse map $T^{-1}$ is a set linear transformation from W onto V.

We also define set linear functional. Let V be a set fuzzy linguistic vector space over a fuzzy linguistic set S. A map f from V on the set S is called a set linguistic linear functional if $f(v) = s$ for every $v \in V$ and s in S, such that $f(sv) = s\,f(v)$ and $f(v) \in S$.

We will illustrate this situation by some simple examples.

*Example 3.29:* Let

$$V = \left\{ \begin{bmatrix} a_1 & a_2 & a_3 \\ a_4 & a_5 & a_6 \\ a_7 & a_8 & a_9 \end{bmatrix}, \begin{pmatrix} a_1 & a_2 & \ldots & a_8 \\ a_9 & a_{10} & \ldots & a_{16} \end{pmatrix}, \begin{bmatrix} a_1 & a_2 & a_3 \\ a_4 & a_5 & a_6 \\ \vdots & \vdots & \vdots \\ a_{13} & a_{14} & a_{15} \end{bmatrix} \middle| a_i \in L, 1 \le i \le 16 \right\}$$

be a set fuzzy linguistic vector space over the fuzzy linguistic set $S = L$.

Define $\eta : V \to L$ by

$$\eta \left( \begin{bmatrix} a_1 & a_2 & a_3 \\ a_4 & a_5 & a_6 \\ a_7 & a_8 & a_9 \end{bmatrix} \right) = \min \{a_1, a_5, a_9\}$$

$$\eta \left( \begin{pmatrix} a_1 & a_2 & \ldots & a_8 \\ a_9 & a_{10} & \ldots & a_{16} \end{pmatrix} \right) = \min \{a_1, a_{16}\}$$

and

$$\eta \left( \begin{bmatrix} a_1 & a_2 & a_3 \\ a_4 & a_5 & a_6 \\ \vdots & \vdots & \vdots \\ a_{13} & a_{14} & a_{15} \end{bmatrix} \right) = \min \{a_2, a_4, a_6, a_8, a_{10}, a_{12}, a_{14}\}.$$

Clearly $\eta$ is a linear functional on V.

Now we proceed onto define the notion of set fuzzy linguistic linear algebra defined on a fuzzy linguistic set $S = L$.

Let V be a set fuzzy linguistic vector space over a fuzzy linguistic set L. If on V we can define a min (or max) operation then we define V to be a set fuzzy linguistic linear algebra over L.

We will first illustrate this situation by some examples.

*Example 3.30:* Let

$$V = \left\{ \begin{bmatrix} a_1 & a_2 & a_3 & a_4 \\ a_5 & a_6 & a_7 & a_8 \\ \vdots & \vdots & \vdots & \vdots \\ a_{21} & a_{22} & a_{23} & a_{24} \end{bmatrix} \;\middle|\; a_i \in L,\ 1 \le i \le 24 \right\},$$

L a fuzzy linguistic space. V is a set fuzzy linguistic linear algebra over the fuzzy linguistic set L.

*Example 3.31:* Let $V = \{(a_1, a_2, \ldots, a_{10}) \mid a_i \in L,\ 1 \le i \le 10\}$, L a fuzzy linguistic space. V is a set fuzzy linguistic linear algebra over the set L.

We have the following interesting theorem.

**THEOREM 3.2:** *Let V be a set fuzzy linguistic linear algebra over the fuzzy linguistic set L. V is a set fuzzy linguistic vector space over L, however a set fuzzy linguistic vector space V defined on L need not in general be a set fuzzy linguistic linear algebra over L.*

*Proof:* Every set fuzzy linguistic linear algebra V over the fuzzy linguistic set L is a set fuzzy linguistic vector space over L by the very definition of set fuzzy linguistic linear algebra.

However to prove the converse we illustrate the situation by an example.

Consider V =

$$\left\{ \begin{pmatrix} a_1 & a_2 & \ldots & a_{10} \\ a_{11} & a_{12} & \ldots & a_{20} \\ a_{21} & a_{22} & \ldots & a_{30} \end{pmatrix}, \begin{bmatrix} a_1 & a_2 & a_3 \\ a_4 & a_5 & a_6 \\ a_7 & a_8 & a_9 \end{bmatrix}, (a_1, a_2, \ldots, a_{25}), \begin{bmatrix} a_1 \\ a_2 \\ a_3 \\ a_4 \end{bmatrix}, \begin{bmatrix} a_1 & a_2 & a_3 & a_4 \\ 0 & a_5 & a_6 & a_7 \\ 0 & 0 & a_8 & a_9 \\ 0 & 0 & 0 & a_{10} \end{bmatrix} \right\}$$

$a_i \in L$, $1 \leq i \leq 30$; L a fuzzy linguistic space} be a set fuzzy linguistic vector space over the fuzzy linguistic set L. Clearly on V we cannot define 'min' or 'max' function. Hence the claim.

*Example 3.32:* Let

$$V = \left\{ \begin{bmatrix} a_1 & a_2 & a_3 \\ a_4 & a_5 & a_6 \\ \vdots & \vdots & \vdots \\ a_{16} & a_{17} & a_{18} \end{bmatrix} \,\middle|\, a_i \in L, \right.$$

L a fuzzy linguistic space, $1 \leq i \leq 18$} be a set fuzzy linguistic vector space, which is a set fuzzy linguistic linear algebra over the fuzzy linguistic set L.

*Example 3.33:* Let

$$V = \left\{ \begin{bmatrix} a_1 \\ a_2 \\ \vdots \\ a_{15} \end{bmatrix} \middle| a_i \in L, 1 \leq i \leq 15, \right.$$

L a fuzzy linguistic space; closed under 'min' operation} be a set fuzzy linguistic linear algebra over L, the fuzzy linguistic set.

We can define set fuzzy linguistic linear subalgebra which is a matter of routine and hence is left as an exercise to the reader. However we give examples of them.

*Example 3.34:* Let

$$V = \left\{ \begin{bmatrix} a_1 & a_2 \\ a_3 & a_4 \\ a_5 & a_6 \\ a_7 & a_8 \\ a_9 & a_{10} \\ a_{11} & a_{12} \\ a_{13} & a_{14} \\ a_{15} & a_{16} \end{bmatrix} \middle| a_i \in L, 1 \leq i \leq 16 \right\}$$

be a set fuzzy linguistic linear algebra over the fuzzy linguistic subset S of L.

Suppose 'm' be the least element in L under min operation. (m can be largest negative value or zero). L is closed under the binary operation, min (or max).

Consider

$$W_1 = \left\{ \begin{bmatrix} a_1 & 0 \\ 0 & a_2 \\ a_3 & 0 \\ 0 & a_4 \\ a_5 & 0 \\ 0 & a_6 \\ a_7 & 0 \\ 0 & a_8 \end{bmatrix} \middle| a_i \in L, 1 \leq i \leq 8 \right\} \subseteq V;$$

$W_1$ is a set fuzzy linguistic linear subalgebra of V over the fuzzy linguistic subset S of L.

$$W_2 = \left\{ \begin{bmatrix} 0 & a_1 \\ a_2 & 0 \\ 0 & a_3 \\ a_4 & 0 \\ 0 & a_5 \\ a_6 & 0 \\ 0 & a_7 \\ a_8 & 0 \end{bmatrix} \middle| a_i \in L, 1 \leq i \leq 8 \right\} \subseteq V;$$

$W_2$ is a set fuzzy linguistic linear subalgebra of V over the fuzzy linguistic subset S of L.

We see $W_1 \cap W_2 = \begin{bmatrix} 0 & 0 \\ 0 & 0 \\ 0 & 0 \\ 0 & 0 \\ 0 & 0 \\ 0 & 0 \\ 0 & 0 \\ 0 & 0 \end{bmatrix}$ and $V = W_1 + W_2$;

that V is a direct sum of set fuzzy linguistic linear subalgebras.

*Example 3.35:* Let

$$V = \left\{ \begin{bmatrix} a_1 & a_2 & a_3 & a_4 \\ a_5 & a_6 & a_7 & a_8 \\ a_9 & a_{10} & a_{11} & a_{12} \\ a_{13} & a_{14} & a_{15} & a_{16} \end{bmatrix} \middle| a_i \in L, 1 \leq i \leq 16 \right\}$$

be a set fuzzy linguistic linear algebra over the fuzzy linguistic set S, $S \subseteq L$.

Consider

$$M_1 = \left\{ \begin{bmatrix} a_1 & 0 & 0 & 0 \\ a_2 & 0 & 0 & 0 \\ a_3 & 0 & 0 & 0 \\ a_4 & 0 & 0 & 0 \end{bmatrix} \middle| a_i \in L, 1 \leq i \leq 4 \right\} \subseteq V,$$

$$M_2 = \left\{ \begin{bmatrix} 0 & a_2 & 0 & 0 \\ 0 & a_3 & 0 & 0 \\ 0 & a_4 & 0 & 0 \\ a_1 & a_5 & 0 & 0 \end{bmatrix} \middle| a_i \in L, 1 \leq i \leq 5 \right\} \subseteq V,$$

$$M_3 = \left\{ \begin{bmatrix} 0 & 0 & a_2 & 0 \\ 0 & 0 & a_3 & 0 \\ 0 & 0 & a_4 & 0 \\ a_1 & 0 & a_5 & 0 \end{bmatrix} \middle| a_i \in L, 1 \le i \le 5 \right\} \subseteq V,$$

and

$$M_4 = \left\{ \begin{bmatrix} 0 & 0 & 0 & a_2 \\ 0 & 0 & 0 & a_3 \\ 0 & 0 & 0 & a_4 \\ a_1 & 0 & 0 & a_5 \end{bmatrix} \middle| a_i \in L, 1 \le i \le 5 \right\} \subseteq V,$$

$M_1$, $M_2$, $M_3$ and $M_4$ are set fuzzy linguistic linear subalgebras of V. We see $M_i \cap M_j \ne (0)$ if $i \ne j$, $1 \le i, j \le 4$.

Further $M_1 + M_2 + M_3 + M_4 \supseteq V$, so V is a pseudo direct sum of set fuzzy linguistic linear subalgebras of V over the set $S \subseteq L$.

Now we have seen examples of direct sum and pseudo direct sum of set fuzzy linguistic linear subalgebras of V over the set $S \subseteq L$.

The notion of set linear transformation, set linear operator and set linear functionals on set fuzzy linguistic linear algebras are a matter of routine and hence is left as an exercise to the reader.

Next we proceed onto define the notion of semigroup fuzzy linguistic vector space over the semigroup S of fuzzy linguistic terms.

**DEFINITION 3.8:** *Let V be a fuzzy linguistic set. S a fuzzy linguistic semigroup with zero under 'min' operation. We call V to be a semigroup fuzzy linguistic vector space over the fuzzy linguistic semigroup S if the following conditions hold good.*

(i)     $\min(s, v) \in V$ for all $s \in S$ and $v \in V$.
(ii)     $\min(0, v) = 0 \in V$ for all $v \in V$ and $0 \in S$.
(iii)     $\min(\min\{s_1, s_2\}, v) = \min\{\min\{s_1, v\}, \min\{s_2, v\}\}$
for all $s_1, s_2 \in S$ and $v \in V$.

We give examples of semigroup fuzzy linguistic vector spaces over fuzzy linguistic semigroups.

*Example 3.36:* Let

$$V = \left\{ (a_1, a_2, a_3), \begin{bmatrix} a_1 \\ a_2 \\ a_3 \\ a_4 \end{bmatrix} \;\middle|\; a_i \in L,\; 1 \leq i \leq 4 \right\}$$

be a semigroup fuzzy linguistic vector space over the fuzzy linguistic semigroup L with '0' as its least element.

*Example 3.37:* Let

$$M = \left\{ \begin{bmatrix} a_1 & a_2 & a_3 & a_4 \\ a_5 & a_6 & a_7 & a_8 \\ a_9 & a_{10} & a_{11} & a_{12} \end{bmatrix}, \begin{bmatrix} a_1 & a_2 & a_3 \\ a_4 & a_5 & a_6 \\ a_7 & a_8 & a_9 \\ a_{10} & a_{11} & a_{12} \\ a_{13} & a_{14} & a_{15} \end{bmatrix}, \begin{bmatrix} a_1 & a_2 & \ldots & a_8 \\ a_9 & a_{10} & \ldots & a_{16} \\ \vdots & \vdots & & \vdots \\ a_{57} & a_{58} & \ldots & a_{64} \end{bmatrix} \right.$$

$$\left. a_i \in L,\; 1 \leq i \leq 64 \right\}$$

be a semigroup fuzzy linguistic vector space over the fuzzy linguistic semigroup L with '0' under min operation.

*Example 3.38:* Let

$$V = \left\{ \begin{bmatrix} a_1 \\ a_2 \\ \vdots \\ a_{10} \end{bmatrix}, \begin{bmatrix} a_1 & a_2 \\ a_3 & a_4 \end{bmatrix}, \begin{bmatrix} a_1 & a_2 & \dots & a_{10} \\ a_{11} & a_{12} & \dots & a_{20} \end{bmatrix} \;\middle|\; a_i \in L,\ 1 \le i \le 20 \right\}$$

be a semigroup fuzzy linguistic vector space over the fuzzy linguistic semigroup L for which '0' is the least element under 'min'. Let V be a semigroup fuzzy linguistic vector space over the semigroup L of fuzzy linguistic terms with '0' under min operation. Let $W \subseteq V$, W a proper subset of V, if W itself is a semigroup fuzzy linguistic vector space over the semigroup L we define W to be a semigroup fuzzy linguistic vector subspace of V over the semigroup L.

We will illustrate this situation by an example.

*Example 3.39:* Let

$$V = \left\{ \begin{bmatrix} a_1 & a_2 & a_3 \\ a_4 & a_5 & a_6 \end{bmatrix}, \begin{bmatrix} a_1 & a_2 \\ a_3 & a_4 \\ a_5 & a_6 \\ a_7 & a_8 \\ a_9 & a_{10} \\ a_{11} & a_{12} \end{bmatrix}, \begin{bmatrix} a_1 & a_2 & a_3 & a_4 & a_5 \\ a_6 & a_7 & a_8 & a_9 & a_{10} \\ a_{11} & a_{12} & a_{13} & a_{14} & a_{15} \\ a_{16} & a_{17} & a_{18} & a_{19} & a_{20} \\ a_{21} & a_{22} & a_{23} & a_{24} & a_{25} \end{bmatrix} \right.$$

$$a_i \in L,\ 1 \le i \le 25 \}$$

be a semigroup fuzzy linguistic vector space over the semigroup L.

$$P = \left\{ \begin{bmatrix} a_1 & 0 & a_2 \\ 0 & a_3 & 0 \end{bmatrix}, \begin{bmatrix} 0 & a_1 \\ a_2 & 0 \\ 0 & a_3 \\ a_4 & 0 \\ 0 & a_5 \\ a_6 & 0 \end{bmatrix}, \begin{bmatrix} a_1 & 0 & 0 & 0 & 0 \\ a_2 & a_3 & 0 & 0 & 0 \\ a_4 & a_5 & a_6 & 0 & 0 \\ a_7 & a_8 & a_9 & a_{10} & 0 \\ a_{11} & a_{12} & a_{13} & a_{14} & a_{15} \end{bmatrix} \right.$$

$$\left. a_i \in L, \ 1 \leq i \leq 15 \right\} \subseteq V$$

is a semigroup fuzzy linguistic vector subspace of V over the fuzzy linguistic semigroup L.

*Example 3.40:* Let

$$V = \left\{ \begin{bmatrix} a_1 & a_2 & a_3 & a_4 \\ a_5 & a_6 & a_7 & a_8 \\ a_9 & a_{10} & a_{11} & a_{12} \\ a_{13} & a_{14} & a_{15} & a_{16} \end{bmatrix}, (a_1, a_2, \ldots, a_{20}), \begin{bmatrix} a_1 & a_2 & \ldots & a_{12} \\ a_{13} & a_{14} & \ldots & a_{24} \end{bmatrix} \right.$$

$$\left. a_i \in L, \ 1 \leq i \leq 24 \right\}$$

be a semigroup fuzzy linguistic vector space over the fuzzy linguistic semigroup L.

Consider

$$W = \left\{ (a_1, a_2, \ldots, a_{20}), \begin{bmatrix} a_1 & a_2 & \ldots & a_{12} \\ a_{13} & a_{14} & \ldots & a_{24} \end{bmatrix} \ \middle| \ a_i \in L, \ 1 \leq i \leq 24 \right\}$$

$\subseteq$ V, W is a semigroup fuzzy linguistic vector subspace of V over the fuzzy linguistic semigroup L.

*Example 3.41:* Let

$$V = \left\{ \begin{bmatrix} a_1 & a_2 & a_3 \\ a_4 & a_5 & a_6 \\ \vdots & \vdots & \vdots \\ a_{28} & a_{29} & a_{30} \end{bmatrix}, \begin{bmatrix} a_1 & a_2 & a_3 & a_4 \\ a_5 & a_6 & a_7 & a_8 \\ a_9 & a_{10} & a_{11} & a_{12} \\ a_{13} & a_{14} & a_{15} & a_{16} \end{bmatrix}, \begin{bmatrix} a_1 & a_2 & \ldots & a_{10} \\ a_{11} & a_{12} & \ldots & a_{20} \\ a_{21} & a_{22} & \ldots & a_{30} \end{bmatrix} \right.$$

$$a_i \in L,\ 1 \leq i \leq 30 \}$$

be a semigroup fuzzy linguistic vector space over the fuzzy linguistic semigroup L.

$$M_1 = \left\{ \begin{bmatrix} a_1 & a_2 & a_3 \\ a_4 & a_5 & a_6 \\ \vdots & \vdots & \vdots \\ a_{28} & a_{29} & a_{30} \end{bmatrix} \;\middle|\; a_i \in L,\ 1 \leq i \leq 30 \right\} \subseteq V,$$

$$M_2 = \left\{ \begin{bmatrix} a_1 & a_2 & a_3 & a_4 \\ a_5 & a_6 & a_7 & a_8 \\ a_9 & a_{10} & a_{11} & a_{12} \\ a_{13} & a_{14} & a_{15} & a_{16} \end{bmatrix} \;\middle|\; a_i \in L,\ 1 \leq i \leq 16 \right\} \subseteq V$$

and

$$M_3 = \left\{ \begin{bmatrix} a_1 & a_2 & \ldots & a_{10} \\ a_{11} & a_{12} & \ldots & a_{20} \\ a_{21} & a_{22} & \ldots & a_{30} \end{bmatrix} \;\middle|\; a_i \in L,\ 1 \leq i \leq 30 \right\} \subseteq V$$

be semigroup fuzzy linguistic vector subspaces of V over the fuzzy linguistic semigroup L.

Clearly $M_i \cap M_j = \phi$ if $i \neq j$, $1 \leq i, j \leq 3$ and $V = M_1 + M_2 + M_3$. Thus V is a direct sum of semigroup fuzzy linguistic vector subspaces of V over the fuzzy linguistic semigroup L.

*Example 3.42:* Let

$$V = \left\{ \begin{bmatrix} a_1 & a_2 & a_3 \\ a_4 & a_5 & a_6 \\ \vdots & \vdots & \vdots \\ a_{28} & a_{29} & a_{30} \end{bmatrix}, \begin{bmatrix} a_1 & a_2 & a_3 & a_4 & a_5 \\ a_6 & \cdots & \cdots & \cdots & a_{10} \\ \vdots & \vdots & \vdots & \vdots & \vdots \\ a_{21} & a_{22} & a_{23} & a_{24} & a_{25} \end{bmatrix}, \begin{bmatrix} a_1 & a_2 & \cdots & a_9 \\ a_{10} & a_{11} & \cdots & a_{18} \\ a_{19} & a_{20} & \cdots & a_{27} \end{bmatrix} \right.$$

$$a_i \in L, \ 1 \leq i \leq 30 \}$$

be a semigroup fuzzy linguistic vector space over the fuzzy linguistic semigroup.

$$W_1 = \left\{ \begin{bmatrix} a_1 & a_2 & a_3 \\ a_4 & a_5 & a_6 \\ \vdots & \vdots & \vdots \\ a_{28} & a_{29} & a_{30} \end{bmatrix} \middle| a_i \in L, \ 1 \leq i \leq 30 \right\} \subseteq V,$$

$$W_2 = \left\{ \begin{bmatrix} 0 & 0 & 0 \\ a_1 & a_2 & a_3 \\ \vdots & \vdots & \vdots \\ a_{25} & a_{26} & a_{27} \end{bmatrix}, \begin{bmatrix} a_1 & a_2 & a_3 & a_4 & a_5 \\ a_6 & \cdots & \cdots & \cdots & a_{10} \\ \vdots & \vdots & \vdots & \vdots & \vdots \\ a_{21} & a_{22} & a_{23} & a_{24} & a_{25} \end{bmatrix} \right.$$

$$a_i \in L, \ 1 \leq i \leq 27 \} \subseteq V$$

and

$$W_3 = \left\{ \begin{bmatrix} 0 & 0 & 0 \\ 0 & 0 & 0 \\ a_1 & a_2 & a_3 \\ \vdots & \vdots & \vdots \\ a_{22} & a_{23} & a_{24} \end{bmatrix}, \begin{bmatrix} a_1 & a_2 & \ldots & a_9 \\ a_{10} & a_{11} & \ldots & a_{18} \\ a_{19} & a_{20} & \ldots & a_{27} \end{bmatrix} \right.$$

$$\left. a_i \in L,\ 1 \leq i \leq 27 \right\} \subseteq V$$

be a semigroup fuzzy linguistic vector subspaces of V over the fuzzy linguistic semigroup L.

Consider

$$W_1 \cap W_2 = \begin{bmatrix} 0 & 0 & 0 \\ a_1 & a_2 & a_3 \\ a_4 & a_5 & a_6 \\ \vdots & \vdots & \vdots \\ a_{25} & a_{26} & a_{27} \end{bmatrix},$$

$$W_1 \cap W_3 = \begin{bmatrix} 0 & 0 & 0 \\ 0 & 0 & 0 \\ a_1 & a_2 & a_3 \\ \vdots & \vdots & \vdots \\ a_{22} & a_{23} & a_{24} \end{bmatrix}$$

and

$$W_2 \cap W_3 = \begin{bmatrix} 0 & 0 & 0 \\ 0 & 0 & 0 \\ a_1 & a_2 & a_3 \\ \vdots & \vdots & \vdots \\ a_{22} & a_{23} & a_{24} \end{bmatrix};$$

further $V \subseteq W_1 + W_2 + W_3$ and hence V is only a pseudo direct sum of semigroup fuzzy linguistic vector subspaces of V over the fuzzy linguistic semigroup L.

We can define semigroup linear operator, semigroup linear transformation and semigroup linear functionals. All these are a matter of routine and we wish to give only examples of them.

*Example 3.43:* Let

$$V = \left\{ \begin{bmatrix} a_1 & a_2 & a_3 \\ a_4 & a_5 & a_6 \end{bmatrix}, \begin{bmatrix} a_1 & a_2 & \ldots & a_6 \\ a_7 & a_8 & \ldots & a_{12} \\ \vdots & \vdots & & \vdots \\ a_{31} & a_{32} & \ldots & a_{36} \end{bmatrix}, \begin{bmatrix} a_1 & a_2 & a_3 \\ a_4 & a_5 & a_6 \\ \vdots & \vdots & \vdots \\ a_{25} & a_{26} & a_{27} \end{bmatrix} \right.$$

$$\left. a_i \in L, 1 \leq i \leq 36 \right\}$$

and

$$W = \left\{ \begin{bmatrix} a_1 & a_2 \\ a_3 & a_4 \end{bmatrix}, \begin{bmatrix} a_1 & a_2 & a_3 & a_4 & \ldots & a_9 \\ a_{10} & a_{11} & a_{12} & a_{13} & \ldots & a_{18} \\ a_{19} & a_{20} & a_{21} & a_{22} & \ldots & a_{27} \\ a_{28} & a_{29} & a_{30} & a_{31} & \ldots & a_{36} \end{bmatrix}, \begin{bmatrix} a_1 & a_2 & a_3 & a_4 \\ a_5 & a_6 & a_7 & a_8 \\ \vdots & \vdots & & \vdots \\ a_{25} & a_{26} & a_{27} & a_{28} \end{bmatrix} \right.$$

$$\left. a_i \in L, 1 \leq i \leq 36 \right\}$$

be two semigroup fuzzy linguistic vector spaces defined over the fuzzy linguistic semigroup. L under 'min' operation with 0 as its least element.

Define a map $\eta : V \to W$ by

$$\eta \left( \begin{bmatrix} a_1 & a_2 & a_3 \\ a_4 & a_5 & a_6 \end{bmatrix} \right) = \begin{bmatrix} a_1 & a_2 \\ a_3 & a_4 \end{bmatrix}$$

$$\eta \left( \begin{bmatrix} a_1 & a_2 & \ldots & a_6 \\ a_7 & a_8 & \ldots & a_{12} \\ \vdots & \vdots & & \vdots \\ a_{31} & a_{32} & \ldots & a_{36} \end{bmatrix} \right) = \begin{bmatrix} a_1 & a_2 & a_3 & a_4 & \ldots & a_9 \\ a_{10} & a_{11} & a_{12} & a_{13} & \ldots & a_{18} \\ a_{19} & a_{20} & a_{21} & a_{22} & \ldots & a_{27} \\ a_{28} & a_{29} & a_{30} & a_{31} & \ldots & a_{36} \end{bmatrix}$$

and

$$\eta \left( \begin{bmatrix} a_1 & a_2 & a_3 \\ a_4 & a_5 & a_6 \\ \vdots & \vdots & \vdots \\ a_{25} & a_{26} & a_{27} \end{bmatrix} \right) = \begin{bmatrix} a_1 & a_2 & a_3 & a_4 \\ a_5 & a_6 & a_7 & a_8 \\ \vdots & \vdots & & \vdots \\ a_{25} & a_{26} & a_{27} & 0 \end{bmatrix}.$$

$\eta$ is a semigroup linear transformation of V into W.

*Example 3.44:* Let

$$V = \left\{ \begin{pmatrix} a_1 & a_2 \\ a_3 & a_4 \end{pmatrix}, (a_1, a_2, a_3, a_4, a_5, a_6), \begin{bmatrix} a_1 & a_2 \\ a_3 & a_4 \\ a_5 & a_6 \\ a_7 & a_8 \\ a_9 & a_{10} \\ a_{11} & a_{12} \end{bmatrix}, \begin{bmatrix} a_1 & a_2 & a_3 & a_4 \\ a_5 & a_6 & a_7 & a_8 \\ a_9 & a_{10} & a_{11} & a_{12} \\ a_{13} & a_{14} & a_{15} & a_{16} \end{bmatrix} \right.$$

$$\left. a_i \in L, 1 \leq i \leq 16 \right\}$$

be a semigroup fuzzy linguistic vector space over a fuzzy linguistic semigroup L with '0' under min operation.

Define $\eta : V \to V$ by

$$\eta\left(\begin{pmatrix} a_1 & a_2 \\ a_3 & a_4 \end{pmatrix}\right) = (a_1, a_2, a_3, 0, 0, a_4)$$

$$\eta((a_1, a_2, a_3, a_4, a_5, a_6)) = \begin{bmatrix} a_1 & 0 \\ 0 & a_2 \\ a_3 & 0 \\ 0 & a_4 \\ a_5 & 0 \\ 0 & a_6 \end{bmatrix}$$

$$\eta\left(\begin{bmatrix} a_1 & a_2 \\ a_3 & a_4 \\ a_5 & a_6 \\ a_7 & a_8 \\ a_9 & a_{10} \\ a_{11} & a_{12} \end{bmatrix}\right) = \begin{bmatrix} a_1 & a_2 & a_3 & a_4 \\ 0 & 0 & a_5 & a_6 \\ a_7 & a_8 & 0 & 0 \\ a_9 & a_{10} & a_{11} & a_{12} \end{bmatrix}$$

and

$$\eta\left(\begin{bmatrix} a_1 & a_2 & a_3 & a_4 \\ a_5 & a_6 & a_7 & a_8 \\ a_9 & a_{10} & a_{11} & a_{12} \\ a_{13} & a_{14} & a_{15} & a_{16} \end{bmatrix}\right) = \left(\begin{bmatrix} a_1 & a_2 \\ a_3 & a_4 \\ a_5 & a_6 \\ a_{11} & a_{12} \\ a_{13} & a_{14} \\ a_{15} & a_{16} \end{bmatrix}\right).$$

$\eta$ is a semigroup linear operator of the semigroup fuzzy linguistic vector space V defined over the fuzzy linguistic semigroup L.

*Example 3.45:* Let

$$V = \left\{ \begin{bmatrix} a_1 & a_2 & a_3 \\ a_4 & a_5 & a_6 \\ \vdots & \vdots & \vdots \\ a_{28} & a_{29} & a_{30} \end{bmatrix}, \begin{pmatrix} a_1 & a_2 & \ldots & a_{10} \\ a_{11} & a_{12} & \ldots & a_{20} \end{pmatrix}, \begin{bmatrix} a_1 & a_2 & a_3 \\ a_4 & a_5 & a_6 \\ a_7 & a_8 & a_9 \end{bmatrix} \right.$$

$$\left. a_i \in L, \, 1 \leq i \leq 30 \right\}$$

be a semigroup fuzzy linguistic vector space over the fuzzy linguistic semigroup L with 'min' operation and 0 as its least element.

Let $f : V \to L$ be a map such that

$$f \left( \begin{bmatrix} a_1 & a_2 & a_3 \\ a_4 & a_5 & a_6 \\ \vdots & \vdots & \vdots \\ a_{28} & a_{29} & a_{30} \end{bmatrix} \right) = \min \{a_1, a_2, \ldots, a_{30}\},$$

$$f \begin{pmatrix} a_1 & a_2 & \ldots & a_{10} \\ a_{11} & a_{12} & \ldots & a_{20} \end{pmatrix} = \min \{a_1, a_2, \ldots, a_{10}\}$$

and

$$f \left( \begin{bmatrix} a_1 & a_2 & a_3 \\ a_4 & a_5 & a_6 \\ a_7 & a_8 & a_9 \end{bmatrix} \right) = \min \{a_1, a_5, a_9\}.$$

f is a linear functional on V.

Now having seen almost all properties associated with semigroup fuzzy linguistic vector spaces defined over fuzzy linguistic semigroup with 'min' operation and '0' as its least element, we now proceed onto define semigroup fuzzy linguistic linear algebra over the fuzzy linguistic semigroup.

Let V be a semigroup fuzzy linguistic vector space over the fuzzy linguistic semigroup; if V is itself a semigroup under min operation with '0' as its minimum element then we define V to be a semigroup fuzzy linguistic linear algebra over the fuzzy linguistic semigroup.

We give some examples.

*Example 3.46*: Let

$$V = \left\{ \begin{bmatrix} a_1 & a_2 & a_3 \\ a_4 & a_5 & a_6 \\ a_7 & a_8 & a_9 \end{bmatrix} \middle| \ a_i \in L, \ 1 \leq i \leq 9 \right\}$$

be a semigroup fuzzy linguistic linear algebra over the fuzzy linguistic semigroup L with 'min' operation and '0' as its least element.

*Example 3.47:* Let

$$V = \left\{ \begin{bmatrix} a_1 & a_2 & ... & a_7 \\ a_8 & a_9 & ... & a_{14} \\ a_{15} & a_{16} & ... & a_{21} \end{bmatrix} \middle| \ a_i \in L, \ 1 \leq i \leq 21 \right\}$$

be a semigroup fuzzy linguistic linear algebra over the fuzzy linguistic semigroup L with min operation and '0' as its least element.

*Example 3.48:* Let

$$V = \left\{ \begin{bmatrix} a_1 & a_2 & a_3 & a_4 \\ a_5 & a_6 & a_7 & a_8 \\ \vdots & \vdots & \vdots & \vdots \\ a_{61} & a_{62} & a_{63} & a_{64} \end{bmatrix} \middle| \; a_i \in L, \; 1 \leq i \leq 64 \right\}$$

be a semigroup fuzzy linguistic linear algebra over the fuzzy linguistic semigroup L with min operation and '0' as its least element.

We can define the notion of semigroup fuzzy linguistic linear subalgebra over a fuzzy linguistic semigroup L and also subsemigroup fuzzy linguistic linear subalgebra over a fuzzy linguistic subsemigroup S of L.

This work is left as an exercise to the reader, but we give some examples of both these concepts.

*Example 3.49:* Let

$$V = \left\{ \begin{bmatrix} a_1 & a_2 & a_3 & a_4 \\ a_5 & a_6 & a_7 & a_8 \\ a_9 & a_{10} & a_{11} & a_{12} \\ a_{13} & a_{14} & a_{15} & a_{16} \end{bmatrix} \middle| \; a_i \in L, \; 1 \leq i \leq 16 \right\}$$

be a semigroup fuzzy linguistic linear algebra over the fuzzy linguistic semigroup L.

Suppose

$$W = \left\{ \begin{bmatrix} a_1 & a_2 & a_3 & a_4 \\ 0 & 0 & 0 & 0 \\ a_5 & a_6 & a_7 & a_8 \\ 0 & 0 & 0 & 0 \end{bmatrix} \,\middle|\, a_i \in L,\ 1 \le i \le 8 \right\} \subseteq V,$$

is a proper subset of V and W a semigroup fuzzy linguistic linear subalgebra of V over L.

*Example 3.50:* Consider the semigroup fuzzy linguistic linear algebra

$$V = \left\{ \begin{bmatrix} a_1 & a_2 & a_3 & a_4 & a_5 & a_6 & a_7 & a_8 \\ a_9 & a_{10} & a_{11} & a_{12} & a_{13} & a_{14} & a_{15} & a_{16} \end{bmatrix} \,\middle|\, a_i \in L, \right.$$

$$1 \le i \le 16 \}$$

over the semigroup L = {0, bad, very bad, good, very good, worst, best} under min operation.

Consider S = {0, bad, very bad, worst} ⊆ L, S is a subsemigroup of L.

Let

$$W = \left\{ \begin{bmatrix} a_1 & a_2 & a_3 & 0 & 0 & 0 & a_7 & 0 \\ 0 & 0 & 0 & a_4 & a_5 & a_6 & 0 & a_8 \end{bmatrix} \,\middle|\, a_i \in L, \right.$$

$$1 \le i \le 8 \} \subseteq V,$$

W is a subsemigroup fuzzy linguistic linear subalgebra of V over the subsemigroup S of L. We see if the fuzzy linguistic semigroup has no proper subsemigroup then we call the semigroup fuzzy linguistic linear algebra to be pseudo simple.

We now give examples of how a semigroup fuzzy linguistic linear algebra can be written as a direct sum and pseudo direct sum of semigroup fuzzy linguistic linear subalgebras.

*Example 3.51:* Let

$$V = \left\{ \begin{bmatrix} a_1 & a_2 & \ldots & a_{10} \\ a_{11} & a_{12} & \ldots & a_{20} \\ a_{21} & a_{22} & \ldots & a_{30} \\ a_{31} & a_{32} & \ldots & a_{40} \end{bmatrix} \middle| \; a_i \in L, \; 1 \leq i \leq 40 \right\}$$

be a semigroup fuzzy linguistic linear algebra over the fuzzy linguistic semigroup L.

Consider

$$W_1 = \left\{ \begin{bmatrix} a_1 & a_2 & \ldots & a_{10} \\ 0 & 0 & \ldots & 0 \\ 0 & 0 & \ldots & 0 \\ 0 & 0 & \ldots & 0 \end{bmatrix} \middle| \; a_i \in L, \; 1 \leq i \leq 10 \right\} \subseteq V,$$

$$W_2 = \left\{ \begin{bmatrix} 0 & 0 & \ldots & 0 \\ a_1 & a_2 & \ldots & a_{10} \\ 0 & 0 & \ldots & 0 \\ 0 & 0 & \ldots & 0 \end{bmatrix} \middle| \; a_i \in L, \; 1 \leq i \leq 10 \right\} \subseteq V,$$

$$W_3 = \left\{ \begin{bmatrix} 0 & 0 & \ldots & 0 \\ 0 & 0 & \ldots & 0 \\ a_1 & a_2 & \ldots & a_{10} \\ 0 & 0 & \ldots & 0 \end{bmatrix} \middle| \; a_i \in L, \; 1 \leq i \leq 10 \right\} \subseteq V,$$

and

$$W_4 = \left\{ \begin{bmatrix} 0 & 0 & ... & 0 \\ 0 & 0 & ... & 0 \\ 0 & 0 & ... & 0 \\ a_1 & a_2 & ... & a_{10} \end{bmatrix} \middle| \; a_i \in L, \; 1 \le i \le 10 \right\} \subseteq V,$$

be semigroup fuzzy linguistic linear subalgebras of V over the fuzzy linguistic semigroup L.

Clearly $V = W_1 + W_2 + W_3 + W_4$ and

$$W_i \cap W_j = \begin{bmatrix} 0 & 0 & ... & 0 \\ 0 & 0 & ... & 0 \\ 0 & 0 & ... & 0 \\ 0 & 0 & ... & 0 \end{bmatrix}, \; i \ne j, \; 1 \le i, j \le 4.$$

Thus V is a direct sum of semigroup fuzzy linguistic linear subalgebras over the fuzzy linguistic semigroup L.

*Example 3.52:* Let

$$V = \left\{ \begin{bmatrix} a_1 & a_2 & a_3 \\ a_4 & a_5 & a_6 \\ \vdots & \vdots & \vdots \\ a_{28} & a_{29} & a_{30} \end{bmatrix} \middle| \; a_i \in L, \; 1 \le i \le 30 \right\}$$

be a semigroup fuzzy linguistic linear algebra over the fuzzy linguistic semigroup L.

Take

$$W_1 = \left\{ \begin{bmatrix} 0 & 0 & 0 \\ a_1 & a_2 & a_3 \\ \vdots & \vdots & \vdots \\ 0 & 0 & 0 \end{bmatrix} \middle| \; a_i \in L, \; 1 \leq i \leq 3 \right\} \subseteq V,$$

$$W_2 = \left\{ \begin{bmatrix} a_1 & a_2 & a_3 \\ 0 & 0 & 0 \\ a_4 & a_5 & a_6 \\ 0 & 0 & 0 \\ \vdots & \vdots & \vdots \\ 0 & 0 & 0 \end{bmatrix} \middle| \; a_i \in L, \; 1 \leq i \leq 6 \right\} \subseteq V,$$

$$W_3 = \left\{ \begin{bmatrix} 0 & 0 & 0 \\ 0 & 0 & 0 \\ 0 & 0 & 0 \\ a_1 & a_2 & a_3 \\ 0 & 0 & 0 \\ \vdots & \vdots & \vdots \\ 0 & 0 & 0 \end{bmatrix} \middle| \; a_i \in L, \; 1 \leq i \leq 3 \right\} \subseteq V,$$

$$W_4 = \left\{ \begin{bmatrix} 0 & 0 & 0 \\ 0 & 0 & 0 \\ 0 & 0 & 0 \\ 0 & 0 & 0 \\ a_1 & a_2 & a_3 \\ a_4 & a_5 & a_6 \\ 0 & 0 & 0 \\ 0 & 0 & 0 \\ 0 & 0 & 0 \\ 0 & 0 & 0 \end{bmatrix} \,\middle|\, a_i \in L,\ 1 \le i \le 6 \right\} \subseteq V,$$

$$W_5 = \left\{ \begin{bmatrix} 0 & 0 & 0 \\ 0 & 0 & 0 \\ 0 & 0 & 0 \\ 0 & 0 & 0 \\ 0 & 0 & 0 \\ 0 & 0 & 0 \\ a_1 & a_2 & a_3 \\ 0 & 0 & 0 \\ 0 & 0 & 0 \\ 0 & 0 & 0 \end{bmatrix} \,\middle|\, a_i \in L,\ 1 \le i \le 6 \right\} \subseteq V,$$

and

$$W_6 = \left\{ \begin{bmatrix} 0 & 0 & 0 \\ \vdots & \vdots & \vdots \\ 0 & 0 & 0 \\ a_1 & a_2 & a_3 \\ a_4 & a_5 & a_6 \\ a_7 & a_8 & a_9 \end{bmatrix} \,\middle|\, a_i \in L,\ 1 \le i \le 9 \right\} \subseteq V,$$

be semigroup fuzzy linguistic linear subalgebras of V over the fuzzy linguistic semigroup L.

$$\text{Cleary } W_i \cap W_j = \begin{bmatrix} 0 & 0 & 0 \\ 0 & 0 & 0 \\ \vdots & \vdots & \vdots \\ 0 & 0 & 0 \end{bmatrix} \text{ if } i \neq j. \ 1 \leq i, j \leq 6 \text{ and}$$

$V = W_1 + W_2 + W_3 + W_4 + W_5 + W_6$.

Thus V is a direct sum of semigroup fuzzy linguistic linear subalgebras of V over L.

*Example 3.53:* Let

$$V = \left\{ \begin{bmatrix} a_1 & a_2 & ... & a_5 \\ a_6 & a_7 & ... & a_{10} \\ a_{11} & a_{12} & ... & a_{15} \\ a_{16} & a_{17} & ... & a_{20} \\ a_{21} & a_{22} & ... & a_{25} \end{bmatrix} \ \middle| \ a_i \in L, 1 \leq i \leq 25 \right\}$$

be a semigroup fuzzy linguistic linear algebra over the fuzzy linguistic semigroup L.

Consider

$$W_1 = \left\{ \begin{bmatrix} a_1 & a_2 & a_3 & a_4 & a_5 \\ a_6 & a_7 & 0 & 0 & 0 \\ 0 & 0 & 0 & 0 & 0 \\ 0 & 0 & 0 & 0 & 0 \\ 0 & 0 & 0 & 0 & 0 \end{bmatrix} \ \middle| \ a_i \in L, 1 \leq i \leq 7 \right\} \subseteq V$$

$$W_2 = \left\{ \begin{bmatrix} 0 & 0 & 0 & 0 & 0 \\ a_1 & a_2 & a_3 & a_4 & a_5 \\ a_6 & a_7 & 0 & 0 & 0 \\ 0 & 0 & 0 & 0 & 0 \\ 0 & 0 & 0 & 0 & 0 \end{bmatrix} \middle| a_i \in L,\ 1 \leq i \leq 7 \right\} \subseteq V$$

$$W_3 = \left\{ \begin{bmatrix} 0 & 0 & 0 & 0 & 0 \\ a_1 & a_2 & 0 & 0 & 0 \\ a_5 & a_6 & a_7 & a_8 & a_9 \\ 0 & 0 & 0 & 0 & 0 \\ 0 & 0 & 0 & 0 & 0 \end{bmatrix} \middle| a_1, a_2, a_5, a_6, a_7, a_8, a_9 \in L \right\} \subseteq V,$$

$$W_4 = \left\{ \begin{bmatrix} 0 & 0 & 0 & 0 & 0 \\ a_1 & a_2 & 0 & 0 & 0 \\ 0 & 0 & 0 & 0 & 0 \\ a_3 & a_4 & a_5 & a_6 & a_7 \\ 0 & 0 & 0 & 0 & 0 \end{bmatrix} \middle| a_i \in L,\ 1 \leq i \leq 7 \right\} \subseteq V$$

and

$$W_5 = \left\{ \begin{bmatrix} 0 & 0 & 0 & 0 & 0 \\ a_1 & a_2 & 0 & 0 & 0 \\ 0 & 0 & 0 & 0 & 0 \\ 0 & 0 & 0 & 0 & 0 \\ a_3 & a_4 & a_5 & a_6 & a_7 \end{bmatrix} \middle| a_i \in L,\ 1 \leq i \leq 7 \right\} \subseteq V$$

be semigroup fuzzy linguistic linear subalgebras of V over the fuzzy linguistic semigroup L.

Clearly $V \subseteq W_1 + W_2 + W_3 + W_4 + W_5$ and

$$W_i \cap W_j \neq \begin{bmatrix} 0 & 0 & 0 & 0 & 0 \\ 0 & 0 & 0 & 0 & 0 \\ 0 & 0 & 0 & 0 & 0 \\ 0 & 0 & 0 & 0 & 0 \\ 0 & 0 & 0 & 0 & 0 \end{bmatrix} \text{ if } i \neq j, 1 \leq i, j \leq 5.$$

Thus V is only pseudo direct sum of semigroup fuzzy linguistic linear subalgebras.

*Example 3.54:* Let

$$V = \left\{ \begin{bmatrix} a_1 & a_2 & a_3 & a_4 & a_5 & a_6 & a_7 & a_8 & a_9 & a_{10} \\ a_{11} & a_{12} & \dots & \dots & \dots & \dots & \dots & \dots & \dots & a_{20} \\ a_{21} & a_{22} & \dots & \dots & \dots & \dots & \dots & \dots & \dots & a_{30} \\ a_{31} & a_{32} & \dots & \dots & \dots & \dots & \dots & \dots & \dots & a_{40} \end{bmatrix} \middle| a_i \in L, \right.$$

$$1 \leq i \leq 40 \}$$

be a semigroup fuzzy linguistic linear algebra over the fuzzy linguistic semigroup L.

Consider

$$W_1 = \left\{ \begin{bmatrix} a_1 & a_5 & 0 & \dots & \dots & \dots & \dots & \dots & \dots & 0 \\ a_2 & a_6 & 0 & \dots & \dots & \dots & \dots & \dots & \dots & 0 \\ a_3 & 0 & 0 & \dots & \dots & \dots & \dots & \dots & \dots & 0 \\ a_4 & 0 & 0 & \dots & \dots & \dots & \dots & \dots & \dots & 0 \end{bmatrix} \middle| a_i \in L, \right.$$

$$1 \leq i \leq 6 \} \subseteq V,$$

$$W_2 = \left\{ \begin{bmatrix} 0 & a_1 & a_5 & 0 & \ldots & \ldots & \ldots & \ldots & 0 \\ 0 & a_2 & a_6 & 0 & \ldots & \ldots & \ldots & \ldots & 0 \\ 0 & a_3 & 0 & 0 & \ldots & \ldots & \ldots & \ldots & 0 \\ 0 & a_4 & 0 & 0 & \ldots & \ldots & \ldots & \ldots & 0 \end{bmatrix} \,\middle|\, a_i \in L, \right.$$

$$\left. 1 \leq i \leq 6 \right\} \subseteq V,$$

$$W_3 = \left\{ \begin{bmatrix} 0 & 0 & a_1 & a_5 & 0 & \ldots & \ldots & \ldots & 0 \\ 0 & 0 & a_2 & a_6 & 0 & \ldots & \ldots & \ldots & 0 \\ 0 & 0 & a_3 & 0 & 0 & \ldots & \ldots & \ldots & 0 \\ 0 & 0 & a_4 & 0 & 0 & \ldots & \ldots & \ldots & 0 \end{bmatrix} \,\middle|\, a_i \in L, \right.$$

$$\left. 1 \leq i \leq 6 \right\} \subseteq V,$$

$$W_4 = \left\{ \begin{bmatrix} 0 & 0 & 0 & a_1 & a_5 & 0 & \ldots & \ldots & 0 \\ 0 & 0 & 0 & a_2 & a_6 & 0 & \ldots & \ldots & 0 \\ 0 & 0 & 0 & a_3 & 0 & 0 & \ldots & \ldots & 0 \\ 0 & 0 & 0 & a_4 & 0 & 0 & \ldots & \ldots & 0 \end{bmatrix} \,\middle|\, a_i \in L, \right.$$

$$\left. 1 \leq i \leq 6 \right\} \subseteq V,$$

$$W_5 = \left\{ \begin{bmatrix} 0 & 0 & 0 & 0 & a_1 & a_5 & 0 & \ldots & 0 \\ 0 & 0 & 0 & 0 & a_2 & a_6 & 0 & \ldots & 0 \\ 0 & 0 & 0 & 0 & a_3 & 0 & 0 & \ldots & 0 \\ 0 & 0 & 0 & 0 & a_4 & 0 & 0 & \ldots & 0 \end{bmatrix} \,\middle|\, a_i \in L, \right.$$

$$\left. 1 \leq i \leq 6 \right\} \subseteq V,$$

$$W_6 = \left\{ \begin{bmatrix} 0 & 0 & 0 & 0 & 0 & a_1 & a_5 & 0 & \ldots & 0 \\ 0 & 0 & 0 & 0 & 0 & a_2 & a_6 & 0 & \ldots & 0 \\ 0 & 0 & 0 & 0 & 0 & a_3 & 0 & 0 & \ldots & 0 \\ 0 & 0 & 0 & 0 & 0 & a_4 & 0 & 0 & \ldots & 0 \end{bmatrix} \middle| a_i \in L, \right.$$

$1 \leq i \leq 6 \} \subseteq V$,

$$W_7 = \left\{ \begin{bmatrix} 0 & 0 & 0 & 0 & 0 & 0 & a_1 & a_5 & 0 & 0 \\ 0 & 0 & 0 & 0 & 0 & 0 & a_2 & a_6 & 0 & 0 \\ 0 & 0 & 0 & 0 & 0 & 0 & a_3 & 0 & 0 & 0 \\ 0 & 0 & 0 & 0 & 0 & 0 & a_4 & 0 & 0 & 0 \end{bmatrix} \middle| a_i \in L, \right.$$

$1 \leq i \leq 6 \} \subseteq V$,

$$W_8 = \left\{ \begin{bmatrix} 0 & 0 & 0 & 0 & 0 & 0 & 0 & a_1 & a_5 & 0 \\ 0 & 0 & 0 & 0 & 0 & 0 & 0 & a_2 & a_6 & 0 \\ 0 & 0 & 0 & 0 & 0 & 0 & 0 & a_3 & 0 & 0 \\ 0 & 0 & 0 & 0 & 0 & 0 & 0 & a_4 & 0 & 0 \end{bmatrix} \middle| a_i \in L, \right.$$

$1 \leq i \leq 6 \} \subseteq V$,

$$W_9 = \left\{ \begin{bmatrix} 0 & 0 & 0 & 0 & 0 & 0 & 0 & 0 & a_1 & a_5 \\ 0 & 0 & 0 & 0 & 0 & 0 & 0 & 0 & a_2 & a_6 \\ 0 & 0 & 0 & 0 & 0 & 0 & 0 & 0 & a_3 & 0 \\ 0 & 0 & 0 & 0 & 0 & 0 & 0 & 0 & a_4 & 0 \end{bmatrix} \middle| a_i \in L, \right.$$

$1 \leq i \leq 6 \} \subseteq V$

and

$$W_{10} = \left\{ \begin{bmatrix} a_5 & 0 & 0 & 0 & 0 & 0 & 0 & 0 & 0 & a_1 \\ a_6 & 0 & 0 & 0 & 0 & 0 & 0 & 0 & 0 & a_2 \\ 0 & 0 & 0 & 0 & 0 & 0 & 0 & 0 & 0 & a_3 \\ 0 & 0 & 0 & 0 & 0 & 0 & 0 & 0 & 0 & a_4 \end{bmatrix} \middle| \; a_i \in L, \right.$$

$$1 \leq i \leq 6 \} \subseteq V$$

be semigroup fuzzy linguistic linear subalgebras of V over the fuzzy linguistic semigroup.

Clearly $V \subseteq W_1 + W_2 + \ldots + W_9 + W_{10}$ and

$$W_i \cap W_j \neq \begin{bmatrix} 0 & 0 & 0 & 0 & 0 & 0 & 0 & 0 & 0 & 0 \\ 0 & 0 & 0 & 0 & 0 & 0 & 0 & 0 & 0 & 0 \\ 0 & 0 & 0 & 0 & 0 & 0 & 0 & 0 & 0 & 0 \\ 0 & 0 & 0 & 0 & 0 & 0 & 0 & 0 & 0 & 0 \end{bmatrix}$$

if $i \neq j$, $1 \leq i, j \leq 10$. V is only a pseudo direct sum of semigroup fuzzy linguistic linear subalgebras of V over L.

Let V be a semigroup fuzzy linguistic linear algebra over the semigroup L.

Suppose $W \subseteq V$ be a proper subset of V be such that W is a semigroup fuzzy linguistic vector subspace of V we call W a pseudo semigroup fuzzy linguistic vector subspace of V over the semigroup L.

We will illustrate this situation by some examples.

*Example 3.55:* Let

$$V = \left\{ \begin{bmatrix} a_1 & a_2 & a_3 & a_4 \\ a_5 & a_6 & a_7 & a_8 \\ a_9 & a_{10} & a_{11} & a_{12} \end{bmatrix} \middle| \; a_i \in L, 1 \leq i \leq 12 \right\}$$

be a semigroup fuzzy linguistic linear algebra over the fuzzy linguistic semigroup L.

Consider

$$W = \left\{ \begin{bmatrix} a_1 & 0 & 0 & 0 \\ 0 & 0 & 0 & 0 \\ 0 & 0 & 0 & 0 \end{bmatrix}, \begin{bmatrix} 0 & 0 & a_2 & a_3 \\ 0 & 0 & 0 & 0 \\ 0 & 0 & 0 & 0 \end{bmatrix}, \begin{bmatrix} 0 & 0 & 0 & 0 \\ 0 & 0 & 0 & 0 \\ a_1 & 0 & 0 & 0 \end{bmatrix}, \begin{bmatrix} 0 & 0 & 0 & 0 \\ 0 & a_1 & a_2 & 0 \\ 0 & 0 & 0 & 0 \end{bmatrix} \right.$$

where $a_1, a_2, a_3 \in L\} \subseteq V$

be a semigroup fuzzy linguistic vector space over L.

Clearly W is a pseudo semigroup fuzzy linguistic vector subspace of V over L.

*Example 3.56:* Let

$$V = \left\{ \begin{bmatrix} a_1 & a_2 & a_3 & a_4 \\ a_5 & a_6 & a_7 & a_8 \\ \vdots & \vdots & \vdots & \vdots \\ a_{45} & a_{46} & a_{47} & a_{48} \end{bmatrix} \; a_i \in L, 1 \leq i \leq 48 \right\}$$

be a semigroup fuzzy linguistic linear algebra over the fuzzy linguistic semigroup L.

Consider

$$P = \left\{ \begin{bmatrix} a_1 & 0 & a_2 & 0 \\ 0 & 0 & 0 & 0 \\ \vdots & \vdots & \vdots & \vdots \\ 0 & 0 & 0 & 0 \end{bmatrix}, \begin{bmatrix} 0 & 0 & 0 & 0 \\ 0 & a_1 & 0 & a_2 \\ 0 & 0 & 0 & 0 \\ \vdots & \vdots & \vdots & \vdots \\ 0 & 0 & 0 & 0 \end{bmatrix}, \begin{bmatrix} 0 & 0 & 0 & 0 \\ 0 & 0 & 0 & 0 \\ a_1 & 0 & 0 & 0 \\ a_2 & 0 & 0 & 0 \\ 0 & 0 & 0 & 0 \\ \vdots & \vdots & \vdots & \vdots \\ 0 & 0 & 0 & 0 \end{bmatrix}, \begin{bmatrix} 0 & 0 & 0 & 0 \\ 0 & 0 & 0 & 0 \\ \vdots & \vdots & \vdots & \vdots \\ 0 & 0 & 0 & 0 \\ a_1 & 0 & a_2 & 0 \end{bmatrix}, \right.$$

$$\left\{ \begin{bmatrix} 0 & 0 & 0 & 0 \\ 0 & 0 & 0 & 0 \\ a_1 & 0 & 0 & 0 \\ 0 & 0 & 0 & a_2 \\ 0 & 0 & 0 & 0 \\ \vdots & \vdots & \vdots & \vdots \\ 0 & 0 & 0 & 0 \end{bmatrix}, \begin{bmatrix} 0 & 0 & 0 & 0 \\ 0 & 0 & 0 & 0 \\ \vdots & \vdots & \vdots & \vdots \\ 0 & 0 & 0 & 0 \\ a_1 & 0 & 0 & a_2 \\ 0 & 0 & 0 & 0 \\ 0 & 0 & 0 & 0 \end{bmatrix} \,\middle|\, a_1, a_2 \in L \right\} \subseteq V$$

be a pseudo semigroup fuzzy linguistic vector subspace of V over the fuzzy linguistic semigroup L.

Now we proceed onto define the notion of polynomials with fuzzy linguistic coefficients.

Let L be a fuzzy linguistic space; $p(x) = \sum_{i=0}^{\infty} a_i x^i$; $a_i \in L$, $p(x)$ is defined as the fuzzy linguistic coefficient polynomial.

Let

$$V = \left\{ \sum_{i=0}^{\infty} a_i x^i \,\middle|\, a_i \in L \right\}$$

be the collection of all fuzzy linguistic coefficient polynomials. We cannot define addition or multiplication as in case of usual polynomials.

Consider $p(x), q(x) \in V$ where

$$p(x) = \sum_{i=0}^{\infty} a_i x^i \text{ and } q(x) = \sum_{i=0}^{\infty} b_i x^i \text{ with } a_i, b_i \in L, \text{ then}$$

$$\min\{p(x), q(x)\} = \sum_{i=0}^{\infty} \min\{a_i, b_i\} x^i = \sum_{i=0}^{\infty} c_i x^i \, ; c_i \in L.$$

This way min operation is performed on V.

Similarly we can define max operation on V.

Thus V can have either min operation or max operation. Using either min or max we find V to be a semigroup of fuzzy linguistic polynomials. Now if we use both min and max the fuzzy linguistic polynomials become a lattice of infinite order.

We can restrict the polynomial to a finite number that is

$$P = \left\{ \sum_{i=0}^{n} a_i x^i, n < \infty, a_i \in L \right\}.$$ Then also P is a semigroup under max (or min) binary operation.

*Example 3.57:* Let

$$V = \left\{ \sum_{i=0}^{8} a_i x^i \,\middle|\, a_i \in L;\ 0 \leq i \leq 8 \right\}$$

be semigroup of a fuzzy linguistic coefficient polynomial under min operation. Clearly V is of finite order.

*Example 3.58:* Let

$$M = \left\{ \sum_{i=0}^{4} a_i x^i \,\middle|\, a_i \in L;\ 0 \leq i \leq 4 \right\}$$

be a fuzzy linguistic coefficient polynomial semigroup under min operation of finite order.

*Example 3.59:* Let

$$T = \left\{ \sum_{i=0}^{\infty} a_i x^i \,\middle|\, a_i \in L \right\}$$

be a fuzzy linguistic coefficient polynomial semigroup under min operation of infinite order.

*Example 3.60:* Let

$$M = \left\{ \sum_{i=0}^{10} a_i x^i \,\middle|\, a_i \in L;\ 0 \leq i \leq 10 \right\}$$

be a fuzzy linguistic coefficient polynomial semigroup of finite order under max operation.

*Example 3.61:* Let

$$T = \left\{ \sum_{i=0}^{3} a_i x^i \,\middle|\, a_i \in L;\ 0 \leq i \leq 3 \right\}$$

be a fuzzy linguistic coefficient polynomial semigroup of finite order under the max operation.

*Example 3.62:* Let

$$T = \left\{ \sum_{i=0}^{\infty} a_i x^i \,\middle|\, a_i \in L \right\}$$

be a fuzzy linguistic coefficient polynomial semigroup under 'max' binary operation and T is an infinite ordered semigroup.

We just show how 'max' or 'min' operations are performed in such polynomials with fuzzy linguistic coefficients.

Let $p(x) = \text{good} + (\text{very bad})\,x + (\text{fair})x^2 + \text{best}\,x^3$ and

$q(x) = \text{better} + (\text{bad})x + (\text{very fair})x^2 + \text{good}\,x^4$ be two fuzzy linguistic polynomials.

Min $\{p(x), q(x)\}$ = min {good, better} + min {very bad, bad}$x$ + min {fair, very fair}$x^2$ + {best, 0}$x^3$ + {0, good}$x^4$.

= better + very bad $x$ + fair $x^2$.

Now for the same p(x) and q(x) the max operation functions in the following way.

Max {p(x), q(x)} = max {good, better} + max {very bad, bad}x + max {fair, very fair}$x^2$ + max {best, 0}$x^3$ + max {0, good}$x^4$

= good + bad x + very fair $x^2$ + best $x^3$ + good $x^4$.

Thus we can say in case of 'min' operation,

deg (p(x).q(x)) = deg (min (p(x), q(x)) ≤ min (deg p(x), deg q(x)).

In case of max operation

deg (p(x).q(x)) = deg (max p(x), q(x))
            = max (deg p(x), deg q(x)).

Thus we can see the difference while using min or max operation.

Now we can use these fuzzy linguistic polynomial coefficients to built set fuzzy linguistic vector spaces and semigroup fuzzy linguistic vector spaces. The definition of these concepts are a matter of routine, so we give only examples of them.

*Example 3.63:* Let L be a any fuzzy linguistic space.

$$P = \left\{ \sum_{i=0}^{2} a_i x^i, \sum_{i=0}^{14} a_i x^i, \sum_{i=0}^{9} a_i x^i \,\middle|\, a_i \in L \right\}$$

is a set fuzzy linguistic vector space of fuzzy linguistic coefficient polynomials over the set L.

*Example 3.64:* Let M = {p(x) = $a_0 + a_1x^3 + a_3x^{12} + a_2x^{10}$, $a_{10} + a_7x^8$, $a_1x^{16} + a_2x^{12} + a_3x^8 + a_4x^4$ | $a_i \in L$; $0 \le i \le 10$} be a set fuzzy linguistic vector space over the fuzzy linguistic set L.

*Example 3.65:* Let

$$M = \left\{ \sum_{i=0}^{13} a_i x^i, \sum_{i=0}^{3} a_i x^i, \sum_{i=0}^{27} a_i x^i \;\middle|\; a_i \in L \right\}$$

be a set fuzzy linguistic vector space over fuzzy linguistic set L.

*Example 3.66:* Let

$$M = \left\{ \sum_{i=0}^{25} a_i x^i \;\middle|\; a_i \in L;\ 0 \le i \le 25 \right\}$$

be a set fuzzy linguistic vector space over fuzzy linguistic set L. We see under min operation M is semigroup. So M is a set fuzzy linguistic linear algebra over the set L of finite order.

*Example 3.67:* Let

$$M = \left\{ \sum_{i=0}^{\infty} a_i x^i \;\middle|\; a_i \in L;\ 0 \le i \le 25 \right\}$$

be a set fuzzy linguistic set linear algebra over the set L. Clearly P is of infinite order.

We can now give examples of them.

*Example 3.68:* Let

$$M = \left\{ \sum_{i=0}^{9} a_i x^i \;\middle|\; a_i \in L;\ 0 \le i \le 9 \right\}$$

be a set fuzzy linguistic linear algebra over the fuzzy linguistic set L.

Consider

$$P = \left\{ \sum_{i=0}^{5} a_i x^i \,\middle|\, a_i \in L;\ 0 \leq i \leq 5 \right\} \subseteq M,$$

P is a set fuzzy linguistic linear subalgebra of M over the fuzzy linguistic set L.

*Example 3.69:* Let

$$V = \left\{ \sum_{i=0}^{20} a_i x^i \,\middle|\, a_i \in L;\ 0 \leq i \leq 20 \right\}$$

be a set fuzzy linguistic linear algebra over the fuzzy linguistic set L.

Consider

$$M = \left\{ \sum_{i=0}^{10} a_i x^i \,\middle|\, a_i \in L;\ 0 \leq i \leq 10 \right\},$$

M is a set fuzzy linguistic linear subalgebra of V over the fuzzy linguistic set L.

Now we can also find vector subspaces of V over the fuzzy linguistic set L.

*Example 3.70:* Let

$$V = \left\{ \sum_{i=2}^{8} a_i x^i, \sum_{i=8}^{10} a_i x^i, \sum_{i=9}^{21} a_i x^i, \sum_{i=20}^{25} a_i x^i \,\middle|\, a_i \in L \right\}$$

be a set fuzzy linguistic vector space of fuzzy linguistic coefficient polynomials over the fuzzy linguistic space L.

Consider

$$L_1 = \left\{ \sum_{i=2}^{8} a_i x^i \,\middle|\, a_i \in L; \ 0 \leq i \leq \infty \right\} \subseteq V,$$

$$L_2 = \left\{ \sum_{i=8}^{10} a_i x^i \,\middle|\, a_i \in L; \ 8 \leq i \leq 10 \right\} \subseteq V,$$

$$L_3 = \left\{ \sum_{i=9}^{21} a_i x^i \,\middle|\, a_i \in L; \ 9 \leq i \leq 21 \right\} \subseteq V$$

and

$$L_4 = \left\{ \sum_{i=20}^{25} a_i x^i \,\middle|\, a_i \in L; \ 20 \leq i \leq 25 \right\} \subseteq V,$$

be set fuzzy linguistic vector subspaces of V over the fuzzy linguistic set L.

Clearly $L_i \cap L_j = \phi$, $i \neq j$; for $1 \leq i, j \leq 4$. Further $V = L_1 + L_2 + L_3 + L_4$, thus V is the direct sum of set fuzzy linguistic subvector spaces over L.

*Example 3.71:* Let

$$V = \left\{ \sum_{i=0}^{8} a_i x^i, \sum_{i=10}^{19} a_i x^i, \sum_{i=18}^{25} a_i x^i, \sum_{i=12}^{28} a_i x^i \,\middle|\, a_i \in L \right\}$$

be a set fuzzy linguistic vector space over the fuzzy linguistic space L.

Consider

$$P_1 = \left\{ \sum_{i=0}^{8} a_i x^i \,\middle|\, a_i \in L \right\} \subseteq V,$$

$$P_2 = \left\{ \sum_{i=10}^{19} a_i x^i \,\middle|\, a_i \in L \right\} \subseteq V,$$

$$P_3 = \left\{ \sum_{i=0}^{6} a_i x^i, \sum_{i=18}^{25} a_i x^i \,\middle|\, a_i \in L \right\} \subseteq V$$

and

$$P_4 = \left\{ \sum_{i=12}^{28} a_i x^i, \sum_{i=18}^{25} a_i x^i \right\} \subseteq V.$$

Clearly $P_i$'s are set fuzzy linguistic vector subspaces of V and $P_i \cap P_j \neq \phi$ or $\{0\}$ if $i \neq j$; $1 \leq i, j \leq 4$.

Further $V \subseteq P_1 + P_2 + P_3 + P_4$, so V is a pseudo direct sum of set fuzzy linguistic vector subspaces of V over L.

*Example 3.72:* Let

$$V = \left\{ \begin{bmatrix} a_1 \\ a_2 \\ a_3 \end{bmatrix}, (a_1, a_2, a_3, a_4, a_5, a_6), \sum_{i=8}^{12} a_i x^i, \begin{bmatrix} a_1 & a_2 & a_3 & a_4 & a_5 \\ a_6 & a_7 & a_8 & a_9 & a_{10} \\ a_{11} & a_{12} & a_{13} & a_{14} & a_{15} \end{bmatrix}, \begin{bmatrix} a_1 & a_2 & a_3 \\ a_4 & a_5 & a_6 \\ a_7 & a_8 & a_9 \\ a_{10} & a_{11} & a_{12} \\ a_{13} & a_{14} & a_{15} \end{bmatrix} \,\middle|\, a_i \in L \right\}$$

be a set fuzzy linguistic vector space over the fuzzy linguistic set L.

Consider

$$M_1 = \left\{ \begin{bmatrix} a_1 \\ a_2 \\ a_3 \end{bmatrix} \middle| a_i \in L, \ 1 \leq i \leq 3 \right\} \subseteq V,$$

$M_2 = \{(a_1, a_2, a_3, a_4, a_5, a_6) \mid a_i \in L, \ 1 \leq I \leq 6\} \subseteq V,$

$$M_3 = \left\{ \sum_{i=8}^{12} a_i x^i \middle| a_i \in L \right\} \subseteq V,$$

$$M_4 = \left\{ \begin{bmatrix} a_1 & a_2 & a_3 & a_4 & a_5 \\ a_6 & a_7 & a_8 & a_9 & a_{10} \\ a_{11} & a_{12} & a_{13} & a_{14} & a_{15} \end{bmatrix} \middle| a_i \in L, \ 1 \leq i \leq 15 \right\} \subseteq V$$

and

$$M_5 = \left\{ \begin{bmatrix} a_1 & a_2 & a_3 \\ a_4 & a_5 & a_6 \\ a_7 & a_8 & a_9 \\ a_{10} & a_{11} & a_{12} \\ a_{13} & a_{14} & a_{15} \end{bmatrix} \middle| a_i \in L, \ 1 \leq i \leq 15 \right\} \subseteq V$$

be set fuzzy linguistic vector subspaces of V over the fuzzy linguistic set L.

Clearly $M_i \cap M_j = \phi$ if $i \neq j$, $1 \leq j, i \leq 5$.

We see $V = M_1 + M_2 + M_3 + M_4 + M_5$, the direct sum of set fuzzy linguistic vector subspaces of V over the fuzzy linguistic set L.

We can define linear transformation and linear operator of set fuzzy linguistic vector spaces over a fuzzy linguistic set L.

We will illustrate this by some examples.

*Example 3.72:* Let

$$V = \left\{ \sum_{i=0}^{20} a_i x^i, \sum_{i=24}^{30} a_i x^i, \sum_{i=100}^{150} a_i x^i \,\middle|\, a_i \in L \right\}$$

be a set fuzzy linguistic vector space over the fuzzy linguistic set L.

Let $f : V \to L$ be a map such that

$$f\left(\sum_{i=0}^{20} a_i x^i\right) = \min \{a_0, a_{20}\}$$

$$f\left(\sum_{i=24}^{30} a_i x^i\right) = \min \{a_{24}, a_{25}, \ldots, a_{30}\}$$

and

$$f\left(\sum_{i=100}^{150} a_i x^i\right) = \min \{a_{100}, a_{110}, a_{120}, a_{130}, a_{140}, a_{150}\}.$$

Clearly f is a set fuzzy linguistic linear functional on V.

*Example 3.73:* Let

$$V = \left\{ \sum_{i=0}^{25} a_i x^i \,\middle|\, a_i \in L \right\}$$

be a semigroup fuzzy linguistic linear algebra over the fuzzy linguistic semigroup L.

Define $f : V \to L$ by

$$f\left(\sum_{i=0}^{25} a_i x^i\right) = \min \{a_0, \ldots, a_{25}\}.$$

Define $g : V \to L$ by

$$g\left(\sum_{i=0}^{25} a_i x^i\right) = \min \{a_5, a_{10}, a_{15}, a_{20}, a_{25}\}.$$

It is easily verified f and g are two distinct semigroup fuzzy linguistic linear functionals on V.

*Example 3.74:* Let

$$V = \left\{\sum_{i=0}^{20} a_i x^i, \sum_{i=95}^{100} a_i x^i, \sum_{i=25}^{30} a_i x^i, \sum_{i=60}^{80} a_i x^i \,\middle|\, a_i \in L\right\}$$

be a set fuzzy linguistic vector space over the fuzzy linguistic space L.

Let T be a set linear operator on V defined by

$$T\left(\sum_{i=0}^{20} a_i x^i\right) = \sum_{i=60}^{80} a_i x^i$$

(by $a_0 \to a_{60}x^{60}$, $a_1 x \to a_{61}x^{61}$, $a_2 x^2 \to$

$a_{62}x^{62}, \ldots, a_{20}x^{20} \to a_{80}x^{80}$)

$$T\left(\sum_{i=25}^{30} a_i x^i\right) = \sum_{i=95}^{100} a_i x^i$$

($a_{25}x^{25} \to a_{95}x^{95}$, $a_{26}x^{26} \to a_{96}x^{96}$, $a_{27}x^{27} \to a_{97}x^{92}$,

$a_{27}x^{28} \to a_{98}x^{98}$, $a_{29}x^{29} \to a_{99}x^{99}$, $a_{30}x^{30} \to a_{100}x^{100}$).

$$T\left(\sum_{i=95}^{100} a_i x^i\right) = \sum_{i=25}^{30} a_i x^i$$

$(a_{95}x^{95} \to a_{25}x^{25}, a_{96}x^{96} \to a_{26}x^{26}, \ldots, a_{100}x^{100} \to a_{30}x^{30})$.

$$T\left(\sum_{i=60}^{80} a_i x^i\right) = \sum_{i=0}^{20} a_i x^i$$

$(a_{60}x^{60} \to a_0x^0, a_{61}x^{61} \to a_1x, a_{62}x^{62} \to a_2x^2, \ldots, a_{80}x^{80} \to a_{20}x^{20})$.

Thus T is a set linear fuzzy linguistic operator on V over the fuzzy linguistic set L.

*Example 3.75:* Let

$$V = \left\{\sum_{i=20}^{36} a_i x^i, \sum_{i=0}^{8} a_i x^i, \sum_{i=40}^{45} a_i x^i, \sum_{i=9}^{15} a_i x^i \,\middle|\, a_i \in L\right\}$$

be a set fuzzy linguistic vector space over the fuzzy linguistic set L.

$$W = \left\{\sum_{i=10}^{16} a_i x^i, \sum_{i=60}^{65} a_i x^i, \sum_{i=19}^{25} a_i x^i, \sum_{i=30}^{38} a_i x^i \,\middle|\, a_i \in L\right\}$$

be a set fuzzy linguistic vector space over the fuzzy linguistic set L.

Define a map T : V $\to$ W by

$$T\left(\sum_{i=20}^{36} a_i x^i\right) = \sum_{i=10}^{16} a_i x^i$$

$(a_{20}x^{20} \to a_{10}x^{10}, a_{21}x^{21} \to a_{11}x^{11}, \ldots)$

$$T\left(\sum_{i=0}^{8} a_i x^i\right) = \sum_{i=30}^{38} a_i x^i$$

$(a_0 x^0 = a_0 \to a_{30} x^{30},\ a_1 x \to a_{31} x^{31},\ \ldots,\ a_8 x^8 \to a_{38} x^{38})$,

$$T\left(\sum_{i=40}^{45} a_i x^i\right) = \sum_{i=60}^{65} a_i x^i$$

$(a_{40} x^{40} \to a_{60} x^{60},\ a_{41} x^{41} \to a_{61} x^{61},\ \ldots,\ a_{45} x^{45} \to a_{65} x^{45})$,

and

$$T\left(\sum_{i=9}^{15} a_i x^i\right) = \sum_{i=19}^{25} a_i x^i$$

$(a_9 x^9 \to a_{19} x^{19},\ a_{10} x^{10} \to a_{20} x^{20},\ \ldots,\ a_{15} x^{15} \to a_{25} x^{25})$.

It is easily verified T is a set fuzzy linguistic linear transformation of V into W.

*Example 3.76:* Let

$$V = \left\{ \sum_{i=20}^{45} a_i x^i \,\middle|\, a_i \in L \right\}$$

and

$$W = \left\{ \sum_{i=0}^{25} a_i x^i \,\middle|\, a_i \in L \right\}$$

be a set fuzzy linguistic linear algebras over the fuzzy linguistic set L.

T : V → W be a map such that

$$T\left(\sum_{i=20}^{45} a_i x^i\right) = \sum_{i=0}^{25} a_i x^i,$$

with $a_{20}x^{20} \to a_0$, $a_{21}x^{21} \to a_1x$, $a_{22}x^{22} \to a_2x^2$,

..., $a_{45}x^{45} \to a_{25}x^{25}$.

Clearly T is a set fuzzy linguistic linear transformation from V to W.

Having seen examples of the concept of polynomials with fuzzy linguistic terms as coefficients we now proceed onto study polynomials with fuzzy matrix linguistic term as coefficients.

Let

$$V = \left\{\sum_{i=20}^{20} a_i x^i \,\middle|\, a_i = (x_1, x_2, \ldots, x_n) \mid x_j \in L; 1 \le j \le n \right.$$

$$\text{and } 0 \le i \le 20\};$$

denote the collection of all fuzzy linguistic row matrix coefficient polynomials.

We will show how 'min' and 'max' operations are performed on V.

Let $p(x) = (x_1^0, x_2^0, x_3^0) + (x_1^1, x_2^1, x_3^1)x + (x_1^2, x_2^2, x_3^2)x^2$
$+ (x_1^4, x_2^4, x_3^4)x^4$ and

$q(x) = (y_1^1, y_2^1, y_3^1)x + (y_1^2, y_2^2, y_3^2)x^2 + (y_1^3, y_2^3, y_3^3)x^3$
$+ (y_1^4, y_2^4, y_3^4)x^4$

be any two fuzzy linguistic row matrix coefficient polynomials.

$x_j^i$, $y_k^t \in L$, $1 \le j, k \le 3$ and $0 \le i \le 4$, $1 \le t \le 4$.

Now we can define min $(p(x), q(x))$

= (min $(x_1^1, y_1^1)$, min $(x_2^1, y_2^1)$, min $(x_3^1, y_3^1)$)x + (min $(x_1^2, y_1^2)$, min $(x_2^2, y_2^2)$, min $(x_3^2, y_3^2)$)$x^2$ + (0, 0)$x^3$ + ((min $(x_1^4, y_1^4)$, min $(x_2^4, y_2^4)$, min $(x_3^4, y_3^4)$)$x^4$ is again a fuzzy linguistic row matrix coefficient polynomial.

Consider $p(x)$ = (good, 0, average, bad) + (0, bad, good, better)x + (0, bad, 0, good)$x^2$ + (bad, average, 0, below average)$x^3$ + (good, good, bad, bad)$x^4$ and

$q(x)$ = (0, good, bad, good) + (bad, bad, good, average)x + (good, bad, best, best)$x^2$ + (average, 0, good, bad)$x^4$ be a fuzzy linguistic row matrix coefficients elements taken form a linguistic space L. (Here in V we take n = 4).

min $(p(x), q(x))$ = min {(good, 0, average, bad), (0, good, bad, good)} + min {(bad, bad, good, average), (0, bad, good, better)} x + min {(0, bad, 0, good), (good, bad, best, best)}$x^2$ + min {(bad, average, 0, below average), (0, 0, 0, 0)}$x^3$ + min {(good, good, bad, bad), (average, 0, good, bad)}$x^4$

= (0, 0, bad, bad) + (0, bad, good, average)x + (0, bad, 0, good)$x^2$ + (0, 0, 0, 0)$x^3$ + (average, 0, bad, bad)$x^4$.

Likewise we can define

max $(p(x), q(x))$ = (good, good, average, good) + (bad, bad, good, better)x + (good, bad, best, best)$x^2$ + (bad, average, 0, below average)$x^3$ + (good, good, good, bad)$x^4$ is the resultant. We see min $(p(x), q(x)) \ne$ max $(p(x), q(x))$.

However we can define operation min or max on the fuzzy linguistic row matrix coefficient polynomials and they under these operations form only a semigroup. When both the operations are performed the collection will be a lattice.

Now we have the following condition on the collection V of all fuzzy linguistic coefficient row matrix polynomial lattices. If V is finite we say V is a finite lattice. If V is infinite V will be an infinite lattice.

We can get several properties associated with them.

**THEOREM 3.9:** *Let*

$$V = \left\{ \sum_{i=0}^{n} a_i x^i \,\middle|\, n < \infty,\, a_i = (t_1, t_2, \ldots, t_m) \right.$$
$$\left. \text{with } t_j \in L,\, 1 \leq j \leq m,\, 0 \leq i \leq n \right\}$$

*be the collection of all fuzzy linguistic row coefficient polynomials. V is a lattice under min, max operation and V is a partially ordered set under the degree of a polynomial.*

This proof is simple and hence is left as an exercise to the reader.

We now give one or two examples before we proceed onto define fuzzy linguistic column matrix coefficient polynomial.

*Example 3.77:* Let

$$V = \left\{ \sum_{i=0}^{12} a_i x^i \,\middle|\, a_i = (s_1, s_2, s_3, s_4, s_5) \text{ with } s_j \in L;\, 1 \leq j \leq 5 \right.$$
$$\left. \text{and } 0 \leq i \leq 12 \right\}$$

be a fuzzy linguistic row matrix coefficient polynomial. Clearly V is a semigroup under 'min'; V is of finite dimension.

$(0,0,0,0,0)$ acts as the minimum element of V.

*Example 3.78:* Let

$$V = \left\{ \sum_{i=0}^{5} a_i x^i \;\middle|\; a_i = (t_1, t_2, \ldots, t_8) \text{ with } t_j \in L \right.$$

(L a fuzzy linguistic space (terms), $1 \leq j \leq 8$ and $0 \leq i \leq 5$} be a fuzzy linguistic row matrix coefficient polynomial semigroup with max operation of finite order.

*Example 3.79:* Let

$$M = \left\{ \sum_{i=0}^{\infty} a_i x^i \;\middle|\; a_i = (d_1, d_2, \ldots, d_{10}); d_j \in L, 1 \leq j \leq 10 \right\}$$

be a fuzzy linguistic row matrix coefficient polynomial semigroup with max (or min) operation. M is of infinite order and has zero divisors (only under min operation) ideals and subrings.

*Example 3.80:* Let

$$P = \left\{ \sum_{i=0}^{4} a_i x^i \;\middle|\; a_i = (d_1, d_2, d_3); d_j \in L, 1 \leq j \leq 3, 0 \leq i \leq 4 \right\}$$

be a fuzzy linguistic row matrix coefficient polynomial semigroup under min operation.

Consider

$$M_1 = \left\{ \sum_{i=0}^{4} a_i x^i \;\middle|\; a_i = (d, 0, 0); d \in L, 0 \leq i \leq 4 \right\} \subseteq P,$$

$$M_2 = \left\{ \sum_{i=0}^{4} a_i x^i \;\middle|\; a_i = (0, d, 0); d \in L, 0 \leq i \leq 4 \right\} \subseteq P$$

and

$$M_3 = \left\{ \sum_{i=0}^{4} a_i x^i \,\middle|\, a_i = (0, 0, d); d \in L, 0 \leq i \leq 4 \right\} \subseteq P$$

be three fuzzy linguistic row matrix coefficient polynomial semigroups under min operation.

Clearly for every $p(x) \in M_1$, and $q(x) \in M_2$ (or $q(x) \in M_3$) we have min $(p(x), q(x)) = (0,0,0)$.

Further $P = M_1 + M_2 + M_3$ and $M_i \cap M_j = (0,0,0)$, $i \neq j$, $1 \leq i, j \leq 3$. Thus P is the direct sum of fuzzy linguistic row matrix coefficient polynomial subsemigroups of P.

Consider

$$N_1 = \left\{ \sum_{i=0}^{4} a_i x^i \,\middle|\, a_i = (d_1, d_2, 0); d_1, d_2 \in L, 0 \leq i \leq 4 \right\} \subseteq P,$$

$$N_2 = \left\{ \sum_{i=0}^{4} a_i x^i \,\middle|\, a_i = (0, d_1, d_2); d_1, d_2 \in L, 0 \leq i \leq 4 \right\} \subseteq P$$

and

$$N_3 = \left\{ \sum_{i=0}^{4} a_i x^i \,\middle|\, a_i = (d_1, 0, d_2); d_1, d_2 \in L, 0 \leq i \leq 4 \right\} \subseteq P$$

be fuzzy linguistic row matrix coefficient polynomial subsemigroups under min operation. Clearly $P \subseteq N_1 + N_2 + N_3$ and $N_i \cap N_j \neq (0, 0, 0)$ if $i \neq j$; $1 \leq i, j \leq 3$. Thus P is only a pseudo direct sum of subsemigroups. Further we see in general for $p(x) \in N_i$ and $q(x) \in N_j$; $i \neq j$. min $(p(x), q(x)) \neq (0, 0, 0)$.

Also each of these subsemigroups are ideals of P, under min operation.

*Example 3.81:* Let

$$S = \left\{ \sum_{i=0}^{8} a_i x^i \;\middle|\; a_i = (d_1, d_2, d_3, d_4, d_5);\; d_j \in L,\; 1 \leq j \leq 5,\; 0 \leq i \leq 8 \right\}$$

be a fuzzy linguistic row matrix coefficient polynomial semigroups under max operation. This has no zero divisors but has only subsemigroups.

Consider

$$P_1 = \left\{ \sum_{i=0}^{8} a_i x^i \;\middle|\; a_i = (0, d_1, d_2, 0, 0);\; d_1, d_2 \in L,\; 0 \leq i \leq 8 \right\} \subseteq S,$$

$$P_2 = \left\{ \sum_{i=0}^{8} a_i x^i \;\middle|\; a_i = (d_1, d_2, 0, 0, 0);\; d_1, d_2 \in L,\; 0 \leq i \leq 8 \right\} \subseteq S,$$

$$P_3 = \left\{ \sum_{i=0}^{8} a_i x^i \;\middle|\; a_i = (d_1, 0, 0, d_2, 0);\; d_1, d_2 \in L,\; 0 \leq i \leq 8 \right\} \subseteq S \text{ and}$$

$$P_4 = \left\{ \sum_{i=0}^{8} a_i x^i \;\middle|\; a_i = (d_1, 0, d_3, 0, d_2);\; d_1, d_2, d_3 \in L,\; 0 \leq i \leq 8 \right\} \subseteq S,$$

proper fuzzy linguistic row matrix polynomial coefficients subsemigroups of S.

We see $S \subseteq P_1 + P_2 + P_3 + P_4 + P_5$, thus S is a pseudo direct sum of fuzzy linguistic row matrix polynomial coefficient subsemigroups.

Now we proceed onto define the notion of fuzzy linguistic column matrix coefficient semigroups under max (or min).

We give examples of them which is as follows.

Let

$$V = \left\{ \sum_{i=0}^{n} a_i x^i \;\middle|\; a_i = \begin{bmatrix} t_1 \\ t_2 \\ \vdots \\ t_m \end{bmatrix} \text{ where } t_j \in L,\; 1 \leq j \leq m \right\}$$

be a collection of all fuzzy linguistic column matrix coefficient polynomials.

V under min (or max) operation is a semigroup. If $n < \infty$, V is of finite order and if $n = \infty$, V is of infinite order; only under min operation V has ideals and zero divisors, under max operation V has no ideals or zero divisors. All these will be described by some examples.

*Example 3.82:* Let

$$V = \left\{ \sum_{i=0}^{8} a_i x^i \;\middle|\; a_i = \begin{bmatrix} t_1 \\ t_2 \\ \vdots \\ t_{20} \end{bmatrix} ;\; t_j \in L,\; 1 \leq j \leq 20,\; 0 \leq i \leq 8 \right\}$$

be a fuzzy linguistic column matrix coefficient polynomial semigroup under 'min'.

We see if $p(x) = \begin{bmatrix} t_1^0 \\ t_2^0 \\ \vdots \\ t_{20}^0 \end{bmatrix} + \begin{bmatrix} t_1^1 \\ t_2^1 \\ \vdots \\ t_{20}^1 \end{bmatrix} x + \begin{bmatrix} t_1^3 \\ t_2^3 \\ \vdots \\ t_{20}^3 \end{bmatrix} x^3 +$

and $q(x) = \begin{bmatrix} x_1^0 \\ x_2^0 \\ \vdots \\ x_{20}^0 \end{bmatrix} + \begin{bmatrix} x_1^1 \\ x_2^1 \\ \vdots \\ x_{20}^1 \end{bmatrix} x^2 + \begin{bmatrix} x_1^3 \\ x_2^3 \\ \vdots \\ x_{20}^3 \end{bmatrix} x^3 ;$

then

$\min(q(x), p(x)) = \begin{bmatrix} \min(x_1^0, t_1^0) \\ \min(x_2^0, t_2^0) \\ \vdots \\ \min(x_{20}^0, t_{20}^0) \end{bmatrix} + \begin{bmatrix} 0 \\ 0 \\ \vdots \\ 0 \end{bmatrix} x^2 + \begin{bmatrix} \min(x_1^3, t_1^3) \\ \min(x_2^3, t_2^3) \\ \vdots \\ \min(x_{20}^3, t_{20}^3) \end{bmatrix} x^3.$

*Example 3.83:* Let

$$V = \left\{ \sum_{i=0}^{5} a_i x^i \,\middle|\, a_i = \begin{bmatrix} t_1 \\ t_2 \\ t_3 \\ t_4 \end{bmatrix} ; t_j \in L, 1 \le j \le 4, \ 0 \le i \le 5 \right\}$$

be a fuzzy linguistic column matrix coefficient polynomial semigroup under max operation.

Clearly V has subsemigroups given by

$$M_1 = \left\{ \sum_{i=0}^{3} a_i x^i \,\middle|\, a_i = \begin{bmatrix} t_1 \\ t_2 \\ t_3 \\ t_4 \end{bmatrix} ; t_j \in L, 1 \le j \le 4, \ 0 \le i \le 3 \right\} \subseteq V.$$

Clearly $M_1$ is not an ideal of V. Further V has no zero divisors.

*Example 3.84:* Let

$$S = \left\{ \sum_{i=0}^{7} a_i x^i \;\middle|\; a_i = \begin{bmatrix} m_1 \\ m_2 \\ m_3 \\ m_4 \\ m_5 \end{bmatrix} ; m_j \in L,\; 1 \leq j \leq 5,\; 0 \leq i \leq 7 \right\}$$

be a fuzzy linguistic column matrix coefficient polynomial semigroup under min operation. S has subsemigroups, ideals and zero divisors. (If on the other hand 'min' operation on S is replaced by max operation, S has no ideals or zero divisors). However S has idempotents.

*Example 3.85:* Let

$$S = \left\{ \sum_{i=0}^{7} a_i x^i \;\middle|\; a_i = \begin{bmatrix} d_1 \\ d_2 \\ d_3 \\ d_4 \\ d_5 \\ d_6 \\ d_7 \end{bmatrix} ; d_j \in L,\; 1 \leq j \leq 7,\; 0 \leq i \leq 7 \right\}$$

be a fuzzy linguistic column matrix coefficient polynomial semigroup under 'min' operation.

Take

$$M_1 = \left\{ \sum_{i=0}^{7} a_i x^i \,\middle|\, a_i = \begin{bmatrix} d_1 \\ d_2 \\ 0 \\ 0 \\ \vdots \\ \vdots \\ 0 \end{bmatrix} ; d_1, d_2 \in L, 0 \leq i \leq 7 \right\} \subseteq P,$$

$$M_2 = \left\{ \sum_{i=0}^{7} a_i x^i \,\middle|\, a_i = \begin{bmatrix} 0 \\ 0 \\ d_1 \\ d_2 \\ 0 \\ 0 \\ 0 \end{bmatrix} ; d_1, d_2 \in L, 0 \leq i \leq 7 \right\} \subseteq P,$$

$$M_3 = \left\{ \sum_{i=0}^{7} a_i x^i \,\middle|\, a_i = \begin{bmatrix} 0 \\ 0 \\ 0 \\ 0 \\ d_1 \\ d_2 \\ 0 \end{bmatrix} ; d_1, d_2 \in L, 0 \leq i \leq 7 \right\} \subseteq P,$$

and

$$M_4 = \left\{ \sum_{i=0}^{7} a_i x^i \,\middle|\, a_i = \begin{bmatrix} 0 \\ 0 \\ 0 \\ 0 \\ 0 \\ 0 \\ d_1 \end{bmatrix} ; d_1, d_2 \in L, 0 \leq i \leq 7 \right\} \subseteq P,$$

be fuzzy linguistic column matrix coefficient polynomial subsemigroups of the semigroup P.

We see $P = M_1 + M_2 + M_3 + M_4$ is a direct sum of subsemigroups and $M_i \cap M_j = \begin{bmatrix} 0 \\ 0 \\ \vdots \\ 0 \end{bmatrix}$, if $i \neq j$, $1 \leq i, j \leq 4$. Also each of the $M_i$'s are also ideals of P; $1 \leq i \leq 4$.

Further for every $p(x) \in M_i$ we have for every $q(x) \in M_j$ for all $j \neq i$ we see

$$\min \{m_i, m_j\} = \begin{bmatrix} 0 \\ 0 \\ \vdots \\ 0 \end{bmatrix}.$$

Thus P has ideals, subsemigroups and zero divisors.

*Example 3.86:* Let

$$V = \left\{ \sum_{i=0}^{\infty} a_i x^i \;\middle|\; a_i = \begin{bmatrix} d_1 \\ d_2 \\ d_3 \\ d_4 \\ d_5 \\ d_6 \end{bmatrix}; d_j \in L, 1 \le j \le 6 \right\}$$

be a fuzzy linguistic column matrix coefficient polynomial semigroup under max.

Consider

$$P_1 = \left\{ \sum_{i=0}^{\infty} a_i x^i \;\middle|\; a_i = \begin{bmatrix} d_1 \\ 0 \\ 0 \\ 0 \\ 0 \\ 0 \\ d_2 \end{bmatrix}; d_1, d_2 \in L \right\} \subseteq V,$$

$$P_2 = \left\{ \sum_{i=0}^{\infty} a_i x^i \;\middle|\; a_i = \begin{bmatrix} d_1 \\ d_2 \\ 0 \\ 0 \\ 0 \\ 0 \end{bmatrix}; d_1, d_2 \in L \right\} \subseteq V,$$

$$P_3 = \left\{ \sum_{i=0}^{\infty} a_i x^i \,\middle|\, a_i = \begin{bmatrix} d_1 \\ 0 \\ d_2 \\ 0 \\ 0 \\ 0 \end{bmatrix} ; d_1, d_2 \in L \right\} \subseteq V,$$

$$P_4 = \left\{ \sum_{i=0}^{\infty} a_i x^i \,\middle|\, a_i = \begin{bmatrix} d_1 \\ 0 \\ 0 \\ d_2 \\ 0 \\ 0 \end{bmatrix} ; d_1, d_2 \in L \right\} \subseteq V$$

and

$$P_5 = \left\{ \sum_{i=0}^{\infty} a_i x^i \,\middle|\, a_i = \begin{bmatrix} d_1 \\ 0 \\ 0 \\ 0 \\ d_2 \\ 0 \end{bmatrix} ; d_1, d_2 \in L \right\} \subseteq V$$

be fuzzy linguistic column matrix subsemigroups of V under 'max' operation.

Clearly $V \subseteq P_1 + P_2 + P_3 + P_4 + P_5$ and

$$P_i \cap P_j = \begin{bmatrix} 0 \\ 0 \\ 0 \\ 0 \\ 0 \\ 0 \\ 0 \end{bmatrix}, \text{ if } i \neq j, 1 \leq i, j \leq 5.$$

Thus V is a pseudo direct sum of subsemigroups. Further no $P_i$ is an ideal of V, i = 1, 2, …, 5.

*Example 3.87:* Let

$$S = \left\{ \sum_{i=0}^{9} a_i x^i \;\middle|\; a_i = \begin{bmatrix} d_1 \\ d_2 \\ d_3 \end{bmatrix}; d_1, d_2, d_3 \in L \right\}$$

be a fuzzy linguistic column matrix coefficient polynomial semigorup under max operation. S has no zero divisors.

**THEOREM 3.10:** *Let*

$$S = \left\{ \sum_{i=0}^{\infty} a_i x^i \;\middle|\; a_i = \begin{bmatrix} d_1 \\ d_2 \\ \vdots \\ d_m \end{bmatrix} \;(or\; e_1, e_2, …, e_n), \right.$$

$$e_j \in L,\; 1 \leq j \leq n),\; d_j \in L,\; 1 \leq j \leq m\}$$

*be a fuzzy linguistic polynomial semigroup under 'max'. S has no ideal and no zero divisors.*

The proof is left as an exercise to the reader.

Now we proceed onto define the notion of fuzzy linguistic m × n (m ≠ n) matrix coefficient polynomial semigroup under max (or min).

We give examples of them.

*Example 3.88:* Let

$$V = \left\{ \sum_{i=0}^{7} a_i x^i \,\middle|\, a_i = \begin{bmatrix} d_{11} & d_{12} & \cdots & d_{1n} \\ d_{21} & d_{22} & \cdots & d_{2n} \\ \vdots & \vdots & & \vdots \\ d_{m1} & d_{m2} & \cdots & d_{mn} \end{bmatrix}; d_{ij} \in L, \right.$$

$$1 \leq i \leq m, 1 \leq j \leq n \}$$

be a fuzzy linguistic $m \times n$ ($m \neq n$) coefficient polynomial semigroup under min. V is a semigroup with ideals and zero divisors.

*Example 3.89:* Let

$$V = \left\{ \sum_{i=0}^{8} a_i x^i \,\middle|\, a_i = \begin{bmatrix} d_1 & d_2 & d_3 \\ d_4 & d_5 & d_6 \\ d_7 & d_8 & d_9 \\ d_{10} & d_{11} & d_{12} \\ d_{13} & d_{14} & d_{15} \end{bmatrix}; d_j \in L, 1 \leq j \leq 15 \right\}$$

be a fuzzy linguistic $5 \times 3$ matrix coefficient polynomial semigroup under max. V has no ideals and zero divisors.

Consider

$$P_1 = \left\{ \sum_{i=0}^{8} a_i x^i \,\middle|\, a_i = \begin{bmatrix} d_1 & d_2 & d_3 \\ 0 & 0 & 0 \\ 0 & 0 & 0 \\ 0 & 0 & 0 \\ d_4 & d_5 & d_6 \end{bmatrix}; d_j \in L, 1 \leq j \leq 6 \right\} \subseteq V,$$

$$P_2 = \left\{ \sum_{i=0}^{8} a_i x^i \,\middle|\, a_i = \begin{bmatrix} 0 & 0 & 0 \\ d_1 & d_2 & d_3 \\ 0 & 0 & 0 \\ d_4 & d_5 & d_6 \\ 0 & 0 & 0 \end{bmatrix} ; d_j \in L, 1 \leq j \leq 6 \right\} \subseteq V$$

and

$$P_3 = \left\{ \sum_{i=0}^{8} a_i x^i \,\middle|\, a_i = \begin{bmatrix} 0 & 0 & 0 \\ 0 & 0 & 0 \\ d_1 & d_2 & d_3 \\ 0 & 0 & 0 \\ 0 & 0 & 0 \end{bmatrix} ; d_j \in L, 1 \leq j \leq 3 \right\} \subseteq V ;$$

be a fuzzy linguistic $5 \times 3$ matrix coefficient polynomial semigroups of V.

Clearly $V = P_1 + P_2 + P_3$ with

$$P_i \cap P_j = \begin{bmatrix} 0 & 0 & 0 \\ 0 & 0 & 0 \\ 0 & 0 & 0 \\ 0 & 0 & 0 \\ 0 & 0 & 0 \end{bmatrix}, \text{ if } i \neq j, 1 \leq i, j \leq 3.$$

Further none of $P_1$, $P_2$ or $P_3$ are ideals of V. V has no zero divisors.

*Example 3.90:* Let

$$V = \left\{ \sum_{i=0}^{19} a_i x^i \;\middle|\; a_i = \begin{bmatrix} d_1 & d_2 & d_3 & d_4 & d_5 & d_6 & d_7 \\ d_8 & d_9 & d_{10} & d_{11} & d_{12} & d_{13} & d_{14} \\ d_{15} & d_{16} & d_{17} & d_{18} & d_{19} & d_{20} & d_{21} \end{bmatrix}; \; d_j \in L, \; 1 \leq j \leq 21 \right\}$$

be a fuzzy linguistic $3 \times 7$ matrix coefficient semigroup under 'min'.

Clearly V has ideals and zero divisors.

Consider

$$H_1 = \left\{ \sum_{i=0}^{19} a_i x^i \;\middle|\; a_i = \begin{bmatrix} d_1 & d_4 & 0 & 0 & 0 & 0 & 0 \\ d_2 & d_5 & 0 & 0 & 0 & 0 & 0 \\ d_3 & d_6 & 0 & 0 & 0 & 0 & 0 \end{bmatrix}, \right.$$

$$\left. d_j \in L, \; 1 \leq j \leq 6 \right\} \subseteq V,$$

$$H_2 = \left\{ \sum_{i=0}^{19} a_i x^i \;\middle|\; a_i = \begin{bmatrix} 0 & 0 & d_1 & d_4 & 0 & 0 & 0 \\ 0 & 0 & d_2 & d_5 & 0 & 0 & 0 \\ 0 & 0 & d_3 & d_6 & 0 & 0 & 0 \end{bmatrix}, \right.$$

$$\left. d_j \in L, \; 1 \leq j \leq 6 \right\} \subseteq V,$$

$$H_3 = \left\{ \sum_{i=0}^{19} a_i x^i \;\middle|\; a_i = \begin{bmatrix} 0 & 0 & 0 & 0 & d_1 & d_4 & 0 \\ 0 & 0 & 0 & 0 & d_2 & d_5 & 0 \\ 0 & 0 & 0 & 0 & d_3 & d_6 & 0 \end{bmatrix}, \right.$$

$$\left. d_j \in L, \; 1 \leq j \leq 6 \right\} \subseteq V$$

and

$$H_4 = \left\{ \sum_{i=0}^{19} a_i x^i \;\middle|\; a_i = \begin{bmatrix} 0 & 0 & 0 & 0 & 0 & 0 & d_1 \\ 0 & 0 & 0 & 0 & 0 & 0 & d_2 \\ 0 & 0 & 0 & 0 & 0 & 0 & d_3 \end{bmatrix}, \right.$$

$$\left. d_j \in L,\, 1 \leq j \leq 3 \right\} \subseteq V$$

be fuzzy linguistic $3 \times 7$ matrix coefficient subsemigroup of V under min operation.

We see

$$H_i \cap H_j = \begin{bmatrix} 0 & 0 & 0 & 0 & 0 & 0 & 0 \\ 0 & 0 & 0 & 0 & 0 & 0 & 0 \\ 0 & 0 & 0 & 0 & 0 & 0 & 0 \end{bmatrix}, \text{ if } i \neq j,\, 1 \leq i, j \leq 4.$$

Further $V = H_1 + H_2 + H_3 + H_4$; thus V is direct sum of fuzzy linguistic $3 \times 7$ matrix coefficient subsemigroups of V. Also each $H_i$ is an ideal of V. We see for every $p(x) \in H_i$ and every $q(x) \in H_j$ ($j \neq i$, $1 \leq i, j \leq 4$).

$$\min(p(x), q(x)) = \begin{bmatrix} 0 & 0 & 0 & 0 & 0 & 0 & 0 \\ 0 & 0 & 0 & 0 & 0 & 0 & 0 \\ 0 & 0 & 0 & 0 & 0 & 0 & 0 \end{bmatrix}.$$

Thus V has zero divisors.

Now consider

$$P_1 = \left\{ \sum_{i=0}^{19} a_i x^i \;\middle|\; a_i = \begin{bmatrix} d_1 & 0 & 0 & 0 & 0 & 0 & d_4 \\ d_2 & 0 & 0 & 0 & 0 & 0 & 0 \\ d_3 & 0 & 0 & 0 & 0 & 0 & d_5 \end{bmatrix}, \right.$$

$$\left. d_j \in L,\, 1 \leq j \leq 5 \right\} \subseteq V,$$

$$P_2 = \left\{ \sum_{i=0}^{19} a_i x^i \,\middle|\, a_i = \begin{bmatrix} 0 & d_1 & 0 & 0 & 0 & 0 & d_4 \\ 0 & d_2 & 0 & 0 & 0 & 0 & d_5 \\ 0 & d_3 & 0 & 0 & 0 & 0 & 0 \end{bmatrix} \right.;$$

$$\left. d_j \in L, 1 \leq j \leq 5 \right\} \subseteq V,$$

$$P_3 = \left\{ \sum_{i=0}^{19} a_i x^i \,\middle|\, a_i = \begin{bmatrix} 0 & 0 & d_1 & 0 & 0 & 0 & 0 \\ 0 & 0 & d_2 & 0 & 0 & 0 & d_4 \\ 0 & 0 & d_3 & 0 & 0 & 0 & d_5 \end{bmatrix} \right.;$$

$$\left. d_j \in L, 1 \leq j \leq 5 \right\} \subseteq V,$$

$$P_4 = \left\{ \sum_{i=0}^{19} a_i x^i \,\middle|\, a_i = \begin{bmatrix} 0 & 0 & 0 & d_1 & 0 & 0 & d_4 \\ 0 & 0 & 0 & d_2 & 0 & 0 & d_5 \\ 0 & 0 & 0 & d_3 & 0 & 0 & 0 \end{bmatrix} \right.;$$

$$\left. d_j \in L, 1 \leq j \leq 5 \right\} \subseteq V,$$

$$P_5 = \left\{ \sum_{i=0}^{19} a_i x^i \,\middle|\, a_i = \begin{bmatrix} 0 & 0 & 0 & 0 & d_1 & 0 & d_4 \\ 0 & 0 & 0 & 0 & d_2 & 0 & d_5 \\ 0 & 0 & 0 & 0 & d_3 & 0 & 0 \end{bmatrix} \right.;$$

$$\left. d_j \in L, 1 \leq j \leq 5 \right\} \subseteq V$$

and

$$P_6 = \left\{ \sum_{i=0}^{19} a_i x^i \,\middle|\, a_i = \begin{bmatrix} 0 & 0 & 0 & 0 & 0 & d_1 & d_4 \\ 0 & 0 & 0 & 0 & 0 & d_2 & 0 \\ 0 & 0 & 0 & 0 & 0 & d_3 & d_5 \end{bmatrix} \right.;$$

$$\left. d_j \in L, 1 \leq j \leq 5 \right\} \subseteq V$$

be fuzzy linguistic $3 \times 7$ matrix coefficient polynomial subsemigroups of the semigroup V.

Clearly

$$P_i \cap P_j \neq \begin{bmatrix} 0 & 0 & 0 & 0 & 0 & 0 & 0 \\ 0 & 0 & 0 & 0 & 0 & 0 & 0 \\ 0 & 0 & 0 & 0 & 0 & 0 & 0 \end{bmatrix}, \text{ if } i \neq j, 1 \leq i, j \leq 6,$$

thus V is only a pseudo direct sum of subsemigroups as $V \subseteq P_1 + P_2 + P_3 + P_4 + P_5 + P_6$.

It is easily verified that $P_i$'s are ideals of V, $1 \leq i \leq 6$.

*Example 3.91:* Let

$$P = \left\{ \sum_{i=0}^{2} a_i x^i \;\middle|\; a_i = \begin{bmatrix} d_1 & d_2 \\ d_3 & d_4 \end{bmatrix}, d_j \in L, 1 \leq j \leq 4 \right\}$$

be a fuzzy linguistic square matrix coefficient polynomial semigroup under max.

Take

$$W_1 = \left\{ \sum_{i=0}^{2} a_i x^i \;\middle|\; a_i = \begin{bmatrix} d_1 & 0 \\ 0 & 0 \end{bmatrix}, d_1 \in L \right\} \subseteq P,$$

$$W_2 = \left\{ \sum_{i=0}^{2} a_i x^i \;\middle|\; a_i = \begin{bmatrix} 0 & d_1 \\ 0 & 0 \end{bmatrix}, d_1 \in L \right\} \subseteq P,$$

$$W_3 = \left\{ \sum_{i=0}^{2} a_i x^i \;\middle|\; a_i = \begin{bmatrix} 0 & 0 \\ d_1 & 0 \end{bmatrix}, d_1 \in L \right\} \subseteq P$$

and

$$W_4 = \left\{ \sum_{i=0}^{2} a_i x^i \;\middle|\; a_i = \begin{bmatrix} 0 & 0 \\ 0 & d_1 \end{bmatrix}, d_1 \in L \right\} \subseteq P$$

as fuzzy linguistic matrix coefficient polynomial subsemigroups of P under max operation. Clearly non of $W_j$'s are ideals of P.

Further $W_i \cap W_j = \begin{bmatrix} 0 & 0 \\ 0 & 0 \end{bmatrix}$, $i \neq j$, $1 \leq i, j \leq 4$. However P has no zero divisors.

$P = W_1 + W_2 + W_3 + W_4$ is a direct sum of subsemigroups of P.

Consider

$$S_1 = \left\{ \sum_{i=0}^{2} a_i x^i \;\middle|\; a_i = \begin{bmatrix} d_1 & 0 \\ d_2 & 0 \end{bmatrix}, d_1, d_2 \in L \right\} \subseteq P,$$

$$S_2 = \left\{ \sum_{i=0}^{2} a_i x^i \;\middle|\; a_i = \begin{bmatrix} d_1 & 0 \\ 0 & d_2 \end{bmatrix}, d_1, d_2 \in L \right\} \subseteq P,$$

$$S_3 = \left\{ \sum_{i=0}^{2} a_i x^i \;\middle|\; a_i = \begin{bmatrix} d_1 & d_3 \\ 0 & d_2 \end{bmatrix}, d_1, d_2, d_3 \in L \right\} \subseteq P,$$

$$S_4 = \left\{ \sum_{i=0}^{2} a_i x^i \;\middle|\; a_i = \begin{bmatrix} d_1 & 0 \\ d_2 & d_3 \end{bmatrix}, d_1, d_2, d_3 \in L \right\} \subseteq P$$

and

$$S_5 = \left\{ \sum_{i=0}^{2} a_i x^i \;\middle|\; a_i = \begin{bmatrix} 0 & d_1 \\ d_2 & d_3 \end{bmatrix}, d_1, d_2, d_3 \in L \right\} \subseteq P$$

be fuzzy linguistic square matrix coefficient polynomial subsemigroups under max of P.

Further $S_i \cap S_j \neq \begin{bmatrix} 0 & 0 \\ 0 & 0 \end{bmatrix}$ if $i \neq j$, $1 \leq i, j \leq 5$ and $P \subseteq S_1 + S_2 + S_3 + S_4 + S_5$; hence P is the pseudo direct sum of polynomial subsemigroups.

None of them are ideals of P.

*Example 3.92:* Let

$$M = \left\{ \sum_{i=0}^{5} a_i x^i \;\middle|\; a_i = \begin{bmatrix} d_1 & d_2 & d_3 & d_4 \\ d_5 & d_6 & d_7 & d_8 \\ d_9 & d_{10} & d_{11} & d_{12} \\ d_{13} & d_{14} & d_{15} & d_{16} \end{bmatrix}; d_j \in L, 1 \leq j \leq 16 \right\}$$

be a fuzzy linguistic square matrix coefficient polynomial semigroup under min.

Take

$$P_1 = \left\{ \sum_{i=0}^{5} a_i x^i \;\middle|\; a_i = \begin{bmatrix} d_1 & 0 & 0 & 0 \\ d_2 & d_3 & 0 & 0 \\ 0 & 0 & 0 & 0 \\ 0 & 0 & 0 & 0 \end{bmatrix}; d_1, d_2, d_3 \in L \right\} \subseteq M,$$

$$P_2 = \left\{ \sum_{i=0}^{5} a_i x^i \;\middle|\; a_i = \begin{bmatrix} 0 & d_1 & d_2 & d_3 \\ 0 & 0 & 0 & 0 \\ 0 & 0 & 0 & 0 \\ 0 & 0 & 0 & 0 \end{bmatrix}; d_1, d_2, d_3 \in L \right\} \subseteq M,$$

$$P_3 = \left\{ \sum_{i=0}^{5} a_i x^i \;\middle|\; a_i = \begin{bmatrix} 0 & 0 & 0 & 0 \\ 0 & 0 & d_1 & d_2 \\ 0 & 0 & 0 & 0 \\ 0 & 0 & 0 & 0 \end{bmatrix}; d_1, d_2 \in L \right\} \subseteq M,$$

$$P_4 = \left\{ \sum_{i=0}^{5} a_i x^i \;\middle|\; a_i = \begin{bmatrix} 0 & 0 & 0 & 0 \\ 0 & 0 & 0 & 0 \\ d_1 & d_2 & 0 & 0 \\ d_3 & 0 & 0 & 0 \end{bmatrix}; d_1, d_2, d_3 \in L \right\} \subseteq M,$$

$$P_5 = \left\{ \sum_{i=0}^{5} a_i x^i \;\middle|\; a_i = \begin{bmatrix} 0 & 0 & 0 & 0 \\ 0 & 0 & 0 & 0 \\ 0 & 0 & 0 & 0 \\ 0 & d_1 & d_2 & d_3 \end{bmatrix}; d_1, d_2, d_3 \in L \right\} \subseteq M$$

and

$$P_6 = \left\{ \sum_{i=0}^{5} a_i x^i \;\middle|\; a_i = \begin{bmatrix} 0 & 0 & 0 & 0 \\ 0 & 0 & 0 & 0 \\ 0 & 0 & d_1 & d_2 \\ 0 & 0 & 0 & 0 \end{bmatrix}; d_1, d_2 \in L \right\} \subseteq M$$

be fuzzy linguistic square matrix coefficient polynomial subsemigroups of M.

Clearly $P_i \cap P_j = \begin{bmatrix} 0 & 0 & 0 & 0 \\ 0 & 0 & 0 & 0 \\ 0 & 0 & 0 & 0 \\ 0 & 0 & 0 & 0 \end{bmatrix}$ if $i \neq j$, $1 \leq i, j \leq 6$.

Thus M is a direct sum of subsemigroups under min operation on M.

All these subsemigroups are also ideals of M under min operation.

Further for every $p(x) \in P_i$, every $q(x) \in P_j$, $j \neq i$, $1 \leq i, j \leq 6$ are such that

$$\min(p(x), q(x)) = \begin{bmatrix} 0 & 0 & 0 & 0 \\ 0 & 0 & 0 & 0 \\ 0 & 0 & 0 & 0 \\ 0 & 0 & 0 & 0 \end{bmatrix}.$$

Thus M has zero divisors. If on M the min binary operation is replaced by max then $P_i$'s are subsemigroups of M, $1 \leq i \leq 6$ but never ideals of M under max. Further M has no zero divisors under max operation.

However under max (or min) M is a direct sum of subsemigroups $P_1, P_2, \ldots, P_6$.

Now we proceed onto define set fuzzy linguistic vector spaces using these fuzzy linguistic matrix coefficient polynomials.

Let $V = \left\{ \sum_{i=0}^{\infty} a_i x^i \;\middle|\; a_i = \text{a row matrix or a column matrix or a} \right.$ square matrix or a rectangular matrix (or used in the mutually exclusive sense) with entries from L$\}$ V is a set fuzzy linguistic vector space over L.

We will illustrate this situation by some examples.

*Example 3.93:* Let

$$V = \left\{ \sum_{i=0}^{2} a_i x^i, \sum_{i=0}^{5} b_i x^i, \sum_{i=3}^{7} d_i x^i \;\middle|\; a_i = (p_1, p_2, p_3) \; b_i = (t_1, t_2, \ldots, t_{10}) \right.$$

and $d_i = (r_1, r_2, \ldots, r_6)$ where $p_j, t_k, r_m \in L$, $1 \leq j \leq 3$, $1 \leq k \leq 10$ and $1 \leq m \leq 6$} be a set fuzzy linguistic row matrix coefficient polynomial vector space over the set L.

*Example 3.94:* Let

$$V = \left\{ \sum_{i=0}^{3} a_i x^i, \sum_{i=2}^{9} b_i x^i, \sum_{i=0}^{8} c_i x^i \,\middle|\, a_i = \begin{bmatrix} d_1 \\ d_2 \\ \vdots \\ d_{20} \end{bmatrix}, \; b_i = (t_1, t_2, \ldots, t_{20}) \right.$$

$$\text{and } c_i = \begin{bmatrix} m_1 & m_2 \\ m_3 & m_4 \\ \vdots & \vdots \\ m_{19} & m_{20} \end{bmatrix} \text{ with } d_j, t_k, m_s \in L; \; 1 \leq j \leq 20,$$

$$1 \leq k \leq 20 \text{ and } 1 \leq s \leq 20\}$$

be a set fuzzy linguistic matrix coefficient polynomial vector space over the set M.

*Example 3.95:* Let

$$W = \left\{ \sum_{i=0}^{19} a_i x^i, \sum_{i=10}^{14} b_i x^i, \sum_{i=0}^{7} d_i x^i, \sum_{i=2}^{5} m_i x^i \right.$$

$$a_i = \begin{bmatrix} b_1 & b_2 & \ldots & b_{10} \\ b_{11} & b_{12} & \ldots & b_{20} \\ b_{21} & b_{22} & \ldots & b_{30} \\ b_{31} & b_{32} & \ldots & b_{40} \end{bmatrix} \quad b_i = \begin{bmatrix} y_1 \\ y_2 \\ \vdots \\ y_{12} \end{bmatrix} \quad d_i = \begin{bmatrix} c_1 & c_2 & \ldots & c_{14} \\ c_{15} & c_{16} & \ldots & c_{28} \end{bmatrix}$$

$$m_i = \begin{bmatrix} a & b \\ c & d \end{bmatrix} \text{ where } a, b, c, d, b_j, y_k, c_l \in L;$$

$$1 \leq j \leq 40, \; 1 \leq k \leq 12 \text{ and } 1 \leq l \leq 28\}$$

be a set fuzzy linguistic matrix coefficient polynomial vector space over the fuzzy linguistic set L.

*Example 3.96:* Let

$$W = \left\{ \sum_{i=0}^{3} a_i x^i, \sum_{i=10}^{5} d_i x^i, \sum_{i=0}^{20} b_i x^i, \sum_{i=2}^{4} c_i x^i \right.$$

$$a_i = \begin{bmatrix} p_1 & b_2 & \ldots & b_{20} \\ p_{21} & b_{22} & \ldots & b_{40} \\ p_{41} & b_{42} & \ldots & b_{60} \end{bmatrix}, \; d_i = \begin{bmatrix} m_1 \\ m_2 \\ \vdots \\ m_{40} \end{bmatrix},$$

$$b_i = \begin{bmatrix} n_1 & n_2 & \ldots & n_8 \\ n_9 & n_{10} & \ldots & n_{16} \\ \vdots & \vdots & & \vdots \\ n_{57} & n_{58} & \ldots & n_{64} \end{bmatrix}, \; c_i = (t_1, t_2, \ldots, t_{40})$$

where $p_j, n_k, m_t, t_p \in L$; $1 \le j \le 60$,

$1 \le k \le 64, \; 1 \le t \le 40, \; 1 \le p \le 40 \}$

be a set fuzzy linguistic matrix coefficient polynomial vector space over the fuzzy linguistic set L.

Let

$$P_1 = \left\{ \sum_{i=0}^{3} a_i x^i \;\middle|\; a_i = \begin{bmatrix} p_1 & b_2 & \ldots & b_{20} \\ p_{21} & b_{22} & \ldots & b_{40} \\ p_{41} & b_{42} & \ldots & b_{60} \end{bmatrix}, p_j \in L, 1 \le j \le 60 \right\} \subseteq M,$$

$$P_2 = \left\{ \sum_{i=10}^{5} d_i x^i \;\middle|\; d_i = \begin{bmatrix} m_1 \\ m_2 \\ \vdots \\ m_{40} \end{bmatrix}, m_j \in L, 1 \le j \le 40 \right\} \subseteq M,$$

$$P_3 = \left\{ \sum_{i=0}^{20} b_i x^i \;\middle|\; b_i = \begin{bmatrix} n_1 & n_2 & \cdots & n_8 \\ n_9 & n_{10} & \cdots & n_{16} \\ \vdots & \vdots & & \vdots \\ n_{57} & n_{58} & \cdots & n_{64} \end{bmatrix} n_j \in L,\ 1 \le j \le 64 \right\} \subseteq M$$

and

$$P_4 = \left\{ \sum_{i=2}^{4} c_i x^i \;\middle|\; c_i = (t_1, t_2, \ldots, t_{40});\ t_j \in L,\ 1 \le j \le 40 \right\} \subseteq M$$

be set fuzzy linguistic vector subspaces of matrix coefficient polynomials of M.

Clearly $P_i \cap P_j \neq \phi$ if $i \neq j$, $1 \le i, j \le 4$, also $M = P_1 + P_2 + P_3 + P_4$, so M is a direct sum of subspaces.

*Example 3.97:* Let

$$S = \left\{ \sum_{i=0}^{3} a_i x^i, \sum_{i=4}^{12} b_i x^i, \sum_{i=10}^{17} c_i x^i, \sum_{i=0}^{20} d_i x^i, \sum_{i=1}^{9} f_i x^i, \sum_{i=0}^{11} e_i x^i \right.$$

$$a_i = (m_1, m_2, \ldots, m_{10}),\ b_i = \begin{bmatrix} t_1 \\ t_2 \\ \vdots \\ t_{12} \end{bmatrix},\ c_i = \begin{bmatrix} n_1 & n_2 & \cdots & n_{10} \\ n_{11} & n_{12} & \cdots & n_{20} \\ \vdots & \vdots & & \vdots \\ n_{91} & n_{92} & \cdots & n_{100} \end{bmatrix},$$

$$d_i = \begin{bmatrix} s_1 & s_2 & \cdots & s_{12} \\ s_{13} & s_{14} & \cdots & s_{24} \\ \vdots & \vdots & & \vdots \\ s_{49} & s_{50} & \cdots & s_{60} \end{bmatrix},\ f_i = \begin{bmatrix} y_1 & y_2 \\ y_3 & y_4 \\ \vdots & \vdots \\ y_{29} & y_{30} \end{bmatrix}$$

$$\text{and } e_i = \begin{bmatrix} g_1 & g_2 & g_3 \\ g_4 & g_5 & g_6 \\ \vdots & \vdots & \vdots \\ g_{28} & g_{29} & g_{30} \end{bmatrix} \text{ where } m_j, t_k, n_s, s_t, y_r \text{ and } g_n \in L,$$

$$1 \leq j \leq 100,\ 1 \leq k \leq 12,\ 1 \leq s \leq 100,$$

$$1 \leq t \leq 60,\ 1 \leq r \leq 30 \text{ and } 1 \leq n \leq 30\}$$

be a set fuzzy linguistic matrix coefficient polynomial vector space over the fuzzy linguistic set L.

Consider

$$P_1 = \left\{ \sum_{i=0}^{3} a_i x^i, \sum_{i=4}^{12} b_i x^i \,\bigg|\, a_i = (m_1, m_2, \ldots, m_{10}), \right.$$

$$b_i = \begin{bmatrix} t_1 \\ t_2 \\ \vdots \\ t_{12} \end{bmatrix} \text{ with } m_j, t_k \in L,\ 1 \leq j \leq 10 \text{ and } 1 \leq k \leq 12 \right\} \subseteq S,$$

$$P_2 = \left\{ \sum_{i=0}^{3} a_i x^i, \sum_{i=10}^{17} c_i x^i \,\bigg|\, a_i = (m_1, m_2, \ldots, m_{10}), \right.$$

$$c_i = \begin{bmatrix} n_1 & n_2 & \ldots & n_{10} \\ n_{11} & n_{12} & \ldots & n_{20} \\ \vdots & \vdots & & \vdots \\ n_{91} & n_{92} & \ldots & n_{100} \end{bmatrix} \text{ with } m_j, n_k \in L,\ 1 \leq j \leq 10$$

$$\text{and } 1 \leq k \leq 100 \right\} \subseteq S,$$

$$P_3 = \left\{ \sum_{i=0}^{3} a_i x^i, \sum_{i=0}^{20} d_i x^i \;\middle|\; a_i = (m_1, m_2, \ldots, m_{10}), \right.$$

$$d_i = \begin{bmatrix} s_1 & s_2 & \ldots & s_{12} \\ s_{13} & s_{14} & \ldots & s_{24} \\ \vdots & \vdots & & \vdots \\ s_{49} & s_{50} & \ldots & s_{60} \end{bmatrix} \text{ with } m_j, s_k \in L,$$

$$1 \leq j \leq 10 \text{ and } 1 \leq k \leq 60 \} \subseteq S,$$

$$P_4 = \left\{ \sum_{i=0}^{3} a_i x^i, \sum_{i=1}^{9} f_i x^i \;\middle|\; a_i = (m_1, m_2, \ldots, m_{10}), \right.$$

$$f_i = \begin{bmatrix} y_1 & y_2 \\ y_3 & y_4 \\ \vdots & \vdots \\ y_{29} & y_{30} \end{bmatrix} \text{ with } m_j, y_k \in L,$$

$$1 \leq j \leq 10 \text{ and } 1 \leq k \leq 30 \} \subseteq S$$

and

$$P_5 = \left\{ \sum_{i=0}^{3} a_i x^i, \sum_{i=0}^{11} e_i x^i \;\middle|\; a_i = (m_1, m_2, \ldots, m_{10}), \right.$$

$$e_i = \begin{bmatrix} g_1 & g_2 & g_3 \\ g_4 & g_5 & g_6 \\ \vdots & \vdots & \vdots \\ g_{28} & g_{29} & g_{30} \end{bmatrix} \text{ with } m_j, g_k \in L,$$

$$1 \leq j \leq 10 \text{ and } 1 \leq k \leq 30 \} \subseteq S$$

be fuzzy linguistic matrix coefficient polynomial vector subspace over the fuzzy linguistic set L.

Clearly $P_i \cap P_j = \left\{ \sum_{i=0}^{3} a_i x^i \,\middle|\, a_i = (m_1, m_2, \ldots, m_{10}), m_j \in L, 1 \leq j \leq 10 \right\} \subseteq S$ if $i \neq j$, $1 \leq i, j \leq 5$} thus $S \subseteq P_1 + P_2 + P_3 + P_4 + P_5$ so S is a pseudo direct sum of vector subspaces of S.

Now we proceed onto give the examples of set fuzzy linguistic linear algebra over a fuzzy linguistic set L. We call a set fuzzy linguistic set vector space V to be a fuzzy linguistic linear algebra if on V we have either the min (or max) closed binary operation which is defined on V.

*Example 3.98:* Let

$$S = \left\{ \sum_{i=0}^{20} a_i x^i \,\middle|\, a_i = \begin{bmatrix} d_1 & d_2 & \ldots & d_{10} \\ d_{11} & d_{12} & \ldots & d_{20} \\ \vdots & \vdots & & \vdots \\ d_{81} & d_{82} & \ldots & d_{90} \end{bmatrix} \text{ where } d_j \in L, 1 \leq j \leq 90 \right\}$$

be a set fuzzy linguistic matrix coefficient polynomial linear algebra over the fuzzy linguistic set L. S under min operation is a fuzzy linguistic polynomial semigroup.

Consider

$$W = \left\{ \sum_{i=0}^{20} a_i x^i \,\middle|\, a_i = \begin{bmatrix} d_1 & d_2 & \ldots & d_{10} \\ 0 & 0 & \ldots & 0 \\ \vdots & \vdots & & \vdots \\ d_{11} & d_{12} & \ldots & d_{20} \end{bmatrix} \text{ where } d_j \in L, \right.$$

$$\left. 1 \leq j \leq 20 \right\} \subseteq S,$$

W is a set fuzzy linguistic linear subalgebra of S over L.

*Example 3.99:* Let

$$S = \left\{ \sum_{i=0}^{20} a_i x^i \,\middle|\, a_i = \begin{bmatrix} d_1 & d_2 & d_3 & d_4 \\ d_5 & d_6 & d_7 & d_8 \\ \vdots & \vdots & & \vdots \\ d_{57} & d_{58} & d_{59} & d_{60} \end{bmatrix} \text{ where } d_j \in L, 1 \le j \le 60 \right\}$$

be a set fuzzy linguistic matrix coefficient polynomial linear algebra over the fuzzy linguistic set L under max operation.

Consider

$$P = \left\{ \sum_{i=0}^{20} a_i x^i \,\middle|\, a_i = \begin{bmatrix} d_1 & d_2 & d_3 & d_4 \\ d_5 & d_6 & d_7 & d_8 \\ d_9 & d_{10} & d_{11} & d_{12} \\ 0 & 0 & 0 & 0 \\ \vdots & \vdots & \vdots & \vdots \\ 0 & 0 & 0 & 0 \end{bmatrix} \text{ where } d_j \in L, \right.$$

$$1 \le j \le 12 \Big\} \subseteq M,$$

P is not a set fuzzy linguistic linear subalgebra of M under max operation, only a set fuzzy linguistic vector subspace of M known as the set fuzzy linguistic pseudo vector subspace of M.

*Example 3.100:* Let

$$S = \left\{ \sum_{i=0}^{3} a_i x^i \,\middle|\, a_i = \begin{bmatrix} d_1 \\ d_2 \\ d_3 \\ d_4 \\ d_5 \end{bmatrix} \text{ where } d_j \in L,\ 1 \le j \le 5 \right\}$$

be a set fuzzy linguistic column matrix coefficient polynomial linear algebra over the set L under min operation V. V has sublinear algebras and zero divisors. However if min operation on V is replaced by max, results are not true identically.

Next we give examples of semigroup fuzzy linguistic vector spaces and semigroup fuzzy linguistic linear algebras using fuzzy linguistic matrix coefficient polynomials.

*Example 3.101:* Let V be a set of fuzzy linguistic matrix coefficient polynomials with entries from L. L be a semigroup under min. (L the fuzzy linguistic set). V be a semigroup fuzzy linguistic matrix coefficient polynomial vector space over the semigroup.

*Example 3.102:* Let V be a set of fuzzy linguistic matrix coefficient polynomials given by

$$V = \left\{ \sum_{i=0}^{8} a_i x^i, \sum_{i=1}^{3} d_i x^i, \sum_{i=1}^{4} g_i x^i \,\middle|\, a_i = (t_1, t_2, ..., t_9), \right.$$

$$d_i = \begin{bmatrix} m_1 & m_2 \\ m_3 & m_4 \\ \vdots & \vdots \\ m_{19} & m_{20} \end{bmatrix} \text{ and } g_i = \begin{bmatrix} n_1 & n_2 & ... & n_{10} \\ n_{11} & n_{12} & ... & n_{20} \\ \vdots & \vdots & & \vdots \\ n_{51} & n_{52} & ... & n_{60} \end{bmatrix}$$

with $t_j, m_k, n_t \in L$, $1 \leq j \leq 9$, $1 \leq k \leq 20$ and $1 \leq t \leq 60\}$,

be a semigroup fuzzy linguistic vector space over the semigroup L under min operation.

*Example 3.103:* Let

$$P = \left\{ \sum_{i=0}^{7} a_i x^i, \sum_{i=8}^{10} d_i x^i, \sum_{i=20}^{5} g_i x^i \,\middle|\, a_i = \begin{bmatrix} m_1 \\ m_2 \\ \vdots \\ m_{20} \end{bmatrix}, d_i = \begin{bmatrix} e_1 & e_2 & e_3 \\ e_4 & e_5 & e_6 \\ e_7 & e_8 & e_9 \end{bmatrix} \right.$$

$$\text{and } g_i = \begin{bmatrix} n_1 & n_2 & \dots & n_{20} \\ n_{21} & n_{22} & \dots & n_{40} \\ n_{41} & n_{42} & & n_{60} \\ n_{61} & n_{62} & \dots & n_{80} \end{bmatrix}$$

with $e_j, m_k, n_t \in L$, $1 \le j \le 9$, $1 \le k \le 20$ and $1 \le t \le 80\}$

be a semigroup fuzzy linguistic matrix coefficient polynomial vector space over the fuzzy linguistic semigroup L.

*Example 3.104:* Let

$$P = \left\{ \sum_{i=0}^{8} a_i x^i, \sum_{i=1}^{3} g_i x^i, \sum_{i=10}^{19} d_i x^i \,\middle|\, a_i = \begin{bmatrix} t_1 \\ t_2 \\ \vdots \\ t_{19} \end{bmatrix}, \right.$$

$$g_i = \begin{bmatrix} d_1 & d_2 & \dots & d_{20} \\ d_{21} & d_{22} & \dots & d_{40} \\ d_{41} & d_{42} & \dots & d_{60} \end{bmatrix} \text{ and } d_k = \begin{bmatrix} p_1 & p_2 & \dots & p_{10} \\ p_{11} & p_{12} & \dots & p_{20} \\ p_{21} & p_{22} & \dots & p_{30} \\ p_{31} & p_{32} & \dots & p_{40} \end{bmatrix}$$

where $t_j, d_k, p_l \in L$, $1 \le j \le 19$, $1 \le k \le 60$ and $1 \le t \le 40\}$

be a semigroup fuzzy linguistic matrix coefficient vector space over the fuzzy linguistic semigroup L.

Let

$$M_1 = \left\{ \sum_{i=0}^{8} a_i x^i \;\middle|\; a_i = \begin{bmatrix} t_1 \\ t_2 \\ \vdots \\ t_{19} \end{bmatrix}, t_j \in L, 1 \le j \le 19 \right\} \subseteq P,$$

$$M_2 = \left\{ \sum_{i=1}^{3} g_i x^i \;\middle|\; g_i = \begin{bmatrix} d_1 & d_2 & \ldots & d_{20} \\ d_{21} & d_{22} & \ldots & d_{40} \\ d_{41} & d_{42} & \ldots & d_{60} \end{bmatrix} d_j \in L, 1 \le j \le 19 \right\} \subseteq P$$

and

$$M_3 = \left\{ \sum_{i=10}^{19} d_i x^i \;\middle|\; g_i = \begin{bmatrix} p_1 & p_2 & \cdots & p_{10} \\ p_{11} & p_{12} & \cdots & p_{20} \\ p_{21} & 9_{22} & \cdots & p_{30} \\ p_{31} & p_{32} & \cdots & p_{40} \end{bmatrix} p_j \in L, 1 \le j \le 40 \right\} \subseteq P;$$

be semigroup fuzzy linguistic matrix coefficient vector subspace over the fuzzy linguistic semigroup L under min operation.

Clearly $M_i \cap M_j = \phi$, $1 \le i, j \le 3$, $P = M_1 + M_2 + M_3$. Thus P is the direct sum of semigroup fuzzy linguistic matrix coefficient vector subspaces over the fuzzy linguistic semigroup L.

*Example 3.105:* Let

$$P = \left\{ \sum_{i=0}^{20} a_i x^i \;\middle|\; a_i = \begin{bmatrix} d_1 & d_2 & \ldots & d_{40} \\ d_{41} & d_{42} & \ldots & d_{80} \\ d_{41} & d_{42} & \ldots & d_{120} \\ d_{121} & d_{122} & \ldots & d_{160} \\ d_{161} & d_{162} & \ldots & d_{200} \end{bmatrix} \text{ with } d_j \in L, 1 \le j \le 200 \right\}$$

be a semigroup fuzzy linguistic linear algebra of matrix coefficient polynomials over the fuzzy linguistic semigroup; 'min' on L. P also under min operation only we give sublinear algebra of P on L.

$$M_1 = \left\{ \sum_{i=0}^{20} a_i x^i \,\middle|\, a_i = \begin{bmatrix} 0 & 0 & ... & 0 \\ d_1 & d_2 & ... & d_{40} \\ d_{41} & d_{42} & ... & d_{80} \\ 0 & 0 & ... & 0 \\ 0 & 0 & ... & 0 \end{bmatrix} \text{ with } d_j \in L, \right.$$

$$\left. 1 \le j \le 80 \right\} \subseteq P,$$

$$M_2 = \left\{ \sum_{i=0}^{20} a_i x^i \,\middle|\, a_i = \begin{bmatrix} d_1 & d_2 & ... & d_{40} \\ 0 & 0 & ... & 0 \\ 0 & 0 & ... & 0 \\ 0 & 0 & ... & 0 \\ 0 & 0 & ... & 0 \end{bmatrix} \text{ with } d_j \in L, \right.$$

$$\left. 1 \le j \le 40 \right\} \subseteq P$$

and

$$M_3 = \left\{ \sum_{i=0}^{20} a_i x^i \,\middle|\, a_i = \begin{bmatrix} 0 & 0 & ... & 0 \\ 0 & 0 & ... & 0 \\ 0 & 0 & ... & 0 \\ d_1 & d_2 & ... & d_{40} \\ d_{41} & d_{42} & ... & d_{80} \end{bmatrix} \text{ with } d_j \in L, \right.$$

$$\left. 1 \le j \le 40 \right\} \subseteq P.$$

We see $M_i \cap M_j = \begin{bmatrix} 0 & 0 & ... & 0 \\ 0 & 0 & ... & 0 \\ 0 & 0 & ... & 0 \\ 0 & 0 & ... & 0 \\ 0 & 0 & ... & 0 \end{bmatrix}$ if $i \neq j$, $1 \leq i, j \leq 3$.

Finally $P = M_1 + M_2 + M_3$ is the direct sum of semigroup fuzzy linguistic linear subalgebra.

*Example 3.106:* Let

$$P = \left\{ \sum_{i=0}^{\infty} a_i x^i \middle| a_i = \begin{bmatrix} d_1 & d_2 & d_3 \\ d_4 & d_5 & d_6 \\ \vdots & \vdots & \vdots \\ d_{28} & d_{29} & d_{30} \end{bmatrix}; d_j \in L, 1 \leq j \leq 30 \right\}$$

be a semigroup fuzzy linguistic linear algebra.

Consider

$$M_1 = \left\{ \sum_{i=0}^{\infty} a_i x^i \middle| a_i = \begin{bmatrix} d_1 & d_2 & d_3 \\ d_4 & d_5 & d_6 \\ 0 & 0 & 0 \\ 0 & 0 & 0 \\ \vdots & \vdots & \vdots \\ 0 & 0 & 0 \end{bmatrix}; d_j \in L, 1 \leq j \leq 6 \right\} \subseteq P,$$

$$M_2 = \left\{ \sum_{i=0}^{\infty} a_i x^i \middle| a_i = \begin{bmatrix} d_1 & d_2 & d_3 \\ 0 & 0 & 0 \\ d_4 & d_5 & d_6 \\ 0 & 0 & 0 \\ \vdots & \vdots & \vdots \\ 0 & 0 & 0 \end{bmatrix}; d_j \in L, 1 \leq j \leq 6 \right\} \subseteq P,$$

$$M_3 = \left\{ \sum_{i=0}^{\infty} a_i x^i \;\middle|\; a_i = \begin{bmatrix} d_1 & d_2 & d_3 \\ 0 & 0 & 0 \\ 0 & 0 & 0 \\ d_4 & d_5 & d_6 \\ 0 & 0 & 0 \\ \vdots & \vdots & \vdots \\ 0 & 0 & 0 \end{bmatrix} ; d_j \in L, 1 \leq j \leq 6 \right\} \subseteq P,$$

$$M_4 = \left\{ \sum_{i=0}^{\infty} a_i x^i \;\middle|\; a_i = \begin{bmatrix} d_1 & d_2 & d_3 \\ 0 & 0 & 0 \\ 0 & 0 & 0 \\ 0 & 0 & 0 \\ d_4 & d_5 & d_6 \\ 0 & 0 & 0 \\ \vdots & \vdots & \vdots \\ 0 & 0 & 0 \end{bmatrix} ; d_j \in L, 1 \leq j \leq 6 \right\} \subseteq P,$$

$$M_5 = \left\{ \sum_{i=0}^{\infty} a_i x^i \;\middle|\; a_i = \begin{bmatrix} d_1 & d_2 & d_3 \\ 0 & 0 & 0 \\ 0 & 0 & 0 \\ 0 & 0 & 0 \\ 0 & 0 & 0 \\ d_4 & d_5 & d_6 \\ 0 & 0 & 0 \\ \vdots & \vdots & \vdots \\ 0 & 0 & 0 \end{bmatrix} ; d_j \in L, 1 \leq j \leq 6 \right\} \subseteq P,$$

$$M_6 = \left\{ \sum_{i=0}^{\infty} a_i x^i \;\middle|\; a_i = \begin{bmatrix} d_1 & d_2 & d_3 \\ 0 & 0 & 0 \\ 0 & 0 & 0 \\ 0 & 0 & 0 \\ 0 & 0 & 0 \\ 0 & 0 & 0 \\ d_4 & d_5 & d_6 \\ 0 & 0 & 0 \\ 0 & 0 & 0 \\ 0 & 0 & 0 \end{bmatrix} ; d_j \in L,\ 1 \le j \le 6 \right\} \subseteq P,$$

$$M_7 = \left\{ \sum_{i=0}^{\infty} a_i x^i \;\middle|\; a_i = \begin{bmatrix} d_1 & d_2 & d_3 \\ 0 & 0 & 0 \\ 0 & 0 & 0 \\ 0 & 0 & 0 \\ 0 & 0 & 0 \\ 0 & 0 & 0 \\ 0 & 0 & 0 \\ d_4 & d_5 & d_6 \\ 0 & 0 & 0 \\ 0 & 0 & 0 \end{bmatrix} ; d_j \in L,\ 1 \le j \le 6 \right\} \subseteq P,$$

$$M_8 = \left\{ \sum_{i=0}^{\infty} a_i x^i \,\middle|\, a_i = \begin{bmatrix} d_1 & d_2 & d_3 \\ 0 & 0 & 0 \\ 0 & 0 & 0 \\ 0 & 0 & 0 \\ 0 & 0 & 0 \\ 0 & 0 & 0 \\ 0 & 0 & 0 \\ 0 & 0 & 0 \\ d_4 & d_5 & d_6 \\ 0 & 0 & 0 \end{bmatrix} ; d_j \in L, 1 \leq j \leq 6 \right\} \subseteq P$$

and

$$M_9 = \left\{ \sum_{i=0}^{\infty} a_i x^i \,\middle|\, a_i = \begin{bmatrix} d_1 & d_2 & d_3 \\ 0 & 0 & 0 \\ 0 & 0 & 0 \\ 0 & 0 & 0 \\ 0 & 0 & 0 \\ 0 & 0 & 0 \\ 0 & 0 & 0 \\ 0 & 0 & 0 \\ 0 & 0 & 0 \\ d_4 & d_5 & d_6 \end{bmatrix} ; d_j \in L, 1 \leq j \leq 6 \right\} \subseteq P,$$

$M_1$, $M_2$, $M$ ,…, $M_9$ are semigroup fuzzy linguistic linear subalgebras of P over the fuzzy linguistic semigroup L under 'min'.

$$M_i \cap M_j = \begin{bmatrix} d_1 & d_2 & d_3 \\ 0 & 0 & 0 \\ 0 & 0 & 0 \\ 0 & 0 & 0 \\ 0 & 0 & 0 \\ 0 & 0 & 0 \\ 0 & 0 & 0 \\ 0 & 0 & 0 \\ 0 & 0 & 0 \\ 0 & 0 & 0 \end{bmatrix} \text{ if } i \neq j, 1 \leq i, j \leq 9.$$

Thus P is only a pseudo direct sum of semigroup fuzzy linguistic linear subalgebras of M over the fuzzy linguistic semigroup L.

Now as in case of usual vector spaces we can in case of set (semigroup) fuzzy linguistic matrix coefficient polynomial vector spaces (linear algebras) define subset (subsemigroup) fuzzy linguistic matrix coefficient polynomial vector subspaces (linear subalgebras), set (semigroup) linear transformation, linear operators and linear functionals define, describe and study the corresponding results.

**Chapter Four**

# FUZZY LINGUISTIC MODELS

In this chapter we introduce the new notion of fuzzy linguistic models and apply them to social problems. We need for this the operations on the fuzzy linguistic matrices. We briefly describe the operations on them. Let L denote a fuzzy linguistic space consisting of fuzzy linguistic terms associated with some problem P.

This chapter has four sections. Section one defines operations on fuzzy linguistic matrices. The new notion of fuzzy linguistic cognitive models is introduced in section two. Fuzzy linguistic relational map model is defined and described in section three. In the final section fuzzy linguistic relation equation is introduced.

**4.1 Operations on Fuzzy Linguistic Matrices**

In this section the authors just introduce the operations on the fuzzy linguistic matrices.

Let $x = (a_1, a_2, a_3, a_4, a_5, a_6)$ be a fuzzy linguistic row matrix (vector) with entries from L.

$$x^t = (a_1, a_2, a_3, a_4, a_5, a_6)^t = \begin{bmatrix} a_1 \\ a_2 \\ a_3 \\ a_4 \\ a_5 \\ a_6 \end{bmatrix}$$

is the fuzzy linguistic column matrix (vector). That is if x = (0, bad, good, worst, fair, very bad) be the fuzzy linguistic row matrix, then transpose of x denoted by

$$x^t = (0, \text{bad, good, worst, fair, very bad})^t = \begin{bmatrix} 0 \\ \text{bad} \\ \text{good} \\ \text{worst} \\ \text{fair} \\ \text{very bad} \end{bmatrix}$$

which is a fuzzy linguistic column matrix. Suppose we have product of a fuzzy linguistic row matrix x can be operated with min or max operation with a fuzzy linguistic column matrix y if and only if x is a 1 × n row matrix then y must be a n × 1 matrix. We show how the 'min' operation looks like.

Suppose

$$y = \begin{bmatrix} a_1 \\ a_2 \\ \vdots \\ a_n \end{bmatrix} \text{ and } x = (b_1, b_2, \ldots, b_n)$$

then min {min {x, y}} = min {min {$a_1$, $b_1$}, min {$a_2$, $b_2$}, …, min {$a_n$, $b_n$}}.

Likewise we can have max max operation; max min operation and min max operation. All these will be illustrated by examples.

Let $x = (x_1, x_2, \ldots, x_{25})$ and $y = (y_1, y_2, \ldots, y_{25})$ with $x_i, y_i \in L$; $1 \leq i \leq 25$. L a fuzzy linguistic space.

max {max (x, y)} = max {max $(x_1, y_1)$, max $\{x_2, y_2\}$, …, max $\{x_{25}, y_{25}\}$}.

Let $x = (x_1, x_2, \ldots x_{10})$ and $y = (y_1, y_2, \ldots, y_{10})$ max (min {x, y}) = max {min $\{x_1, y_1\}$, min $\{x_2, y_2\}$, …, min $\{x_{10}, y_{10}\}$}.

Likewise min (max {x, y}) = min {max $\{x_1, y_1\}$, max $\{x_2, y_2\}$, …, max $\{x_{10}, y_{10}\}$}.

We will see how for a typical problem the four operations give different sets of answers.

Let

$$x = \begin{bmatrix} good \\ bad \\ fair \\ 0 \\ v.fair \\ best \\ bad \\ 0 \\ good \\ better \end{bmatrix} \text{ and }$$

y = (bad, 0, good, very bad, 0, best, fair, best, better, 0)

be two fuzzy linguistic matrices.

min {min (y, x)}
= min {bad, 0, fair, 0, 0, best, bad, 0, better, 0} = 0

max {min (y, x)}
= max {bad, 0, fair, 0, 0, best, bad, 0, better, 0} = best

min {max (y, z)}
= min {good, bad, good, very bad, v.fair, best, fair, best, good, better}

= very bad and

max (max {y, x})

= max {good, bad, good, very bad, very fair, best fair, best, good, better}

= best.

Suppose we are interested in finding min {x, y} then we have

min {x, y} =

$$\begin{bmatrix} bad & 0 & good & very\ bad & 0 & good & fair & good & better & 0 \\ bad & 0 & bad & very\ bad & 0 & bad & bad & bad & bad & 0 \\ bad & 0 & fair & very\ bad & 0 & fair & fair & fair & fair & 0 \\ 0 & 0 & 0 & 0 & 0 & 0 & 0 & 0 & 0 & 0 \\ bad & 0 & v.fair & very\ bad & 0 & v.fair & v.fair & v.fair & v.fair & 0 \\ bad & 0 & good & very\ bad & 0 & best & fair & best & better & 0 \\ bad & 0 & bad & very\ bad & 0 & bad & bad & bad & bad & 0 \\ 0 & 0 & 0 & 0 & 0 & 0 & 0 & 0 & 0 & 0 \\ bad & 0 & good & very\ bad & 0 & good & fair & good & better & 0 \\ bad & 0 & better & very\ bad & 0 & better & fair & better & better & 0 \end{bmatrix}.$$

This is the way a fuzzy linguistic column matrix of order n × 1 is multipliced with a 1 × n fuzzy linguistic row matrix with min operation.

We can also find max {x, y} =

$$\begin{bmatrix} \text{good} & \text{bad} & \text{good} & \text{very bad} & \text{v.fair} & \text{best} & \text{fair} & \text{best} & \text{good} & \text{better} \\ \text{bad} & \text{bad} & \text{good} & \text{bad} & \text{bad} & \text{best} & \text{fair} & \text{best} & \text{better} & \text{bad} \\ \text{fair} & \text{fair} & \text{good} & \text{fair} & \text{fair} & \text{best} & \text{fair} & \text{best} & \text{good} & \text{better} \\ \text{bad} & 0 & \text{good} & \text{v.bad} & \text{fair} & \text{best} & \text{fair} & \text{best} & \text{better} & \text{bad} \\ \text{v.fair} & \text{v.fair} & \text{good} & \text{v.fair} & \text{v.fair} & \text{best} & \text{v.fair} & \text{best} & \text{fair} & \text{fair} \\ \text{best} & \text{best} & \text{best} & \text{best} & \text{best} & \text{best} & \text{best} & \text{best} & \text{better} & 0 \\ \text{bad} & \text{bad} & \text{good} & \text{bad} & \text{bad} & \text{best} & \text{fair} & \text{best} & \text{better} & \text{v.fair} \\ \text{bad} & 0 & \text{good} & \text{v.bad} & 0 & \text{best} & \text{fair} & \text{best} & \text{best} & \text{best} \\ \text{good} & \text{good} & \text{good} & \text{good} & \text{good} & \text{best} & \text{good} & \text{best} & \text{good} & \text{good} \\ \text{better} & \text{better} & \text{good} & \text{better} & \text{better} & \text{best} & \text{better} & \text{best} & \text{better} & \text{better} \end{bmatrix}.$$

We see clearly max {x, y} ≠ min {x, y}.

Further 0 dominates in min and best dominates in max. When an expert wants to boost the results on the positive side he can use the max operation. If min operation is used it gives the worst state of affairs.

Thus we see we get different results for these four types of operations. According to need one can use any one of the operations.

Now we find the transpose of a m × n fuzzy linguistic matrix M.

Let M =

$$M = \begin{bmatrix} \text{good} & 0 & \text{best} & \text{bad} & \text{fair} & \text{best} \\ 0 & \text{good} & \text{bad} & \text{worst} & 0 & \text{very bad} \\ \text{bad} & 0 & \text{fair} & \text{very fair} & \text{bad} & \text{better} \\ \text{worst} & \text{bad} & 0 & \text{good} & \text{better} & 0 \\ 0 & \text{best} & \text{fair} & 0 & \text{worst} & \text{best} \\ \text{good} & \text{better} & \text{fair} & \text{best} & \text{worst} & 0 \end{bmatrix}.$$

Now the transpose of this fuzzy linguistic 6 × 6 matrix M denoted

$$M^t = \begin{bmatrix} \text{good} & 0 & \text{bad} & \text{worst} & 0 & \text{good} \\ 0 & \text{good} & 0 & \text{bad} & \text{best} & \text{better} \\ \text{best} & \text{bad} & \text{fair} & 0 & \text{fair} & \text{fair} \\ \text{bad} & \text{worst} & \text{very fair} & \text{good} & 0 & \text{best} \\ \text{fair} & 0 & \text{bad} & \text{better} & \text{worst} & \text{worst} \\ \text{best} & \text{very bad} & \text{better} & 0 & \text{best} & 0 \end{bmatrix}.$$

Likewise if

$$N = \begin{bmatrix} \text{best} & 0 & \text{bad} \\ \text{fair} & \text{bad} & 0 \\ 0 & \text{very fair} & \text{good} \\ \text{good} & 0 & \text{best} \\ \text{best} & \text{worst} & 0 \\ 0 & \text{best} & \text{good} \\ \text{good} & \text{better} & \text{best} \\ 0 & \text{bad} & 0 \end{bmatrix}$$

is any fuzzy linguistic 8 × 3 matrix, to find $N^T$.

$$N^T = \begin{bmatrix} \text{best} & \text{fair} & 0 & \text{good} & \text{best} & 0 & \text{good} & 0 \\ 0 & \text{bad} & \text{very fair} & 0 & \text{worst} & \text{best} & \text{better} & \text{bad} \\ \text{bad} & 0 & \text{good} & \text{best} & 0 & \text{good} & \text{best} & 0 \end{bmatrix}.$$

We see $N^T$ is a $3 \times 8$ fuzzy linguistic matrix.

Now we can find the product of two rectangular linguistic matrices M and N if M is a $n \times t$ matrix then N must be a $t \times m$ matrix then only MN is defined how ever NM is not defined in this case and MN is a $n \times m$ fuzzy linguistic matrix.

However product of any $n \times n$ matrix with itself is always defined.

We will illustrate these two situations by some examples.

Let

$$M = \begin{bmatrix} \text{good} & \text{bad} & 0 & \text{best} & \text{worst} \\ \text{bad} & 0 & \text{good} & 0 & \text{best} \\ 0 & \text{good} & \text{bad} & \text{worst} & 0 \\ \text{fair} & \text{best} & 0 & \text{bad} & \text{good} \\ \text{good} & \text{bad} & \text{good} & 0 & \text{bad} \end{bmatrix}$$

$$= \min(\min(M, M)) = \begin{bmatrix} 0 & 0 & 0 & 0 & 0 \\ 0 & 0 & 0 & 0 & 0 \\ 0 & 0 & 0 & 0 & 0 \\ 0 & 0 & 0 & 0 & 0 \\ 0 & 0 & 0 & 0 & 0 \end{bmatrix}.$$

We see if '0' occurs atleast once in every row and atleast once in every column then $\min(\min(M, M)) = (0)$.

Now we find

$$\max(\min(M, M)) = \begin{bmatrix} good & best & worst & good & good \\ good & good & good & bad & worst \\ bad & worst & good & bad & good \\ good & bad & good & fair & best \\ good & good & bad & good & bad \end{bmatrix}.$$

We now find the value of

$$\max(\min(M, M)) = \begin{bmatrix} best & best & best & best & best \\ best & best & best & best & best \\ good & good & good & good & good \\ best & best & best & best & best \\ good & good & good & good & good \end{bmatrix}.$$

$$\min(\max(M, M)) = \begin{bmatrix} 0 & bad & bad & worst & 0 \\ bad & 0 & bad & 0 & bad \\ bad & bad & 0 & bad & worst \\ 0 & fair & bad & worst & 0 \\ bad & bad & 0 & bad & bad \end{bmatrix}.$$

We see max {max (M, M)} gives an extreme or better values and min {min {M, M}} gives an extreme low values, where as min max non negative values and max min more positive values.

Such four types of operations can be used as per need of the problem.

Now if M is a n × m matrix and N is a m × t matrix we can find min (min {M, N}), min (max {M, N}), max (max (M, N)} and max (min {M, N}).

We will illustrate these four types of operations.

Let

$$M = \begin{bmatrix} \text{low} & 0 & \text{high} & 0 & \text{very low} & \text{very high} \\ 0 & \text{high} & 0 & \text{low} & \text{high} & 0 \\ \text{high} & 0 & \text{medium} & 0 & \text{low} & \text{low} \\ \text{low} & \text{low} & 0 & \text{medium} & 0 & \text{medium} \end{bmatrix}$$

and

$$N = \begin{bmatrix} \text{low} & 0 & \text{medium} & 0 \\ \text{high} & \text{low} & 0 & \text{high} \\ 0 & \text{very high} & \text{low} & 0 \\ \text{medium} & 0 & \text{high} & \text{low} \\ \text{very low} & \text{medium} & 0 & \text{high} \\ \text{high} & \text{very low} & \text{low} & \text{medium} \end{bmatrix}$$

be two fuzzy linguistic matrices associated with temperature of an experiment.

$$\text{Clearly min (min }\{M, N\}) = \begin{bmatrix} 0 & 0 & 0 & 0 \\ 0 & 0 & 0 & 0 \\ 0 & 0 & 0 & 0 \\ 0 & 0 & 0 & 0 \end{bmatrix}.$$

$$\max \{\min (M, N)\} = \begin{bmatrix} \text{high} & \text{high} & \text{low} & \text{medium} \\ \text{high} & \text{medium} & \text{low} & \text{high} \\ \text{low} & \text{medium} & \text{medium} & \text{low} \\ \text{medium} & \text{low} & \text{medium} & \text{medium} \end{bmatrix}$$

$$\max\{\max(M, N)\} = \begin{bmatrix} \text{very high} & \text{very high} & \text{very high} & \text{very high} \\ \text{high} & \text{very high} & \text{high} & \text{high} \\ \text{high} & \text{very high} & \text{high} & \text{high} \\ \text{high} & \text{very high} & \text{high} & \text{high} \end{bmatrix}.$$

$$\text{Now } \min\{\max(M, N)\} = \begin{bmatrix} \text{low} & 0 & 0 & \text{low} \\ 0 & 0 & \text{low} & 0 \\ \text{low} & 0 & \text{low} & \text{low} \\ 0 & \text{low} & 0 & 0 \end{bmatrix}.$$

We see this four types of operations gives four types of 4 × 4 fuzzy linguistic matrices. Now we can also find linguistic row matrix with a square or a rectangular fuzzy linguistic matrix which is compatible.

Consider X = (fast, slow, very slow, just fast, very fast) to be the fuzzy linguistic row matrix.

Take

$$M = \begin{bmatrix} \text{slow} & \text{medium} & \text{fast} & \text{slow} \\ \text{very slow} & \text{slow} & \text{medium} & \text{very slow} \\ \text{fast} & \text{slow} & \text{fast} & \text{just fast} \\ \text{fast} & \text{medium} & \text{slow} & \text{fast} \\ \text{just fast} & \text{slow} & \text{very slow} & \text{fast} \end{bmatrix}$$

be a 5 × 4 fuzzy linguistic matrix. Consider min {min {X, M}} = (very slow, very slow, very slow, very slow).

Consider max (min {X, M}) = (just fast, medium, fast, fast).

Consider min (max {X, M}) = (slow, slow, medium, slow).

To find max (max {X, M}) = (very fast, very fast, very fast, very fast).

Now one can compare the results and use the appropriate operations needed.

Now we can also find other types of operations, which is left as an exercise to the reader.

Now we proceed onto describe the models.

## 4.2 Fuzzy Linguistic Cognitive Models

In this section the notion of Fuzzy Linguistic Cognitive Models (FLCM) are introduced and described.

For reference about Fuzzy Cognitive Maps and its working please refer [3-5, 11].

Let us consider $C_1, C_2, \ldots, C_n$ to be some n attributes / concepts associated with the problem in hand. Now instead of having the 'on' or 'off' state of the concepts $C_1, \ldots, C_n$ as in FCM we with these concepts associate linguistic states like; always, never, often, very often, not that often, much, very much, not that much, medium, large etc. So at one time the state $C_i$ may be only one of these states and '0', when it is none of these states. Clearly these states are not values between [0, 1] but they are linguistic states and not numbers; that is why we say the dynamical system which we are defining is called as the Fuzzy Linguistic Cognitive Maps (FLCM) model.

Suppose at an instant the attribute / concept $C_i$ is in the state 'often' and its impact on the another attribute $C_j$ ($i \neq j$) increases the effect 'often' in $C_j$ (that is increase in $C_i$ increases ($C_j$)) then we map in the linguistic graph with $C_1, \ldots, C_n$ as nodes, the vertex $C_i$ to $C_j$ as 'positive often'. If decrease in $C_i$ decrease $C_j$ we mark as 'positive often'. If increase (decrease) in $C_i$ decrease (increase) $C_j$ we map as 'negative often'. However $C_i$ cannot have impact on $C_i$ so $C_i$ to $C_i$ is 0. Now our state vectors will be $X = (a_1, a_2, \ldots, a_n)$ where $a_i \in$ {linguistic variables associated with the problem} $\cup \{0\} = L$; o (L) $< \infty$, $1 \leq i \leq n$. Now using the linguistic graph we get the linguistic matrix.

If G be some fuzzy linguistic graph and M the corresponding fuzzy linguistic n × n matrix and if $x = (a_1, \ldots, a_n)$ be the state vector then we define M to be the fuzzy linguistic dynamical system; now if $M = (m_{ij})$; $m_{ij} \in L$, $1 \leq i, j \leq n$ then $xM = (b_1, \ldots, b_n)$ where $b_i \in L$; $1 \leq i \leq n$.

Let $xM = x_1$ after updating (for updating refer [Kasko WBV]); now we find $x_1M = x_2$ (say) we find $x_2M$ we continue finding such vectors and finally we arrive at a fixed point or a limit cycle since L is a finite set.

This fixed point or the limit cycle is defined as the fuzzy linguistic hidden pattern. If $x_1 = (t_1, \ldots, t_n)$ is the fixed point (or limit cycle) we see the direct impact of $x = (a_1, \ldots, a_n)$ on the problem and we can intrepet it.

We will illustrate this first by an example.

*Example 4.2.1:* Suppose we are interested in studying the child labour problem. Let $(C_1, C_2, \ldots, C_6)$ be six attributes / concepts associated with it.

$C_1$ - Child Labour
$C_2$ - Good Teacher
$C_3$ - School Drop out
$C_4$ - Poverty
$C_5$ - Public encouraging child labour
$C_6$ - Broken Family.

These concepts need not be explained as they are self explanatory we take the following fuzzy linguistic variables as the state vectors. Let L = {0, often, + often, –often, very much, much, not that much, little, very little, more etc} (Any other L can be got by the expert as per his / her need).

So the state vectors as well as the related fuzzy linguistic matrices take their values from the set L. Also the vertices of the fuzzy linguistic graph take their values from L.

Now using the experts opinion we have the following fuzzy linguistic graph.

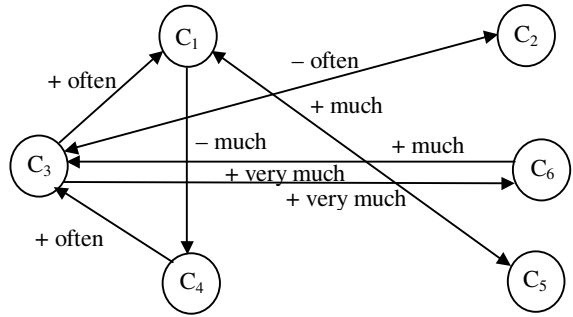

Now using this linguistic graph we have the following fuzzy linguistic matrix M; where M is the dynamical system associated with the problem of child labour.

$$M = \begin{array}{c} \\ C_1 \\ C_2 \\ C_3 \\ C_4 \\ C_5 \\ C_6 \end{array} \begin{array}{c} \begin{array}{cccccc} C_1 & C_2 & C_3 & C_4 & C_5 & C_6 \end{array} \\ \left[ \begin{array}{cccccc} 0 & 0 & 0 & +\text{often} & +\text{very much} & 0 \\ 0 & 0 & -\text{often} & 0 & 0 & 0 \\ +\text{often} & -\text{much} & 0 & 0 & 0 & +\text{much} \\ 0 & 0 & +\text{often} & 0 & 0 & 0 \\ +\text{very much} & 0 & 0 & 0 & 0 & 0 \\ 0 & 0 & +\text{much} & 0 & 0 & 0 \end{array} \right] \end{array}.$$

Now our state vector takes values from L. The min and max operations on L are defined as follows.

Let $\min\{0, a_i\} = 0$ and $\max\{0, a_i\} = a_i$ for all $a_i \in L$.

$a_i \in L$ min $\{a_i, a_i\} = a_i$ and max $\{a_i, a_i\} = a_i$.
min $\{a_i - a_i\} = -a_i$ and max $\{a_i, -a_i\} \in a_i$ for all $a_i \in L$.
For instance min {often, very often} = often
max {often, very often} = very often
min {much, often} = often,
max {much, often} = much
min {−much, often} = − much
and max {−much, often} = often.
Like this operations on L are performed.

Now suppose we want to study the effect of state vector "Child labour often occurs" under the condition 'often' $\in L$ and all other states are in the off state or zero state. To find the effect of $x = $ (+ often, 0, 0, 0, 0, 0) on the dynamical system M. $xM = $ (0, 0, 0, + often, + often, 0) after updating the state vector as the first concept was in the on state with '+ often' in the final results as well as in every step it should continue to remain in the '+ often' state. So let $x_1 = $ (+, often, 0, + often, often, 0).

Now we find $x_1M$ using as before min $\{$min $(a_i, m_{ij})\}$ where $x_1 = (a_1, a_2, \ldots, a_6)$ and $M = (m_{ij})$; $m_{ij}, a_i \in L$, $1 \le i, j \le 6$.

Thus $x_1M = $ (+ often, 0, 0, + often, + often, 0) ×

$$\begin{bmatrix} 0 & 0 & 0 & +\text{often} & +\text{very much} & 0 \\ 0 & 0 & -\text{often} & 0 & 0 & 0 \\ +\text{often} & -\text{much} & 0 & 0 & 0 & +\text{much} \\ 0 & 0 & +\text{often} & 0 & 0 & 0 \\ +\text{very much} & 0 & 0 & 0 & 0 & 0 \\ 0 & 0 & +\text{much} & 0 & 0 & 0 \end{bmatrix}$$

= (+ often, 0, + often, + often, + often, 0) leading to a fixed point.

Thus the hidden pattern clearly presents if 'often' child labour is in vogue it is 'often' poverty. $C_4$ and $C_5$, public encouraging child labour.

It is pertinent to mention here that '+' can be replaced by positive and '–' by negative. How ever an expert can choose without '+'ve or '–ve' with 0 as the least which is described in the following.

Let L be a fuzzy linguistic set associated with the n fuzzy linguistic attributes $C_1, C_2, \ldots, C_n$ of a problem P where each $C_i$ takes every fuzzy linguistic term in L at one time or the other; for i=1, 2, …, n.

Now we see $C_i$ may at one time take a fuzzy linguistic term from L and we study that influence on other attributes. Further if $C_i$ is influencing the attribute $C_j$ then ($i \neq j$) we see the influence or association is again a fuzzy linguistic term from L. The state of $C_i$ at any time can be 0 or any other fuzzy linguistic term from L.

Thus if $(a_1, a_2, \ldots, a_n) = X$ denotes the state fuzzy linguistic vector where $a_i$'s are fuzzy linguistic terms which the attributes $C_i$ enjoys at that time is $a_i$ with $a_i \in L$; $1 \leq i \leq n$. Now the influence of each of the fuzzy concepts on each other $C_1, C_2, \ldots, C_n$ is denoted by the fuzzy linguistic graph similar to the directed graph of the fuzzy cognitive maps. These graphs of FLCMs take edge values from L and the edge values are not +1 or –1 or values from [0, 1], they are basically linguistic terms from L.

Now associated with this fuzzy linguistic graph we obtain a n × n fuzzy linguistic matrix M whose entries are from L and the diagonal elements are zero.

Thus if

$$M = \begin{bmatrix} 0 & x_{12} & x_{13} & \ldots & x_{1n} \\ x_{21} & 0 & x_{23} & \ldots & x_{2n} \\ \vdots & \vdots & \vdots & & \vdots \\ x_{n1} & x_{n2} & x_{n3} & \ldots & 0 \end{bmatrix} \text{ where } x_{ij} \in L$$

(some of the $x_{ij}$'s can be zero).

M is defined as the fuzzy linguistic dynamical system associated with the problem.

If $X = (a_1, \ldots, a_n)$ where $a_i \in L$ then we can find min (max (X, M)) or max (min (X, M)) or min min (X, M) or max (max (X, M)) according to the needs of the expert and the nature of the problem at hand.

We find max (min (X, M)) (say) if X, is the value after updating that is what ever the attribute enjoyed in the beginning should continue till the end. Only the effect of other attributes which where in the 0 state are found then we once again find max (min (X, M)) if $X_2$ is the value after updating all the fuzzy linguistic states except 0's we again find max (min ($X_2$, M)) and so on until we arrive at a fixed fuzzy linguistic vector or a fuzzy linguistic vector as a limit cycle; this is always achieved as L is finite we call this fuzzy linguistic fixed point or fuzzy linguistic limit cycle as a fuzzy linguistic hidden pattern of the fuzzy linguistic dynamical system M.

We describe the working of the Fuzzy Linguistic Cognitive Model (FLCM) by the following example.

*Example 4.2.2:* Consider the problem of finding the eight transit system which includes the level of service and the convenience factors. We have the following eight attributes.

- $C_1$ - Frequency of the service along a route
- $C_2$ - In-vehicle travel time along the route
- $C_3$ - Travel fare along the route
- $C_4$ - Speed of the vehicles along the route
- $C_5$ - Number of intermediate points in the route.
- $C_6$ - Waiting time
- $C_7$ - Number of transfers in the route
- $C_8$ - Congestion in the vehicle.

The fuzzy linguistic terms associated with $C_1, C_2, \ldots, C_8$ as well as the problem are L = {0, often, always, a little, much, very much, usually, some times}. (It is pertinent to mention here more number of fuzzy linguistic terms or lesser number of fuzzy linguistic terms can be associated which is the solely the wish of the expert). We now give the fuzzy linguistic graph whose vertices are $C_1, C_2, \ldots, C_8$ and the edges take values from L are as follows.

Frequency of the vehicle along the route often depends on the in-vehicle travel time along the route; for if we have 'often' that is the frequency is more we see the in-vehicle travel time will also be less as the crowd will be less. Likewise the other attributes are all related in a similar way. We now give the corresponding fuzzy linguistic matrix of the fuzzy linguistic graph.

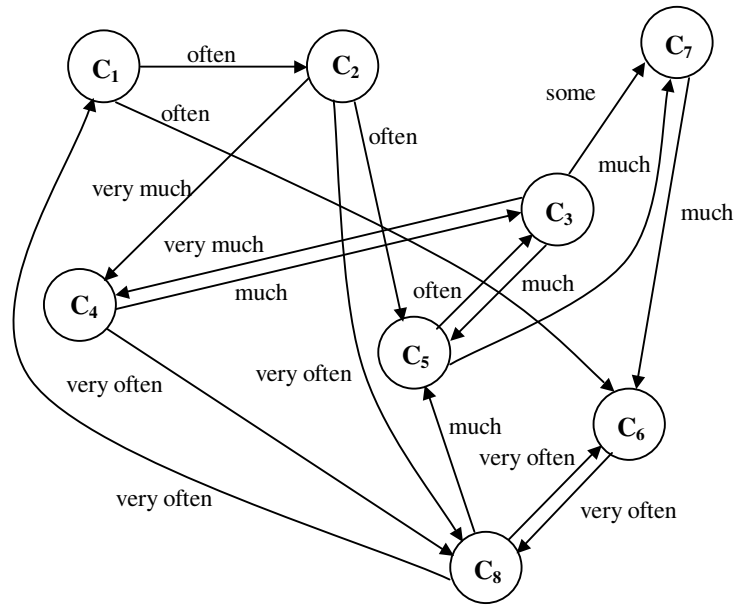

$$M = \begin{array}{c} \\ C_1 \\ C_2 \\ C_3 \\ C_4 \\ C_5 \\ C_6 \\ C_7 \\ C_8 \end{array} \begin{array}{c} \begin{array}{cccccccc} C_1 & C_2 & C_3 & C_4 & C_5 & C_6 & C_7 & C_8 \end{array} \\ \left[ \begin{array}{cccccccc} 0 & \text{often} & 0 & 0 & 0 & \text{often} & 0 & 0 \\ 0 & 0 & 0 & \text{very much} & \text{often} & 0 & 0 & \text{very often} \\ 0 & 0 & 0 & \text{very much} & \text{much} & 0 & \text{some} & 0 \\ 0 & 0 & \text{much} & 0 & 0 & 0 & 0 & \text{very often} \\ 0 & 0 & \text{often} & 0 & 0 & 0 & \text{much} & 0 \\ 0 & 0 & 0 & 0 & 0 & 0 & 0 & \text{very often} \\ 0 & 0 & 0 & 0 & \text{much} & 0 & 0 & 0 \\ \text{very often} & 0 & 0 & 0 & \text{much} & \text{very often} & 0 & 0 \end{array} \right] \end{array}.$$

M is defined as the fuzzy linguistic dynamical system associated with the problem. Now we find the hidden pattern of any fuzzy linguistic state vector X on M, where X = (often, 0, some, 0, often, some, 0, often) is the fuzzy linguistic state vector supplied by an expert. To study the effect of X on M using max {max {X, M}} after updating is (often, often, some, very much, often, some, much, often) = $X_1$.

Now max {max {$X_1$, M}} after updating is say $X_2$;

$X_2$ = (often, very much, some, very much, often, some, very much, often).

Now we find max {max ($X_2$, M)} after updating we get say $X_3$;

$X_3$ = (often, very much, some, very much, often, some, very much, often).

We see the fuzzy linguistic hidden pattern of the fuzzy linguistic state vector X on the fuzzy linguistic dynamical system M is a fixed point given by

$X_3$ = (often, very much, some, very much, often, some, very much, often).

Thus if the experts takes the fuzzy linguistic states that frequency of the vehicles along the route is 'often', and the travel fare is 'some' i.e.; not much or very much and the number of intermediate points in the route is often (not very often or much or very much) and 'some' time is the waiting time with congestion in the vehicle is often we see; the fuzzy linguistic resultant is that the invehicle travel time along the route is very much, speed along the route is also very much and the number of transfers in that route is very much. Thus we have illustrated the model with max {max} operation.

Now suppose another expert wants to use for the same fuzzy linguistic state vector X and the same dynamical system, to study the effect using the 'max min' operation.

Now max {min (X, M)} after updating we get say $X_1$,

$X_1$ =(often, often, some, 0, often, some, often, often).

Now we find max {min (X, M)} after updating we get say $X_2$;

$X_2$ = (often, often, some, 0, often, some often, often) which is a fuzzy linguistic hidden pattern of the fuzzy linguistic state vector X under 'max min' operation.

Now when we use 'max min' operation on the fuzzy linguistic state vector X the effect of speed of the vehicle is '0' that is the fuzzy linguistic position of $C_4$ remains unaffected further the number of transfer in that route is only often.

Thus an expert may prefer 'max min' operation on the fuzzy dynamical system instead of the 'max max' operation.

This is only an illustration of how the fuzzy linguistic cognitive model functions.

We can use 'min min' or 'min max' operation also. It is at the liberty of the expert to use any of the operation to study the problem.

### 4.3 Fuzzy Linguistic Relational Map Model

Here we introduce yet another new fuzzy linguistic model known as the fuzzy linguistic relational map model. This model is more a generalization of the fuzzy linguistic cognitive model. For with any problem P we have an association of the attributes or concepts which are related with fuzzy linguistic terms. If the attributes / concepts can be divided into two disjoint spaces say D and R where D is the fuzzy linguistic attributes known as the fuzzy linguistic domain space of the problem and R the fuzzy linguistic range space of the problems are let L be the fuzzy linguistic space which takes all the fuzzy linguistic terms associated both with R and D.

We see if $(D_1, D_2, \ldots, D_m)$ are the concepts in the fuzzy linguistic domain space D and $(R_1, R_2, \ldots, R_n)$ are the fuzzy linguistic concepts of the fuzzy linguistic range space R, then the fuzzy linguistic graph with $\{D_1, D_2, \ldots, D_m$ and $R_1, R_2, \ldots, R_n\}$ as the vertices and edges will be between $D_i$'s and $R_j$'s only $1 \leq i \leq m$, $1 \leq j \leq n$ and the fuzzy linguistic edges will take its fuzzy linguistic values from the fuzzy linguistic space L.

Associated with this fuzzy linguistic graph we will get the fuzzy linguistic matrix. The related fuzzy linguistic matrix will be known as the fuzzy linguistic dynamical system associated with the fuzzy linguistic relational map model.

We will have a pair of fuzzy linguistic state vectors one associated with the fuzzy linguistic domain space D and the other associated with the fuzzy linguistic range space R. As in case of fuzzy linguistic cognitive models we can get the hidden pattern which will be a pair of fuzzy linguistic state vectors associated with one fuzzy linguistic state vector from D or R. As in case of fuzzy linguistic cognitive models the fuzzy linguistic hidden pattern can be a fixed point or a limit cycle.

We can use any of the four types of operations 'max, max' or {min, max} or {max, min} or {min, min} to find the fuzzy linguistic hidden pattern.

We will illustrate this situation by an example.

***Example 4.3.1:*** Let us study the employee - employer fuzzy linguistic relational model. Suppose we have the following fuzzy linguistic concepts / attributes associated with the employee taken as the domain space

- $D_1$ - Pay with allowances and bonus
- $D_2$ - Only pay to employee
- $D_3$ - Pay with allowance to employee
- $D_4$ - Best performance by the employee
- $D_5$ - Average performance by the employee
- $D_6$ - Employee works for more number for hours
- $D_7$ - Employee works for less number of hours.

Suppose the following nodes / concepts are taken as the range space of the employer.

- $R_1$ - maximum profit to the employer
- $R_2$ - Only profit to the employer
- $R_3$ - Neither profit nor less to the employer
- $R_4$ - loss to the employer
- $R_5$ - Heavy loss to the employer

The fuzzy linguistic terms associated with the fuzzy linguistic domain and range spaces be taken as L.

L = {0, gain, loss, no loss no gain, just gain, just loss, gain, heavy loss, good gain}.

Any expert can suggest / take any other fuzzy linguistic terms as L. It is solely left for the expert to choose whatever he / she wishes to work with. We give the associated fuzzy linguistic graph of the problem using both the fuzzy linguistic domain space and range space.

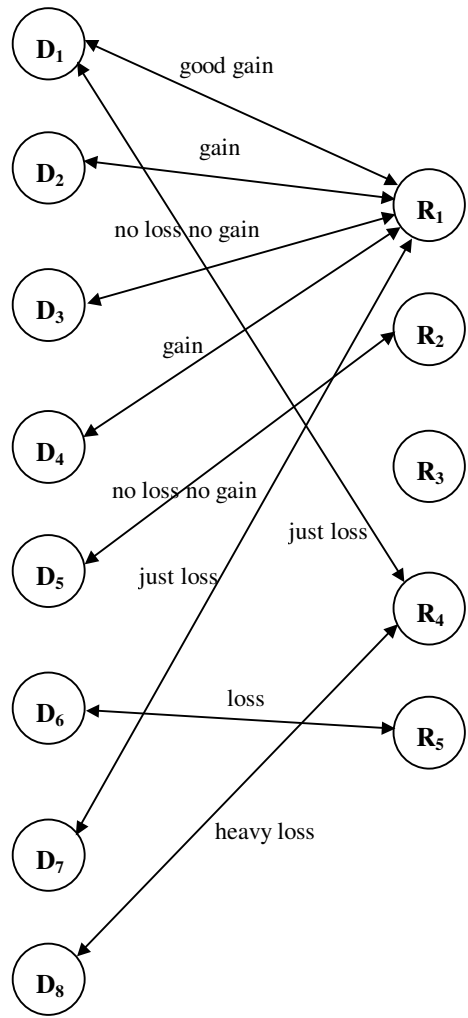

We give the associated fuzzy linguistic matrix N of the fuzzy linguistic graph.

$$N = \begin{array}{c} \\ D_1 \\ D_2 \\ D_3 \\ D_4 \\ D_5 \\ D_6 \\ D_7 \\ D_8 \end{array} \begin{bmatrix} R_1 & R_2 & R_3 & R_4 & R_5 \\ \text{good gain} & 0 & 0 & \text{just loss} & 0 \\ \text{gain} & 0 & 0 & 0 & 0 \\ \text{no loss no gain} & 0 & 0 & 0 & 0 \\ \text{gain} & 0 & 0 & 0 & 0 \\ 0 & \text{no loss no gain} & 0 & 0 & 0 \\ 0 & 0 & 0 & 0 & \text{loss} \\ \text{just loss} & 0 & 0 & 0 & 0 \\ 0 & 0 & 0 & \text{heavy loss} & 0 \end{bmatrix}.$$

Now we find the resultant of any fuzzy linguistic vector on N, the dynamical system associated with the problem.

Let $X$ = (gain, 0, loss, 0, gain, loss, 0, 0) be the given fuzzy linguistic state vector. The effect of $X$ on the fuzzy linguistic dynamical system N is as follows:

Max {max $(X, N)$} = ((good gain, gain, gain, gain, gain)} = $Y$

We find max (max $\{Y, N^T\}$) = (gain, gain, loss, gain, gain, loss, gain, gain) = $X_1$ (say) after updating.

We find max (max $\{X_1, N\}$)

= (good gain, gain, gain, gain) = $Y_1 = Y$.

Thus the fuzzy linguistic hidden pattern is a fixed pair and the employer always get gain when the attribute $D_1$ enjoys the fuzzy linguistic state 'gain' with 'loss' for the attribute $D_3$; 'gain' by average performance by the employee - $D_4$ and the fuzzy linguistic term loss if the performance of the employee is poor.

Under these conditions if the employer gives pay with allowance and bonus to the worker the employer gets 'gain'.

This is only an illustrative example and interested researchers can apply these models to real world problems.

### 4.4 Fuzzy Linguistic Relation Equation

In this section we for the first time introduce the notion of fuzzy linguistic relation equations and describe some of the properties related with them.

Let $L_1$ and $L_2$ be any two fuzzy linguistic sets that is both $L_1$ and $L_2$ contain fuzzy linguistic terms or $L_1 = L_2$ otherwise. Let R be a fuzzy linguistic relation that is to each fuzzy linguistic term of $L_1$ two or more fuzzy linguistic terms in $L_2$ are assigned, that is to each $x_1 \in L_1$ to the domain space ($L_1$ is the domain linguistic space) we associate a $y_2 \in L_2$ and here the degree of membership is not a value between [0, 1] but a fuzzy linguistic value from L.

This we will first illustrate by some example.

***Example 4.4.1:*** Let $L_1 = \{x_1, x_2, x_3, x_4, x_5, x_6, x_7\}$ be seven fuzzy linguistic terms and $L_2 = \{y_1, y_2, y_3, y_4, y_5, y_6\}$ be some six fuzzy linguistic terms. These two fuzzy linguistic sets are fuzzy linguistically related by mapped by R where R takes its fuzzy linguistic membership values as fuzzy linguistic terms from a fuzzy linguist set L.

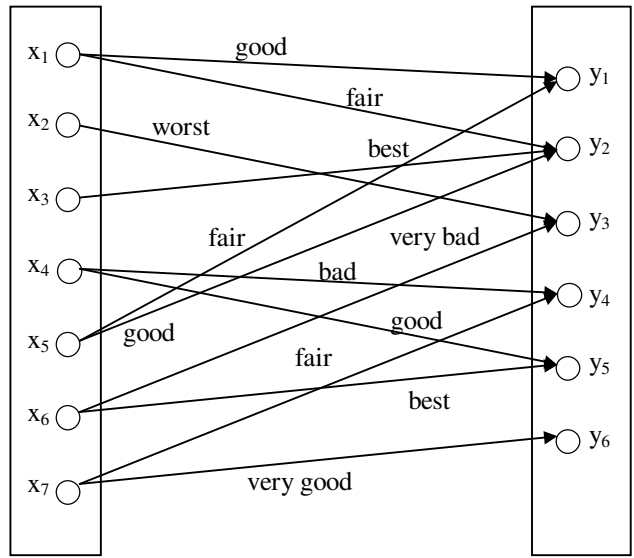

The corresponding fuzzy linguistic binary relation. The fuzzy linguistic membership matrix R is as follows.

$$R = \begin{array}{c} \\ x_1 \\ x_2 \\ x_3 \\ x_4 \\ x_5 \\ x_6 \\ x_7 \end{array} \begin{bmatrix} y_1 & y_2 & y_3 & y_4 & y_5 & y_6 \\ good & fair & 0 & 0 & 0 & 0 \\ 0 & 0 & worst & 0 & 0 & 0 \\ 0 & fair & 0 & 0 & 0 & 0 \\ 0 & 0 & 0 & bad & good & 0 \\ best & good & 0 & 0 & 0 & 0 \\ 0 & 0 & v.bad & 0 & best & 0 \\ 0 & 0 & 0 & fair & 0 & v.good \end{bmatrix}$$

Thus for fuzzy linguistic sets $L_1$ and $L_2$ we can have a fuzzy linguistic membership function R. The representation by the diagram will be known as the fuzzy linguistic sagittal diagram and the fuzzy linguistic matrix will be known as the fuzzy linguistic membership matrix where the membership values are also fuzzy linguistic terms.

We can compose two fuzzy linguistic binary relation provided both of them take their fuzzy linguistic membership values from the same fuzzy linguistic set L.

We will just illustrate how such compositions are made.

Suppose

$$P = \begin{bmatrix} good & bad & fair & best \\ bad & fair & good & good \\ 0 & good & fair & good \\ good & bad & good & fair \end{bmatrix}$$

and

$$Q = \begin{bmatrix} good & bad & good & fair & 0 & best \\ bad & good & best & good & bad & 0 \\ best & fair & best & 0 & good & bad \\ 0 & fair & good & bad & good & good \end{bmatrix}$$

be two fuzzy linguistic membership matrices of a fuzzy linguistic binary relation with the values taken from the same fuzzy linguistic space L.

Now

$$P \circ Q = \begin{bmatrix} good & bad & fair & best \\ bad & fair & good & good \\ 0 & good & fair & good \\ good & bad & good & fair \end{bmatrix} \circ$$

$$\begin{bmatrix} good & bad & good & fair & 0 & best \\ bad & good & best & good & bad & 0 \\ best & fair & best & 0 & good & bad \\ 0 & fair & good & bad & good & good \end{bmatrix}$$

= R where R is found in the following way.

$$\text{If } R = \begin{bmatrix} r_1 & r_2 & r_3 & r_4 & r_5 & r_6 \\ r_7 & r_8 & r_9 & r_{10} & r_{11} & r_{12} \\ r_{13} & r_{14} & r_{15} & r_{16} & r_{17} & r_{18} \\ r_{19} & r_{20} & r_{21} & r_{22} & r_{23} & r_{24} \end{bmatrix} \text{ then}$$

$r_1$ = max [min (good, good), min (bad, bad), min (fair, best), min (best, 0)]
= max [good, bad, fair, 0] = good

$r_2$ = max [min (good, bad), min (bad, good), min (fair, fair), min (best, fair)}
= max {bad, bad, fair, fair} = fair

$r_3$ = max {min {good, good}, min {bad, best}, min {fair, best}, min {best, good}}
= max {good, bad, fair, good} = good

$r_4$ = max {min {good, fair}, min {bad, good}, min {fair, 0}, min {best, bad}}
= max {fair, bad, 0, bad } = fair.

$r_5$ = max {min {good, 0}, min {bad, bad},
       min {fair, good}, min {best, good}}
     = max {0, bad, fair, good} = good

$r_6$ = max {min {good, best}, min {bad, 0},
       min {fair, bad}, min {best, good}}
     = max {good, 0, bad, good} = good

$r_7$ = max {min {bad, good}, min {fair, bad},
       min {good, best}, min {good, 0}}
     = max {bad, bad, good, 0} = good

$r_8$ = max {min {bad, bad}, min {fair, good},
       min {good, fair}, min {good, fair}}
     = max {bad, fair, fair, fair} = fair

$r_9$ = max {min {bad, good}, min {fair, bad},
       min {good, best}, min {good, good}}
     = max {bad, fair, good, good} = good

$r_{10}$ = max {min {bad, fair}, min {fair, good},
        min {good, 0}, min {good, bad}}
      = max {bad, fair, 0, bad} = fair

$r_{11}$ = max {min {bad, 0}, min {fair, bad},
        min {good, good}, min {good, good}}
      = max {0, bad, good, good} = good

$r_{12}$ = max {min {bad, best}, min {fair, 0},
        min {good, bad}, min {good, good}}
      = max {bad, 0, 0, bad, good} = good

$r_{13}$ = max {min {0, good}, min {good, bad},
        min {fair, best}, min {good, 0}}
      = max {0, bad, fair, 0} = fair

$r_{14}$ = max {min {0, bad}, min {good, good},
        min {fair, fair}, min {good, fair}}
      = max {0, good, fair, fair} = good

$$r_{15} = \max \{\min \{0, good\}, \min \{good, best\},$$
$$\min \{fair, best\}, \min \{good, good\}\}$$
$$= \max \{0, good, fair, good\} = good$$

$$r_{16} = \max \{\min \{0, fair\}, \min \{good, good\},$$
$$\min \{fair, 0\}, \min \{good, bad\}\}$$
$$= \max \{0, good, 0, bad\} = good$$

$$r_{17} = \max \{\min \{0, 0\}, \min \{good, bad\}, \min \{fair, good\},$$
$$\min \{good, good\}\}$$
$$= \max \{0, bad, fair, good\} = good$$

$$r_{18} = \max \{\min \{0, best\}, \min \{good, 0\}, \min \{fair, bad\},$$
$$\min \{good, good\}\}$$
$$= \max \{0, 0, bad, good\} = good$$

$$r_{19} = \max \{\min \{good, good\}, \min \{bad, bad\},$$
$$\min \{good, best\}, \min \{fair, 0\}\}$$
$$= \max \{good, bad, good, 0\} = good$$

$$r_{20} = \max \{\min \{good, bad\}, \min \{bad, good\},$$
$$\min \{good, fair\}, \min \{fair, fair\}\}$$
$$= \max \{bad, bad, fair, fair\} = fair$$

$$r_{21} = \max \{\min \{good, good\}, \min \{bad, best\},$$
$$\min \{good, best\}, \min \{fair, good\}\}$$
$$= \max \{good, bad, good, fair\} = good$$

$$r_{22} = \max \{\min \{good, fair\}, \min \{bad, good\},$$
$$\min \{good, 0\}, \min \{fair, bad\}\}$$
$$= \max \{fair, bad, 0, bad\} = fair$$

$$r_{23} = \max \{\min \{good, 0\}, \min \{bad, bad\},$$
$$\min \{good, good\}, \min \{fair, good\}\}$$
$$= \max \{0, bad, good, fair\} = good$$

$$r_{24} = \max \{\min \{good, best\}, \min \{bad, 0\},$$
$$\min \{good, bad\}, \min \{bad, 0\}\}$$
$$= \max \{good, 0, bad, fair\} = good$$

$$\text{Now } R = \begin{bmatrix} \text{good} & \text{fair} & \text{good} & \text{fair} & \text{good} & \text{good} \\ \text{good} & \text{fair} & \text{good} & \text{fair} & \text{good} & \text{good} \\ \text{fair} & \text{good} & \text{good} & \text{good} & \text{good} & \text{good} \\ \text{good} & \text{fair} & \text{good} & \text{fair} & \text{good} & \text{good} \end{bmatrix}.$$

Now we can also use the operation of min {max}.

$$\text{For instance if } P \circ Q = t = \begin{bmatrix} T_1 & T_2 & \ldots & T_6 \\ T_7 & T_8 & \ldots & T_{12} \\ T_{13} & T_{14} & \ldots & T_{18} \\ T_{19} & T_{20} & \ldots & T_{24} \end{bmatrix}$$

using min max.

$T_1$ = min {max {good, good}, max {bad, bad}, max {fair, best}, max {best, 0}}
 = min {good, bad, best, best} = bad.

We see $r_1 \neq T_1$.

$T_2$ = min {max {bad, good}, max {fair, bad}, max {good, best}, max {good, 0}}
 = min {good, fair, best, good} = fair and so on.

Thus we need not have $T = R$.

So depending on the problem at hand we have to choose the operations.

We can as in case of usual binary relations in case of fuzzy linguistic relations have the concept of ternary fuzzy linguistic relation from two fuzzy linguistic binary relations, which results in a fuzzy linguistic binary relation.

This is illustrated by the following diagram.

Let $L_1$, $L_2$ and $L_3$ be three fuzzy linguistic sets with the following fuzzy linguistic sagittal diagram.

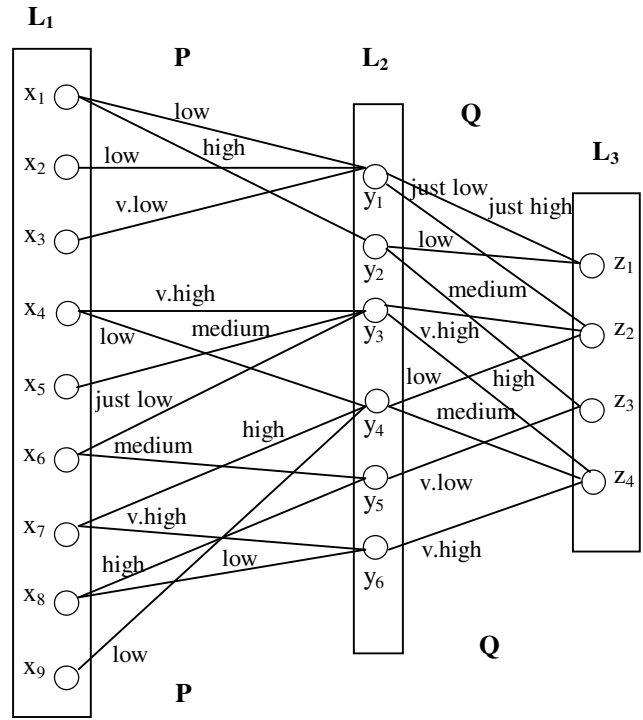

Let L the fuzzy linguistic set be {0, low, very low, high, medium, very high, just low, just high, just medium}.

We have S = P o Q.

We can also as in case of binary relation on a single set we can find fuzzy linguistic binary relation using a fuzzy linguistic set X and L is the fuzzy linguistic set which works on X. L = {0, fast, very fast, very slow, just slow, just fast, medium}.

We have the following fuzzy linguistic binary relation which is described by the fuzzy linguistic sagittal diagram and the related fuzzy linguistic membership matrix.

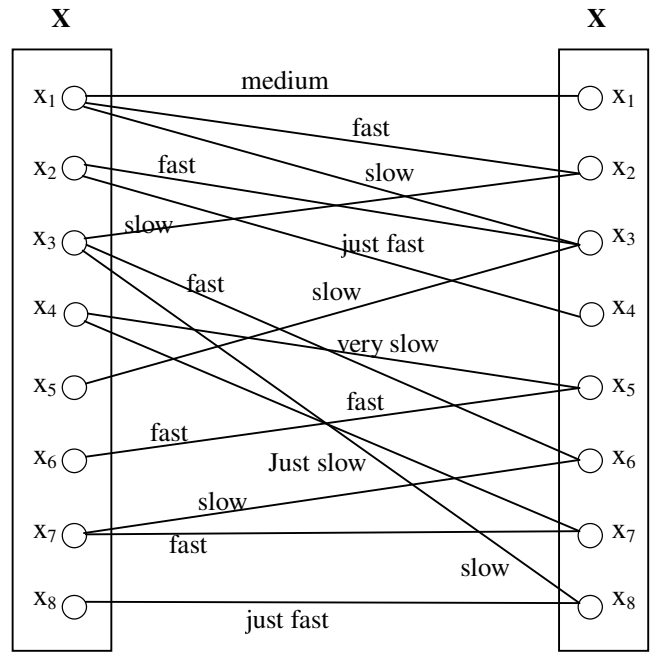

The fuzzy linguistic membership matrix associated with X is

|   | $x_1$ | $x_2$ | $x_3$ | $x_4$ | $x_5$ | $x_6$ | $x_7$ | $x_8$ |
|---|---|---|---|---|---|---|---|---|
| $x_1$ | medium | fast | slow | 0 | 0 | 0 | 0 | 0 |
| $x_2$ | 0 | 0 | fast | just fast | 0 | 0 | 0 | 0 |
| $x_3$ | 0 | slow | 0 | 0 | 0 | fast | 0 | slow |
| $x_4$ | 0 | 0 | 0 | 0 | very slow | 0 | fast | 0 |
| $x_5$ | 0 | 0 | slow | 0 | 0 | 0 | 0 | 0 |
| $x_6$ | 0 | 0 | 0 | 0 | fast | 0 | 0 | 0 |
| $x_7$ | 0 | 0 | 0 | 0 | 0 | slow | fast | 0 |
| $x_8$ | 0 | 0 | 0 | 0 | 0 | 0 | 0 | just fast |

.

Now we can describe how the membership fuzzy linguistic functions works.

Let P (X, Y), Q (Y, Z) and R (X, Z) be defined on three fuzzy linguistic set X, Y and Z using the single fuzzy linguistic set l. Suppose P, Q and R be the fuzzy linguistic membership functions associated with them.

$$\text{Suppose } P \circ Q = R \quad \ldots \text{ I}$$

that is max min $(a_{ij}, b_{jk}) = r_{ik}$ where $a_{ij}$ are elements of P, $b_{jk}$ elements of Q and $r_{ik}$ are elements of R (it is to be noted $a_{ij}$, $b_{jk}$, $r_{ik}$ are only fuzzy linguistic terms).

We can solve I given the fuzzy linguistic matrices Q and R we want to find P. This is not possible for we are yet to define inverse for fuzzy linguistic matrices what best can be done if Q corresponds to fuzzy linguistic membership function from Y to Z then we define $Q^{-1}$ as the fuzzy linguistic membership function from Z to Y, sometimes we see we cannot find the inverse of Q. To solve I we have to use the new notion of feed forward fuzzy linguistic neural network. Hence we proceed onto define the new notion of fuzzy linguistic neural networks as well as the fuzzy linguistic layer neural networks [12, 14-5].

In case of fuzzy linguistic neural networks we see the weights are from the fuzzy linguistic set.

In the feed forward neural network with m inputs and only one layer with n-neurons if we replace the weights $w_{ij}$ by the fuzzy linguistic terms, where $\{x_1, \ldots, x_n\}$ and $\{y_1, \ldots, y_m\}$ are also fuzzy linguistic terms we get the fuzzy linguistic neural network [12, 14-5].

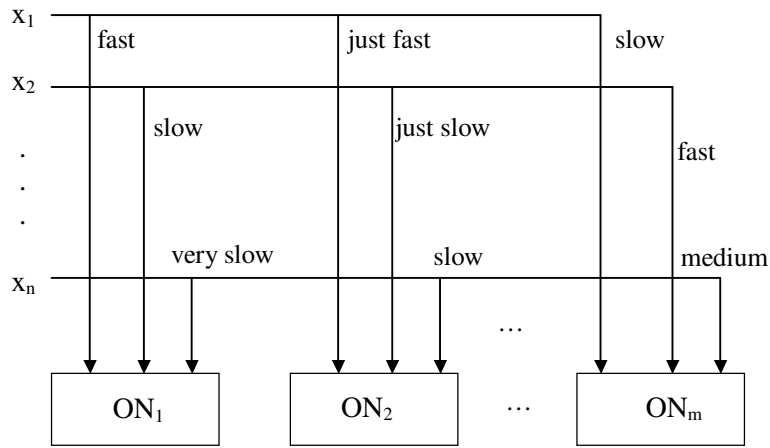

We can get the fuzzy linguistic feed forward neural network with two layers. Let $\{x_0, x_1, \ldots, x_n\}$ and $\{y_1, y_2, \ldots, y_m\}$ be fuzzy linguistic input and $\{y_1, y_2, \ldots, y_m\}$ be the expected fuzzy linguistic outputs $w_{ij}$'s be fuzzy linguistic (weights) terms; $1 \leq i \leq m$, $1 \leq j \leq n$. It is pertinent to mention here instead of using weights $l(w_{ij})$ from the interval [0, 1], we are using fuzzy linguistic terms like {low, high, very low, very high, just low, medium, 0, so on} or {fast, slow, very fast, 0, etc} or {more, less, etc} depending upon the problem at hand. Keeping these information we have the feed forward neutral network with two layers.

It is pertinent to mention here we give the fuzzy linguistic weights $l(w_{ij})$ as well as $l(W_{ij})$ only by trial and error method.

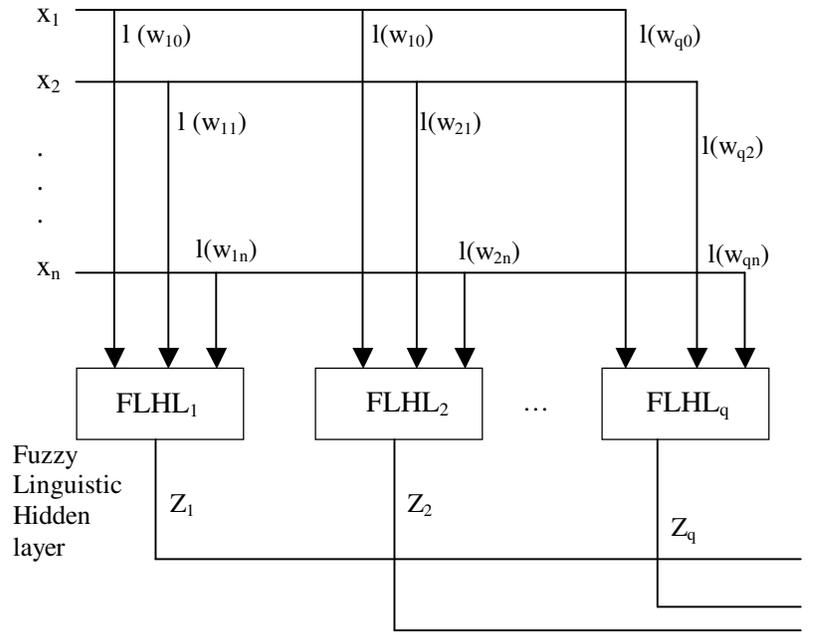

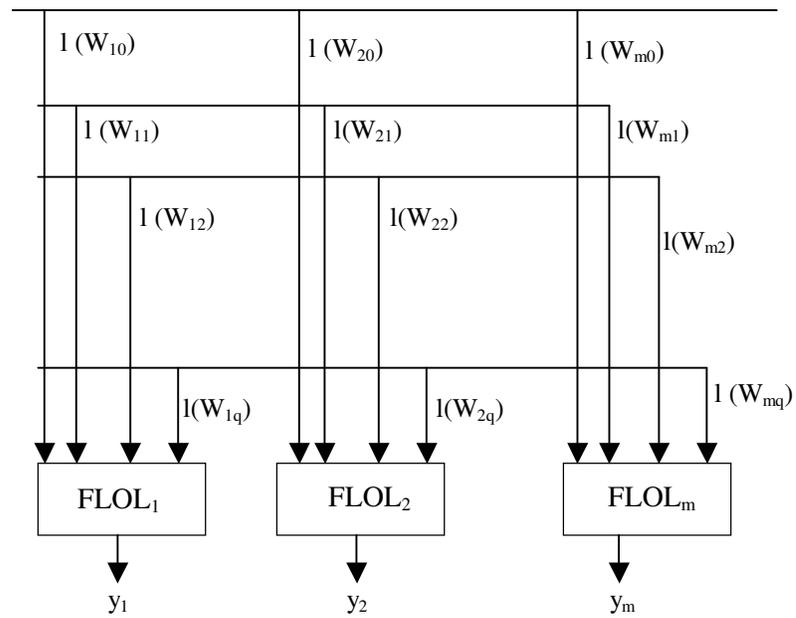

Fuzzy Linguistic Output layer

Further we can have likewise 3 fuzzy linguistic hidden layers or n fuzzy linguistic hidden layers and only one output fuzzy linguistic layer. While using the fuzzy linguistic feed forward neural network the following conditions are imposed.

1. All the matrices P, Q and R are fuzzy linguistic membership matrices have fuzzy linguistic terms associated with the problem.

2. These fuzzy membership matrices can be represented by the fuzzy linguistic relation equation P o Q = R
(We can have max or min or any one of the functions can be used as per need of the problem).

Instead of the activation function we have the linguistic function $f_2$ taking its values from the linguistic space L. We assume our space L is a non over lapping fuzzy linguistic space. This assumption is basically made to make a partition of these fuzzy linguistic term (space). Since we want L to be an orderable space we take L to be a chain lattice connected topological fuzzy linguistic space.

Now we instead of working with [0, 1], work with L which is just like [0, 1].

Throughout we will have our membership matrices to be fuzzy linguistic membership matrices and all entries and final output layer will only be fuzzy linguistic terms.

These fuzzy linguistic relational equations model can be used in the analysis of social problems, medical diagnosis, problems pertaining to chemical industries and so on.

For the criteria is we can transform the values in the interval [0, 1] into fuzzy linguistic terms in L and vice versa provided L is a chain lattice connected fuzzy linguistic topological space.

**Chapter Five**

# SUGGESTED PROBLEMS

In this chapter we suggest over 100 problems for the reader some of which are difficult and are research problems.

1.  Give some examples of fuzzy linguistic spaces.

2.  Show the lattice associated with the fuzzy linguistic space need not always be a chain lattice.

3.  Obtain some nice properties about fuzzy linguistic spaces.

4.  Give an example of a fuzzy linguistic topological space by defining a topology on it.

5.  Give an example of a lattice connected fuzzy linguistic topological space.

6.  Give an example of a lattice disconnected fuzzy linguistic topological space.

7.  Show that a fuzzy linguistic space cannot always be associated with a lattice.

8. What is a over lapping fuzzy linguistic space?

9. Give an example of a fuzzy linguistic over lapping space.

10. Construct a fuzzy linguistic topological space associated with the bonded labour problem.

11. Draw the fuzzy linguistic graph associated with the bonded labour problem.

12. Is the fuzzy linguistic graph given by you a connected fuzzy linguistic graph in problem 11.

13. Give an example of a fuzzy linguistic topological space which is graph connected.

14. Does there exist a fuzzy linguistic topological space which is not graph connected?

15. Is every fuzzy linguistic topological space discrete?

16. Obtain some interesting properties about fuzzy linguistic topological space which are graph connected.

17. Give an example of a fuzzy linguistic topological space which is totally graph disconnected.

18. Does there exist a fuzzy linguistic topological space which is not lattice connected?

19. Give an example of a fuzzy linguistic topological space L which is discrete. Is every fuzzy linguistic topological subspace of a fuzzy linguistic topological space L discrete?

20. Develop some interesting properties about fuzzy linguistic topological discrete spaces.

21. Determine some interesting features of a fuzzy linguistic topological space which is not lattice connected.

22. Obtain some interesting properties about graph connected fuzzy linguistic topological spaces.

23. Is every strongly graph connected fuzzy linguistic topological space a graph connected fuzzy linguistic topological space?

24. Let

$$V = \left\{ \begin{bmatrix} a_1 \\ a_2 \\ a_3 \\ a_4 \end{bmatrix}, (a_1, a_2, a_3), \begin{bmatrix} a_1 & a_2 \\ a_3 & a_4 \end{bmatrix} \middle| a_i \in L, 1 \leq i \leq 4 \right\}$$

(L a fuzzy linguistic space of order 10) be a set fuzzy linguistic vector space over the set L.

(i) Find the number of elements in V.

(ii) Write V as a direct sum of set fuzzy linguistic vector subspaces.

(iii) Find a set linear operator on V.

25. Let

$$V = \left\{ \begin{bmatrix} a_1 & a_2 & \ldots & a_{10} \\ a_{11} & a_{12} & \ldots & a_{20} \\ a_{21} & a_{22} & \ldots & a_{30} \\ a_{31} & a_{32} & \ldots & a_{40} \\ a_{41} & a_{42} & \ldots & a_{50} \end{bmatrix} \middle| a_i \in L, 1 \leq i \leq 50 \right\}$$

be a semigroup fuzzy linguistic vector space over the fuzzy linguistic semigroup S.

(i) Write V as a direct sum of semigroup fuzzy linguistic vector subspaces.

(ii) Find a semigroup fuzzy linguistic linear operator T on V so that $T^{-1}$ exists.

26. Let

$$V = \left\{ (a_1, a_2, a_3, a_4, a_5), \begin{bmatrix} a_1 \\ a_2 \\ \vdots \\ a_{12} \end{bmatrix}, \begin{bmatrix} a_1 & a_2 & \ldots & a_{12} \\ a_{13} & a_{14} & \ldots & a_{24} \\ a_{25} & a_{26} & \ldots & a_{36} \\ a_{37} & a_{38} & \ldots & a_{48} \\ a_{49} & a_{50} & \ldots & a_{60} \\ a_{61} & a_{62} & \ldots & a_{72} \end{bmatrix} \right.$$

$a_i \in L, 1 \leq i \leq 72\}$ be a set fuzzy linguistic vector space over L.

(i) Write V as a pseudo direct sum of set fuzzy linguistic vector subspaces.

(ii) Find a set linear T on V.

(iii) Find a set linear functional on V.

27. Study the set of all maps f from V into L where V is the set fuzzy linguistic vector space and L is a fuzzy linguistic set denoted by L (V, L).

28. Let V =

$$\left\{ \begin{pmatrix} a_1 & a_2 \\ a_3 & a_4 \end{pmatrix}, \begin{bmatrix} a_1 & a_2 & a_3 & a_4 \\ a_5 & a_6 & a_7 & a_8 \\ a_9 & a_{10} & a_{11} & a_{12} \\ a_{13} & a_{14} & a_{15} & a_{16} \end{bmatrix}, \begin{pmatrix} a_1 & a_2 & \ldots & a_{20} \\ a_{21} & a_{22} & \ldots & a_{40} \\ a_{41} & a_{42} & \ldots & a_{60} \end{pmatrix}, \begin{bmatrix} a_1 & a_2 & a_3 \\ a_4 & a_5 & a_6 \\ \vdots & \vdots & \vdots \\ a_{28} & a_{29} & a_{30} \end{bmatrix} \right\}$$

$a_i \in L$, L a fuzzy linguistic set} be a set fuzzy linguistic vector space over L.

   (i) Find a set linear operator on V.

   (ii) Find a set linear functional on V.

   (iii) Write V as a direct sum of set fuzzy linguistic vector subspaces.

   (iv) Write V as a pseudo direct sum of set fuzzy linguistic vector subspaces.

29. If V be a set fuzzy linguistic vector space over a fuzzy linguistic set L. Find the structure of $\text{Hom}_L$ (V, V).

30. Find the structure enjoyed by the set of all set fuzzy linguistic linear operators on a set fuzzy linguistic vector space V over L; that is what is the structure of $\text{Hom}_L$ (V, V)?

31. Let V be a semigroup fuzzy linguistic vector space over the fuzzy linguistic semigroup L. W another (W ≠ V) semigroup fuzzy linguistic vector space over the same fuzzy linguistic semigroup L. Find $\text{Hom}_L$ (V, W).

32. Study the algebraic structure enjoyed by L (V, L) where
$$V = \left\{ \sum_{i=10}^{20} a_i x^i \;\middle|\; a_i \in L, \; 1 \leq i \leq 20 \right\}, \text{ a semigroup fuzzy}$$
linguistic linear algebra over the fuzzy linguistic semigroup L.

33. Study the algebraic structure of all set fuzzy linguistic linear operators on a set fuzzy linguistic vector space V over L.

34. Let
$$V = \left\{ \begin{bmatrix} a_1 & a_2 & a_3 & a_4 \\ a_5 & a_6 & a_7 & a_8 \\ a_9 & a_{10} & a_{11} & a_{12} \end{bmatrix} \;\middle|\; a_i \in L, \; 1 \leq i \leq 12 \right\}$$

be a set fuzzy linguistic matrix linear algebra over the fuzzy linguistic set L under 'min' operation.

(i) Find subspaces of V over L.

(ii) Write V as a direct sum of subspaces.

(iii) Write V as a pseudo direct sum of subspaces.

(iv) Find for a subspace W of V, $W^\perp$ so that $V = W \oplus W^\perp$.

35. Let $V = \left\{ \sum_{i=0}^{8} a_i x^i \;\middle|\; a_i \in L \right\}$ be a set fuzzy linguistic vector space over the fuzzy linguistic set L.

(i) Find subspaces if any of V.

(ii) Does V contain subset vector subspaces for any fuzzy linguistic set L?

(iii) Can V be written as a direct sum of subspaces?

(iv) Find an invertible set fuzzy linguistic linear operator on V.

36. Let V = {($a_1$, $a_2$, $a_3$, $a_4$, $a_5$),

$$\begin{bmatrix} d_1 \\ d_2 \\ \vdots \\ d_{10} \end{bmatrix}, \begin{bmatrix} e_1 & e_2 & \ldots & e_7 \\ e_8 & e_9 & \ldots & e_{14} \\ e_{15} & e_{16} & \ldots & e_{21} \end{bmatrix}, \begin{bmatrix} m_1 & m_2 & m_3 & m_4 & m_5 & m_6 \\ m_7 & m_8 & \ldots & \ldots & \ldots & m_{12} \\ m_{13} & m_{14} & \ldots & \ldots & \ldots & m_{18} \\ m_{19} & m_{20} & \ldots & \ldots & \ldots & m_{24} \\ m_{25} & m_{26} & \ldots & \ldots & \ldots & m_{30} \end{bmatrix}$$

$a_i$, $d_j$, $e_k$, $m_s \in L$, $1 \leq i \leq 5$, $1 \leq j \leq 5$, $1 \leq k \leq 21$ and $1 \leq s \leq 30$} be a set fuzzy linguistic matrix vector space over the fuzzy linguistic set L.

(i) Find a set linear operator T on V so that $T^{-1}$ exists.

(ii) Write V as a direct sum of subvector spaces.

(iii) Write V as a pseudo direct sum of set fuzzy linguistic vector subspaces.

(iv) Define a set fuzzy linear functional on V.

(v) Is it possible to write V as V = W+$W^\perp$? (W ≠ V). Justify your answer.

37. For any set fuzzy linguistic vector space V over the fuzzy linguistic set L find L (V, L) = {all set fuzzy linguistic linear functionals}.

(i) What is the algebraic linguistic structure enjoyed by L (V,L)?

(ii) If V is finite what will be the dimension of L (V, L) (no. of elements in V is n and |L| = m).

(iii) Is L(V, L) a set fuzzy linguistic vector space over L?

38. Study the algebraic structure enjoyed by the collection of all set fuzzy linguistic linear operator $Hom_L(V,V)$ of a set fuzzy linguistic vector space V with $|V| < \infty$ over a finite fuzzy linguistic set L.

    (i) If the vector space is replaced by a linear algebra what is the additional property enjoyed by $Hom_L(V, V)$?

39. Obtain some interesting features enjoyed by fuzzy linguistic cognitive models.

40. Give an example of a fuzzy linguistic overlapping topological space.

41. Using a fuzzy linguistic overlapping topological space study the Klin's temperature problem of cement industries.

42. What are the special features enjoyed by the fuzzy linguistic relation equations model?

43. Describe the fuzzy linguistic relational model with a social problem.

44. Describe the fuzzy linguistic feed forward neural network model with 3 layers.

45. Is it possible to transfer all fuzzy models into fuzzy linguistic models where L the fuzzy linguistic space is a chain lattice connected fuzzy linguistic topological discrete space?

46. Give any nice applications of linguistic fuzzy relational maps model.

47. Let

$$M = \left\{ \begin{bmatrix} a_1 \\ a_2 \\ a_3 \\ a_4 \\ a_5 \end{bmatrix}, \sum_{i=0}^{20} a_i x^i, (a_1, a_2, a_3, a_4) \sum_{i=1}^{5} d_i x^i \right.$$

$$d_i = \begin{bmatrix} a_1 & a_2 & \dots & a_{10} \\ a_{11} & a_{12} & \dots & a_{20} \\ a_{21} & a_{22} & \dots & a_{30} \end{bmatrix}; a_i \in L, 0 \leq i \leq 30 \right\}$$

be a set fuzzy linguistic vector space defined over the fuzzy linguistic set L = {0, low, very low, medium, just low, high, very high, just high, lowest, highest}.

(i) Is M finite or infinite? (That is does M contain infinite number of elements or only finite number of elements).

(ii) Write M as a direct sum of fuzzy linguistic set vector subspaces.

(iii) Find a fuzzy linguistic linear operator on M.

(iv) Find a fuzzy linguistic linear functional on M.

(v) What is the algebraic structure of L (M, L)?

(vi) Find the algebraic structure of $Hom_L(M,M)$.

48. Obtain some interesting features enjoyed by set fuzzy linguistic linear algebra with matrix coefficient polynomials.

49. Let

$$P = \left\{ \sum_{i=0}^{20} a_i x^i \,\middle|\, a_i = \begin{bmatrix} d_1 & d_2 & d_3 \\ d_4 & d_5 & d_6 \\ \vdots & \vdots & \vdots \\ d_{40} & d_{41} & d_{42} \end{bmatrix} ; d_j \in L; 1 \leq j \leq 4 \right\}$$

be a semigroup fuzzy linguistic linear algebra over L.

(i) Find semigroup fuzzy linguistic linear subalgebras.

(ii) Write P as a direct sum of semigroup fuzzy linguistic linear subalgebras.

(iii) Does every semigroup fuzzy linguistic linear subalgebra W of V have $W^\perp$ such that

$$W \cap W^\perp = \begin{bmatrix} 0 & 0 & 0 \\ 0 & 0 & 0 \\ \vdots & \vdots & \vdots \\ 0 & 0 & 0 \end{bmatrix} ?$$

50. Let

$$M = \left\{ \begin{bmatrix} a_1 & a_2 & \cdots & a_{20} \\ a_{21} & a_{22} & \cdots & a_{40} \\ \vdots & \vdots & & \vdots \\ a_{81} & a_{82} & \cdots & a_{100} \end{bmatrix} \,\middle|\, d_i = \begin{bmatrix} a_1 & a_2 & \cdots & a_{10} \\ a_{11} & a_{12} & \cdots & a_{20} \\ a_{21} & a_{22} & \cdots & a_{30} \end{bmatrix} ; \right.$$

$a_i \in L$, $1 \leq i \leq 100$} be a semigroup fuzzy linguistic linear algebra over the fuzzy linguistic space L.

(i) Write M as a direct sum of semigroup fuzzy linguistic linear subalgebras.

(ii) Write M as a pseudo direct sum of semigroup fuzzy linguistic linear subalgebras.

(iii) Write M as $W + W^\perp$. (W a proper semigroup fuzzy linguistic linear subalgebra of M over L).

(iv) Does there exist a semigroup fuzzy linguistic linear operator on M?

51. Obtain some interesting properties about

$$S = \left\{ \sum_{i=0}^{\infty} a_i x^i \;\middle|\; a_i \in L \right\}.$$

(i) Is it a lattice under '$\cup$' and '$\cap$' (i.e., max and min)?

(ii) Can S have sublattice?

(iii) Is S a monoid under 'min'?

52. Let

$$M = \left\{ \sum_{i=0}^{2} a_i x^i \;\middle|\; a_i = (d_1, d_2, d_3, d_4, d_5);\, d_j \in L;\, 1 \leq j \leq 5 \right\}$$

be a semigroup under 'min'.

(i) Can M have ideals?

(ii) Can M have zero divisors?

(iii) Can M have subsemigroups which are not ideals?

53. Study questions (i) and (iii) for M in problem (52) in which min operation is replaced by 'max' operation.

54. Let

$$P = \left\{ \sum_{i=0}^{7} a_i x^i \ \middle| \ a_i = \begin{bmatrix} d_1 & d_2 & d_3 \\ d_4 & d_5 & d_6 \\ d_7 & d_8 & d_9 \\ d_{10} & d_{11} & d_{12} \end{bmatrix} ; d_j \in L; 0 \leq j \leq 12 \right\}$$

be a fuzzy linguistic semigroup under max operation.

(i) Can P have ideals?

(ii) Can P have subsemigroups?

(iii) Can P have zero divisors?

(iv) If $|L| = 12$ say what is the order of P?

(v) Can we write P as a direct sum of subsemigroups?

55. Let

$$T = \left\{ \sum_{i=0}^{20} a_i x^i \ \middle| \ a_i = \begin{bmatrix} d_1 & d_2 & d_3 \\ d_4 & d_5 & d_6 \\ \vdots & \vdots & \vdots \\ d_{28} & d_{29} & d_{30} \end{bmatrix} ; d_j \in L; 1 \leq j \leq 30 \right\}$$

be a fuzzy linguistic semigroup under 'min'.

(i) Find ideals of T.

(ii) Write T as $S + S^\perp$ where S is a subsemigroup under min.

(iii) Write T as a pseudo direct sum of subsemigroups.

(iv) Does T have subsemigroups which are not ideals?

56. Let

$$P = \left\{ \sum_{i=0}^{2} a_i x^i, \sum_{i=0}^{5} d_i x^i \;\middle|\; a_i = (p_1, p_2, \ldots, p_{15}), \right.$$

$$d_j = \begin{pmatrix} q_1 & q_2 & \ldots & q_{10} \\ q_{11} & q_{12} & \ldots & q_{20} \end{pmatrix} \text{ where } p_t, q_k \in L, \ 1 \le t \le 15 \text{ and }$$

$1 \le k \le 20\}$ be a set fuzzy linguistic vector space over the set L.

(i) Write P as a direct sum.

(ii) Write P as a pseudo direct sum.

(iii) Find $\text{Hom}_L(P, P)$.

(iv) Find $L(P, L)$.

57. Let

$$N = \left\{ \sum_{i=0}^{7} a_i x^i \;\middle|\; a_i = \begin{pmatrix} d_1 & d_2 & d_3 & d_4 \\ d_5 & \ldots & \ldots & d_8 \\ d_9 & \ldots & \ldots & d_{12} \\ d_{13} & \ldots & \ldots & d_{16} \\ d_{17} & \ldots & \ldots & d_{20} \end{pmatrix} \text{ with } d_j \in L, \right.$$

$1 \le j \le 20\}$ be a fuzzy linguistic semigroup linear algebra over the fuzzy linguistic semigroup N 'min'.

(i) Write N as a direct sum.

(ii) Write N as pseudo direct sum.

(iii) Let $M \subseteq N$ be a semigroup fuzzy linguistic linear subalgebra; does there always exist a $M^\perp$ such that $N = M + M^\perp$.

(iv) If $M = \left\{ \sum_{i=0}^{4} a_i x^i \mid a_i = \begin{pmatrix} d_1 & d_2 & d_3 & d_4 \\ d_5 & \ldots & \ldots & d_8 \\ d_9 & \ldots & \ldots & d_{12} \\ d_{13} & \ldots & \ldots & d_{16} \\ d_{17} & \ldots & \ldots & d_{20} \end{pmatrix} \; d_j \in L, \right.$

$1 \leq j \leq 20 \} \subseteq N$. Does $M^\perp$ exist?

(v) If min operation on N is replaced by max operation study question (i) to (iv).

58. Let

$$V = \left\{ \sum_{i=0}^{3} a_i x^i \mid a_i = \begin{pmatrix} d_1 \\ d_2 \\ \vdots \\ d_{20} \end{pmatrix} \text{ with } d_j \in L, 1 \leq j \leq 20 \right\}$$

be a set fuzzy linguistic linear algebra over the fuzzy linguistic set L.

(i) Write V as a direct sum.

(ii) For $M = \left\{ \sum_{i=0}^{3} a_i x^i \mid a_i = \begin{bmatrix} d_1 \\ d_2 \\ d_3 \\ d_4 \\ 0 \\ \vdots \\ 0 \end{bmatrix} \text{ with } d_j \in L, 1 \leq j \leq 4 \right\}$

$\subseteq V$ find $M^\perp$ and write $V = M + M^\perp$.

$$\text{If } P = \left\{ \sum_{i=0}^{3} a_i x^i \;\middle|\; a_i = \begin{bmatrix} 0 \\ 0 \\ 0 \\ 0 \\ 0 \\ d_1 \\ d_2 \\ d_3 \\ 0 \\ 0 \\ \vdots \\ 0 \end{bmatrix} \text{ with } d_j \in L,\ 1 \le j \le 3 \right\}$$

find $P \subseteq V$ so that $M + P = V$.
Can P be a total complement of M?

(iii) Can V have subset fuzzy linguistic linear subalgebras?

59. Does there exists a set fuzzy linguistic linear algebra V which has no subset fuzzy linear subalgebra for any L, L a fuzzy linguistic set over which V is defined.

60. Define a fuzzy linguistic set L so that L is a semigroup under min. Let S be a semigroup fuzzy linguistic linear algebra over L. Characterize those subsemigroup fuzzy linguistic linear subalgebras of V defined over L.

61. Let L be a fuzzy linguistic set with $|L| = 8$.

$$V = \left\{ \sum_{i=0}^{3} a_i x^i, (a_1, a_2, a_3) \;\middle|\; a_i \in L,\ 0 \le i \le 3 \right\} \text{ be a set fuzzy}$$

linguistic vector space over L.

(i) Find $|V|$.

(ii) Write V as a direct sum.

(iii) How many set fuzzy linguistic operators are in $\text{Hom}_L (V,V)$?

(iv) What is the cardinality of $L(V,L)$?

62. Let $M = \left\{ \sum_{i=0}^{4} a_i x^i \mid a_i = (d_1, d_2) \text{ where } d_1, d_2 \in L \right.$ where L is a fuzzy linguistic set with four fuzzy linguistic terms including 0$\}$ be a set fuzzy linguistic linear algebra over L

(i) Find $|M|$.

(ii) Find $|L(M, L)|$.

(iii) Find $|\text{Hom}_L (M, M)|$.

(iv) Write M as a direct sum.

(v) Write M as sum of an orthogonal space and its complement.

63. Find some applications of lattice chain connected fuzzy linguistic topological spaces L with $|L| = 9$.

64. Find some applications of set fuzzy linguistic linear algebra V over a fuzzy linguistic set L.

65. Let L = {good, bad, very good, fair, best, very bad, 0, very fair, just good, just bad, better} be a fuzzy linguistic set.

(i) Is L a connected topological space?

(ii) Is L a over lapping topological space?

(iii) What is the metric that can be defined on L?

66. Give a real world example of a feed forward fuzzy linguistic neural network with one layer.

67. Use the feed forward fuzzy linguistic neural network with two layers in a problem and describe the problem.

68. Enumerate the advantages of using fuzzy linguistic models for social problems.

69. Can fuzzy linguistic models be used in problems related to industries? Justify / substantiate your claim.

70. Can fuzzy linguistic models be used as mathematical models? Justify.

71. Obtain some applications of fuzzy linguistic relational maps model.

72. Suppose L is a fuzzy linguistic set.

    (i) In how many ways can a fuzzy linguistic metrics on L be defined?

    (ii) Does the number of fuzzy linguistic metrics on L depend only on L?

73. Give an example of L; L a fuzzy linguistic set on which atleast 3 fuzzy linguistic metrics can be defined.

74. Is it always possible to translate a fuzzy membership function to a fuzzy linguistic membership function?

75. Can we translate a fuzzy linguistic membership function into a fuzzy membership function?

76. Describe by a real world problem the applications of FLCM.

77. Give some possible application of set fuzzy linguistic vector spaces.

78. Distinguish between the set fuzzy linguistic vector spaces and set fuzzy linguistic linear algebras defined over a fuzzy linguistic set L.

79. Distinguish between the semigroup fuzzy linguistic vector spaces and set fuzzy linguistic vector spaces defined over a fuzzy linguistic set L.

80. What is the advantage of using semigroup fuzzy linguistic linear algebra in the place of semigroup fuzzy linguistic vector space?

81. Let
$$V = \left\{ \sum_{i=0}^{10} a_i x^i, \sum_{i=0}^{4} d_i x^i \,\middle|\, d_j = \begin{bmatrix} g_1 & g_2 & g_3 & \cdots & g_{10} \\ g_{11} & g_{12} & g_{13} & \cdots & g_{20} \\ g_{21} & g_{22} & g_{23} & \cdots & g_{30} \\ g_{31} & g_{32} & g_{33} & \cdots & g_{40} \\ g_{41} & g_{42} & g_{43} & \cdots & g_{50} \end{bmatrix} \right\}$$
with $g_k, a_i \in L$, $0 \le j \le 4$; $1 \le k \le 50$ and $0 \le i \le 10$} be a set fuzzy linguistic vector space over the fuzzy linguistic set L.

   (i) Write V as direct sum.

   (ii) Is $V = \left\{ \sum_{i=0}^{10} a_i x^i \,\middle|\, a_i \in L \right\} \subseteq L$ a pseudo set fuzzy linguistic linear algebra under 'min' operation?

   (iii) Write V as a pseudo direct sum.

   (iv) Find $\text{Hom}_L (V, V)$.

   (v) What is the algebraic structure enjoyed by L(V, L)?

82. Let
$$V = \left\{ \begin{bmatrix} a_1 & a_2 & \cdots & a_8 \\ a_9 & a_{10} & \cdots & a_{16} \\ a_{17} & a_{18} & \cdots & a_{24} \\ a_{25} & a_{26} & \cdots & a_{32} \\ a_{33} & a_{34} & \cdots & a_{40} \\ a_{41} & a_{42} & \cdots & a_{48} \end{bmatrix} \,\middle|\, a_i \in L, \right.$$

L a fuzzy linguistic semigroup under min; $|L| = 12$, $1 \leq i \leq 48$} be a semigroup fuzzy linguistic linear algebra over L.

(i) What is $|V|$?

(ii) Write $V = W + W^\perp$.

(iii) Find zero divisors in V.

(iv) If 'min' is replaced by max can V have zero divisors?

(v) If 'min' is replaced by 'max' in V can V be written as $V = W + W^\perp$?

83. What is the advantage of using 'min' instead of 'max' in semigroup fuzzy linguistic linear algebras V defined over the fuzzy linguistic semigroup L with min?

84. Describe some special properties enjoyed by set fuzzy linguistic spaces.

85. Let $M = \left\{ \sum_{i=0}^{6} a_i x^i \;\middle|\; a_i = \begin{bmatrix} d_1 \\ d_2 \\ d_3 \\ \vdots \\ d_{15} \end{bmatrix} \right.$ with $d_j \in L$, $1 \leq j \leq 15\}$ be a set fuzzy linguistic linear algebra under min operation over the fuzzy linguistic set L.

(i) If $|L| = 7$ find $|M|$.

(ii) Write $W + W^\perp = M$.

(iii) Write M as a direct sum.

(iv) Write M as a pseudo direct sum.

86. Let

$$M = \left\{ \sum_{i=0}^{7} a_i x^i \;\middle|\; a_i = \begin{bmatrix} d_1 & d_2 & d_3 \\ d_4 & d_5 & d_6 \\ d_7 & d_8 & d_9 \end{bmatrix} \text{ with } d_j \in L,\ 1 \le j \le 9 \right\}$$

be a semigroup fuzzy linguistic linear algebra over L. $|L| = 7$. How many W's exists in M such that $W + W^\perp = M$.

87. Write down the algebraic structure of $\text{Hom}_L(M, M)$ in problem (86).

88. Let

$$M = \left\{ \sum_{i=0}^{8} a_i x^i \;\middle|\; a_i = \begin{bmatrix} d_1 & d_2 & d_3 \\ d_4 & d_5 & d_6 \\ d_7 & d_8 & d_9 \end{bmatrix} \text{ with } d_j \in L,\ 1 \le j \le 9 \right\}$$

and

$$N = \left\{ \sum_{i=0}^{17} a_i x^i \;\middle|\; a_i = (d_1, d_2, \ldots, d_9),\ d_j \in L,\ 1 \le j \le 9 \right\}$$ be a

semigroup fuzzy linguistic linear algebras over L.

   (i) Find $T : M \to N$ so that T is a semigroup fuzzy linguistic linear transformation from M to N.

   (ii) Find the algebraic structure enjoyed by $\text{Hom}_L(M, N)$.

   (iii) What is the relation between L(M, L) and L(N, L)?

89. If M and N be any two set fuzzy linguistic linear algebras defined over a fuzzy linguistic set L. What can be said about L(M, L) and L(N, L) if $|M| = |N|$?

90. Obtain the algebraic structure enjoyed by $\text{Hom}_L(V, W)$; $V \ne W$, V and W set fuzzy linguistic vector spaces over a fuzzy linguistic set L.

91. Give an example of a lattice connected fuzzy linguistic topological space where lattice is not a chain lattice.

92. Does there exists a fuzzy linguistic topological space which is graph connected and the graph is a complete graph? Justify your claim.

93. What is the difference between a lattice connected fuzzy linguistic topological space and graph connected fuzzy linguistic topological space?

94. Derive a difference between a lattice connected fuzzy linguistic topological space and a graph disconnected fuzzy linguistic topological space.

95. Give an example of fuzzy linguistic over lapping topological space with 11 elements.

96. Let $M = \left\{ \sum_{i=0}^{9} a_i x^i, \sum_{i=1}^{3} n_i x^i, \sum_{i=2}^{7} m_i x^i \middle| a_i = \begin{bmatrix} d_1 & d_2 \\ d_3 & d_4 \\ \vdots & \vdots \\ d_{19} & d_{20} \end{bmatrix}, n_j = \begin{bmatrix} t_1 & t_2 & \cdots & t_{20} \\ t_{21} & t_{22} & \cdots & t_{40} \\ t_{41} & t_{42} & \cdots & t_{60} \end{bmatrix} \text{ and } m_t = \begin{bmatrix} p_1 & p_2 & \cdots & p_7 \\ p_8 & p_9 & \cdots & p_{14} \\ p_{15} & p_{16} & \cdots & p_{21} \\ \vdots & \vdots & & \vdots \\ p_{43} & p_{44} & \cdots & p_{49} \end{bmatrix} \text{ where }$

$t_s, d_t, p_r \in L; 1 \leq s \leq 60, 1 \leq r \leq 49$ and $1 \leq t \leq 20$} be a set fuzzy linguistic vector space over the set L.

(i) Find $\text{Hom}_L (M, M)$.

(ii) Write M as a direct sum.

(iii) Write M as a pseudo direct sum.

97. If M is a semigroup fuzzy linguistic vector space over the fuzzy linguistic semigroup L under 'min' operation (M as

in problem 96). Find the special features enjoyed by M as a semigroup fuzzy linguistic vector space.

98. Let $P = \left\{ \sum_{i=0}^{15} a_i x^i \mid a_i = \begin{pmatrix} d_1 & d_2 & d_3 & d_4 \\ d_5 & \ldots & \ldots & d_8 \\ d_9 & \ldots & \ldots & d_{12} \\ d_{13} & \ldots & \ldots & d_{16} \end{pmatrix} d_j \in L, 1 \leq j \leq 16, |L| = 5 \right\}$ be a semigroup fuzzy linguistic linear algebra with 'min' operation.

   (i) Find |P|.

   (ii) Write $P = W + W^\perp$.

   (iii) Write P as pseudo direct sum.

99. For problem (98) by appropriately selecting the fuzzy linguistic terms find for P subsemigroup fuzzy linguistic linear subalgebras.

100. Does there exist a semigroup fuzzy linguistic linear algebra which has no subsemigroup fuzzy linguistic linear subalgebras?

101. Does there exist a set fuzzy linguistic linear algebra V over a fuzzy linguistic set L such that V has no subset fuzzy linguistic linear subalgebras.

102. Suppose L in problem 101 is a fuzzy linguistic chain lattice can we say every singleton set with zero, that is $S = \{a_i, 0\}$ ($a_i \neq 0$, $a_i \in L$) is such that V can have subset fuzzy linguistic linear subalgebra of V over L.

# FURTHER READING

# ABOUT THE AUTHORS

**Dr.W.B.Vasantha Kandasamy** is an Associate Professor in the Department of Mathematics, Indian Institute of Technology Madras, Chennai. In the past decade she has guided 13 Ph.D. scholars in the different fields of non-associative algebras, algebraic coding theory, transportation theory, fuzzy groups, and applications of fuzzy theory of the problems faced in chemical industries and cement industries. She has to her credit 646 research papers. She has guided over 68 M.Sc. and M.Tech. projects. She has worked in collaboration projects with the Indian Space Research Organization and with the Tamil Nadu State AIDS Control Society. She is presently working on a research project funded by the Board of Research in Nuclear Sciences, Government of India. This is her $64^{th}$ book.

On India's 60th Independence Day, Dr.Vasantha was conferred the Kalpana Chawla Award for Courage and Daring Enterprise by the State Government of Tamil Nadu in recognition of her sustained fight for social justice in the Indian Institute of Technology (IIT) Madras and for her contribution to mathematics. The award, instituted in the memory of Indian-American astronaut Kalpana Chawla who died aboard Space Shuttle Columbia, carried a cash prize of five lakh rupees (the highest prize-money for any Indian award) and a gold medal.
She can be contacted at vasanthakandasamy@gmail.com
Web Site: http://mat.iitm.ac.in/home/wbv/public_html/
or http://www.vasantha.in

**Dr. Florentin Smarandache** is a Professor of Mathematics at the University of New Mexico in USA. He published over 75 books and 200 articles and notes in mathematics, physics, philosophy, psychology, rebus, literature. In mathematics his research is in number theory, non-Euclidean geometry, synthetic geometry, algebraic structures, statistics, neutrosophic logic and set (generalizations of fuzzy logic and set respectively), neutrosophic probability (generalization of classical and imprecise probability). Also, small contributions to nuclear and particle physics, information fusion, neutrosophy (a generalization of dialectics), law of sensations and stimuli, etc. He got the 2010 Telesio-Galilei Academy of Science Gold Medal, Adjunct Professor (equivalent to Doctor Honoris Causa) of Beijing Jiaotong University in 2011, and 2011 Romanian Academy Award for Technical Science (the highest in the country). Dr. W. B. Vasantha Kandasamy and Dr. Florentin Smarandache got the 2011 New Mexico Book Award for Algebraic Structures. He can be contacted at smarand@unm.edu

**Rev. Fr. K. Amal, SJ** is presently pursuing his Ph.D. in Mathematics in the Madurai Kamaraj University, Tamil Nadu, India. Having completed his bachelor's degree in mathematics at the St. Joseph's College, Trichy, he entered the Society of Jesus in July 1977. In 1985, he completed his master's degree in mathematics in Loyola College, Chennai and his master's degree in sociology in 1992. Ordained in 1988, he started his service as a parish priest in Mallihapuram, a Dalit parish, in the Chengelpet district. He later served the All India Catholic University Federation (AICUF ), a student movement, for 15 long years till 2003; in the latter seven-and-a-half years he occupied the position of National Director, AICUF. During his tenure, the movement witnessed newer areas of operation and scaled greater heights in its expansion and animation works, particularly with Dalit and Tribal students. Transferred to the Indian Social Institute, Bangalore, in 2004, as its Superior, he was invited to be the Director of Works of the Human Rights Desk and the Training Unit. Besides this, he has been an Advisor to the **Asia-Pacific** Chapter of the International Movement of Catholic Students (IMCS) since June 2003. In this capacity he had traveled far and wide in the Asia Pacific region and elsewhere too. He has organized, participated, animated and served as a resource person in many workshops, seminars and meetings. He can be contacted at [kamalsj@gmail.com](kamalsj@gmail.com)